\documentclass{amsart}
\usepackage{amsmath,amsthm,amssymb,amscd,url,enumerate}
\usepackage[utf8]{inputenc}
\usepackage{hyperref}
\usepackage{pdflscape}
\usepackage{caption}
\usepackage[noadjust]{cite}
\usepackage{xypic}
\usepackage{tikz}
\usetikzlibrary{decorations.markings,arrows}
\usepackage{tikz-cd}
\usepackage{amsfonts}
\usepackage{amsmath}
\usepackage{graphicx}
\usepackage{xcolor}
\usepackage{mathrsfs}
\usepackage{subcaption}
\usepackage[shortlabels,inline]{enumitem}
\usepackage{comment}
\usepackage{mathtools}
\usepackage{fullpage}

\theoremstyle{plain} 
\newtheorem{theorem}{Theorem}[section]

\newtheorem{proposition}[theorem]{Proposition}
\newtheorem{lemma}[theorem]{Lemma}
\newtheorem{corollary}[theorem]{Corollary}

\newtheorem{mainthm}{Main Theorem}

\newtheorem{comp}[mainthm]{Comparison}

\theoremstyle{definition}

\newtheorem{example}[theorem]{Example}
\newtheorem{definition}[theorem]{Definition}
\newtheorem{remark}[theorem]{Remark}

\DeclareMathOperator{\End}{\operatorname{End}}
\DeclareMathOperator{\Aut}{\operatorname{Aut}}
\DeclareMathOperator{\Cl}{\operatorname{Cl}}
\DeclareMathOperator{\Hom}{\operatorname{Hom}}

\DeclareMathOperator{\val}{val}
\DeclareMathOperator{\ord}{ord}

\usepackage[marginparwidth=7mm, marginparsep=1mm,]{geometry} 

\newcommand{\ZZ}{{\mathbb{Z}}}
\newcommand{\QQ}{{\mathbb{Q}}}
\newcommand{\Ocal}{{\mathcal{O}}}
\newcommand{\Ical}{{\mathcal{I}}}
\newcommand{\OO}{{\mathcal{O}}}
\newcommand{\Fq}{{\mathbb{F}_q}}

\newcommand{\FF}{{\mathbb{F}}}
\newcommand{\ff}{{\mathfrak{f}}}
\newcommand{\qq}{{\mathfrak{q}}}
\newcommand{\PP}{{\mathfrak{P}}}
\newcommand{\pp}{{\mathfrak{p}}}
\newcommand{\frl}{\mathfrak{l}}
\newcommand{\frL}{\mathfrak{L}}
\newcommand{\mm}{\mathfrak{m}}
\newcommand{\Fpbar}{{\overline{\mathbb{F}}_p}}

\newcommand{\dmin}{d_{\min}}
\newcommand{\Gasc}{G^\text{asc}}
\newcommand{\Gdesc}{G^\text{desc}}
\newcommand{\coloneqq}{=}

\title{Isogeny graphs of abelian varieties and singular ideals in orders}
\author{Sarah Arpin}
\address{{Virginia Polytechnic Institute and State University, Blacksburg, VA 24061}}
\email{{sarpin@vt.edu}}

\author{Stefano Marseglia}
\address{Mathematical Institute, Utrecht University, P.O. Box 80010, 3508 TA, Utrecht, The Netherlands}
\address{Laboratoire GAATI, Université de la Polynésie française, BP 6570 -- 98702 Faaa, Polynésie française}
\address{(current) Laboratoire Jean Alexandre Dieudonné, Université Côte Azur, 06108 Nice Cedex 2, France}
\email{stefano.marseglia@univ-cotedazur.fr}

\author{Caleb Springer}
\address{Center for Communications Research, Princeton, NJ 08540}
\email{c.springer.math@gmail.com}


\date{\today}

\begin{document}
\begin{abstract}
Famously, Kohel proved that isogeny graphs of ordinary elliptic curves are beautifully structured objects, now called volcanos.
We prove graph structural theorems for abelian varieties of any dimension with commutative endomorphism ring and containing a fixed locally Bass order, leveraging an ideal-theoretic perspective on isogeny graphs. 
This generalizes previous results, which relied on restrictive additional assumptions, such as maximal real multiplication, ordinary, and absolutely simple (Brooks, Jetchev, Wesolowski 2017). 
In~particular, our work also applies to non-simple and non-ordinary isogeny classes.
To obtain our results, we first prove a structure theorem for the lattice of inclusion of the overorders of a locally Bass order in an \'etale algebra which is of independent interest. 
This analysis builds on a careful study of local singularities of the orders.
We include several examples of volcanoes and isogeny graphs exhibiting unexpected properties ultimately due to our more general setting.
\end{abstract}

\maketitle

\section{Introduction}\label{sec:intro}
Isogeny graphs of elliptic curves are structures which arrange the elements of an isogeny class in a way that is both intrinsically beautiful and helpful for solving algorithmic problems,  see \cite{Sutherland_2013} for an expository treatment.
For both flavors of elliptic curves over fields of characteristic $p$, namely the ordinary and supersingular cases, isogeny graphs have been widely studied and are well-understood \cite{Waterhouse,Pizer,Kohel94}.
In particular, the connected components of $\ell$-isogeny graphs of ordinary elliptic curves defined over a finite field of characteristic $p$ are \textit{isogeny volcanoes} \cite{Kohel94}, a term first coined by Fouquet and Morain \cite{FouquetMorain02}. 
The same is true for supersingular elliptic curves defined over a prime field $\mathbb{F}_p$.

In certain special cases, the structure theorems for the $\ell$-isogeny graphs of elliptic curves have been generalized to higher dimensional abelian varieties.
In particular, absolutely simple ordinary abelian varieties with so-called locally maximal real multiplication are considered in \cite{BrooksJetchevWesolowski17,Martindale18,IonicaThome20}.
See Section~\ref{subsec:related_work} for a recollection of the setup of \cite{BrooksJetchevWesolowski17}, whose main results imply the ones of \cite{Martindale18} and \cite{IonicaThome20}.
The ultimate goal of our paper is to provide a structure theorem for isogeny graphs of abelian varieties in wider generality. 
This includes the removal of the assumption of being ordinary, analogously to the elliptic curve case.
Along the way, we also investigate the behavior of singular ideals in lattices of orders in \'etale algebras, which is of independent interest.

The starting point of our investigation is the observation that the hypotheses of previous results imply that the endomorphism rings of the abelian varieties are (locally) \emph{Gorenstein orders}, see Definition~\ref{def:Gor_Bass}.
After investigating the Gorenstein property directly and providing a structure theorem for the minimal overorders of a (locally) Gorenstein order, we exploit an ideal-theoretic viewpoint on isogeny graphs to produce the most general structure theorem for isogeny graphs of non-polarized abelian varieties with commutative endomorphism algebras currently known.

\subsection{Main contributions}

Recall that an order $R$ in a finite product of number field is called a \emph{Bass order} at some maximal ideal $\frl$ if all its $\frl$-overorders are Gorenstein at the maximal ideals above $\frl$, see Definitions~\ref{def:frl_oo} and~\ref{def:Gor_Bass}.
Our first main result, stated below as Main Theroem~\ref{mainthm:Bass_cond}, is a classification of the overorders of an order which is Bass at $\frl$.
It refines \cite[Cor.~5.25]{HofmannSircana20}.
Additionally, this result allows us to generalize, the main result of \cite{ChoHongLee24_arXiv} where the authors provide an explicit classification of the overorders of $R$ when $R$ is Bass at every maximal ideal.
We do not require the order to be an integral domain, and, in fact, we even relax the assumption of being Bass (at every maximal ideal) to a weaker one.
This is Corollary~\ref{cor:number_of_ladders}.
We stress that our proof is considerably shorter.
Before stating Main Theorem~\ref{mainthm:Bass_cond}, we need to introduce the concept of \emph{$\frl$-multiplicator ladder of an order $R$}, see Definition~\ref{def:multiplicatorladder}.
We say that a chain of inclusions
\[ R=R_d \subsetneq R_{d-1} \subsetneq \ldots \subsetneq R_0 \]
of overorders of $R$ is the $\frl$-multiplicator ladder of $R$ if the following conditions hold:
\begin{itemize}
    \item $R=R_d, R_{d-1}, \ldots, R_0$ are all the $\frl$-overorders of $R$;
    \item for each $i=1,\ldots,d$, the order $R_i$ has a unique maximal ideal $\frL_i$ above $\frl$ and $\frL_i$ is singular;
    \item for each $i=1,\ldots,d$, the order $R_{i-1}$ equals $(\frL_i:\frL_i)$ and is the unique minimal $\frl$-overorder of $R_i$.
\end{itemize}

\begin{mainthm}\label{mainthm:Bass_cond}
    Let $R$ be an order and let $\frl$ be a maximal ideal of $R$.
    Assume that $R$ is Bass at $\frl$.
    Let $T$ be an overorder of $R$.
    Then there exists a unique overorder $\OO$ of $R$ such that:
    \begin{enumerate}[(i)]
        \item $\OO_\frl = R_\frl$;
        \item $\OO$ has an $\frL$-multiplicator ladder $\OO=\OO_d\subsetneq \OO_{d-1} \subsetneq \ldots \subsetneq \OO_0$ where $\frL=\frl\OO$ is the unique singular maximal ideal of $\OO$ above $\frl$;
        \item there exists a unique index $0\leq i\leq d$ such that $T=\OO_i$.
    \end{enumerate}
    Moreover, for $i=0,\ldots,d$, the conductor $\ff_{\OO_i}=(\OO_i:\OO_K)$ satisfies $\ff_{\OO_i}=\frL^i\ff_{\OO_0} = \frl^i\ff_{\OO_0}$ and $\OO_i = \OO + \frl^i\ff_{\OO_0} = R +\frl^i\ff_{\OO_0}$.
\end{mainthm}
An expanded and slightly more general version of Main Theorem~\ref{mainthm:Bass_cond} is found below as Theorem~\ref{thm:all_ladders}.
The key ingredient of the proof is a careful analysis of the maximal ideals above $\frl$.
The statements in Main Theorem~\ref{mainthm:Bass_cond} are classical for an order $R$ in a quadratic imaginary field.

\bigskip

Main Theorem~\ref{mainthm:Bass_cond} is the backbone of the proof of our second main contribution, namely Main Theorem~\ref{mainthm:B}, which is about isogeny graphs of abelian varieties defined over finite fields.
We introduce here only the notation needed to state the theorem.
Let $A_0$ be an abelian variety defined over a finite field $\Fq$ with commutative $\Fq$-endomorphism algebra $K$.
Let $\pi$ be the element of $K$ representing the Frobenius endomorphism of $A_0$.
For every abelian variety $A$ defined over $\Fq$ which is $\Fq$-isogenous to $A_0$, we identify $\End_{\Fq}(A)$ with an order in $K$ by sending the Frobenius endomorphism of $A$ to $\pi$.

Now, let $R$ be an order in $K$ and $\frl$ be a maximal ideal of $R$ above a rational prime $\ell$, coprime with $q$.
Let $\Ical_R$ be the set of abelian varieties $A$ which are $\Fq$-isogenous to $A_0$ and such that $R \subseteq \End(A)$.
For each fractional $R$-ideal $I\subseteq \End(A)$ we define $A[I]$ as 
the (scheme-theoretic) intersection of $\ker(\alpha)$ where $\alpha$ runs over all the elements of $I\cdot \End(A)$, the extension of $I$ to $\End(A)$.

We say that an $\Fq$-isogeny $\varphi: A \to B$ is an \emph{(ascending or horizontal) $\frl$-isogeny} from $A$ if 
$\ker\varphi \cong R/\frl$ as $R$-modules, and
if $\varphi$ factors as an $\frL$-multiplication followed by an isomorphism, that is, 
$\varphi: A \to A/A[\frL] \cong B$, where $\frL$ is a maximal ideal of $\OO_i=\End(A)$ above $\frl$.
For such an $\frl$-isogeny $\varphi:A \to B$, if $\End(A) \subsetneq \End(B)$, we define also a \emph{(descending) $\frl$-isogeny} $B\to C$ by the commutativity of a certain diagram, see Definition~\ref{def:l_isogenies}.

We define the \emph{$(R,\frl)$-isogeny graph} (see Definition~\ref{def:R_frl_graph}) as the directed graph where:
\begin{itemize}
	\item The vertices are $\Fq$-isomorphism classes of abelian varieties in $\Ical_R$.
    \item For every $\frl$-isogeny from $A$ to $B$ for $A,B\in\Ical_R$ up to pre- and post-composition with automorphisms (see Appendix~\ref{app:auts}), we place a directed edge from the vertex of the isomorphism class of $A$ to the vertex of the isomorphism class of $B$.
\end{itemize}

Assume that $R$ is Bass at $\frl$.
Let $G$ be a connected component of the $(R,\frl)$-isogeny graph, see Lemma~\ref{lem:stronglyconnected}.
We identify $A\in \Ical_R$ with the vertex it represents in the graph.
In Lemma~\ref{le:graph_ladder}.\ref{it:component_height} and Proposition~\ref{prop:same_num_AV_along_ladder}.\ref{prop:same_num_AV_along_ladder:yesEnds}, we show that
there is a unique $\frl$-multiplicator ladder $\OO_d\subsetneq \dots \subsetneq \OO_0$ in the set of overorders of $R$ and an integer $\dmin$ satisfying $d\geq \dmin\geq 0$
such that $\OO_{d_{\min}},\dots, \OO_0$ are precisely the endomorphism rings of the abelian varieties in $G$.
For each $0\leq i\leq d_{\min}$, let $G_i$ be the subgraph of $G$ whose vertices are isomorphism classes of abelian varieties at level $i$, meaning $\End(A) = \OO_i$.
We refer to $G_0$ as the \emph{surface} of $G$.
Set $\delta_\frl=-1$ (resp.~$0$, $1$) if $\frl$ is inert in $\OO_0$ (resp.~ramified or split), see Definition~\ref{def:min_ext_primes}. 
\begin{mainthm}\label{mainthm:B}
    Assume that $R$ is Bass at $\frl$.
    \begin{enumerate}[(a)]
        \item \label{mainthm:B:it:surface} 
        Let $\Cl_\frl(\OO_0)$ be  the subgroup of $\Cl(\OO_0)$ generated by the maximal $\OO_0$-ideals lying over $\frl$.
        The vertices of $G_0$ consists of a single orbit of $\Cl_\frl(\OO_0)$.
        \begin{itemize}
            \item If $\delta_\frl =-1$ then $G_0$ is a totally disconnected graph, that is, the set of edges is empty.
            \item If $\delta_\frl =0$ or $1$ then $G_0$ is isomorphic to the (directed) Cayley graph of $\Cl_\frl(\OO_0)$ with generators the maximal ideals of $\OO_0$ above $\frl$.
            In particular, $G_0$ is connected.
        \end{itemize}
        
        \item \label{mainthm:B:it:edges_up_surface} If $i = 0$ and $A\in G_0$, then there are $\frac{ \#(R/\frl)-\delta_\frl}{[\OO_{0}^\times : \OO_1^\times]}$ vertices in $G_1$ whose (unique) ascending edge has target $A$.
        There are no edges from $G_j$ to $G_0$ for $j > 1$. 
        
        \item \label{mainthm:B:it:edges_up_lower} If $1 \leq i < \dmin$ and $A\in G_i$, then there are $\frac{ \#(R/\frl)}{[\OO_{i}^\times : \OO_{i+1}^\times]}$ vertices in $G_{i+1}$ whose (unique) ascending edge has target $A$. 
        These are the only edges with targets in $G_i$.
                
        \item \label{mainthm:B:it:one_desc_per_asc} The in-degree is equal to the out-degree for each vertex on $G_i$ for $1\leq i \leq \min\{\dmin, d-1\}$.  Moreover, if $\dmin = d$, then the vertices at level $d$ have out-degree $1$ and in-degree $0$.
    \end{enumerate}
\end{mainthm}
This result can be found later in the text as Theorem~\ref{thm:graph_structure}.
Main Theorem~\ref{mainthm:B} can be considered as a generalization of \cite[Theorem~4.3]{BrooksJetchevWesolowski17}.
We make this statement precise in Section~\ref{subsec:related_work} below.
The small discrepancy between Theorem~\ref{mainthm:B}.\ref{mainthm:B:it:surface} and \cite[Theorem 4.3.(ii)]{BrooksJetchevWesolowski17} is discussed in Remark~\ref{rmk:discrepancyBJW}.
Main Theorem~\ref{mainthm:B} can easily be turned into an algorithm to compute $(R,\frl)$-isogeny graphs.
We use this algorithm to compute several examples. 
An implementation is available at \url{https://github.com/stmar89/R_ell_isogeny_graphs/}.

\bigskip

Our third main contribution, given below as Main Theorem~\ref{mainthm:C}, gives a characterization of when a connected component $G$ of an $(R,\frl)$-isogeny graph is an \emph{$r$-volcano} for some positive integer $r$, see Definition~\ref{def:r-volcano}.
Strictly speaking $G$ cannot be a volcano since $G$ is a directed graph, while volcanoes are undirected.
So, we define the \emph{undirected ascending graph} $\Gasc$ (resp.~\emph{undirected descending graph} $\Gdesc$) associated to $G$ as the graph with the same vertices as $G$ and whose edges are ascending (resp.~descending) edges of $G$ without their direction, and whose horizontal edges are determined as follows: draw $k$ undirected edges connecting a pair of vertices $(v_1,v_2)$ whenever the directed graph contains both $k$ horizontal edges of the form $(v_1,v_2)$ and $k$ horizontal edges of the form $(v_2,v_1)$ for $v_1\neq v_2$, and if a vertex has a loop then draw this loop as undirected (it may only be traversed in one direction).
We say that $G$ is an $r$-volcano if and only if $\Gasc=\Gdesc$ and $\Gasc$ is an $r$-volcano, see Definitions~\ref{def:asc_desc_G}  and \ref{def:G_volcano}.

To simplify the exposition, we will give here a characterization for when $G$ is an $r$-volcano only when $\dmin>0$, that is, when $G$ does not consist only of a surface. 
The complete statements can be found in Theorems~\ref{thm:volcano_inertramified} and \ref{thm:volcano_split}.
\begin{mainthm}\label{mainthm:C}
    If $\dmin>0$ then $G$ is an $r$-volcano if and only if the following three conditions hold:
    \begin{itemize}
        \item $\frl\OO_{\dmin-1}$ is principal;
        \item\label{lem:gascvolcano_it1} $r +1 =\frac{\#(R/\frl)-\delta_\frl}{[\OO_0^\times:\OO_1^\times]}+r_0$, where $r_0=\delta_\frl+1$;
        \item\label{lem:gascvolcano_it2} $ r + 1 = \frac{\#(R/\frl)}{[\Ocal_i^\times : \Ocal_{i+1}^\times]} + 1$ for $1 \leq i < \dmin$.
    \end{itemize}
\end{mainthm}
Main Theorem~\ref{mainthm:C} is a more general version of the last part of \cite[Theorem~4.3]{BrooksJetchevWesolowski17}.
Again, more details about this comparison are given in Section~\ref{subsec:related_work} below.

\subsection{Comparison with \texorpdfstring{\cite[Theorem 4.3]{BrooksJetchevWesolowski17}}{[BJW17, Theorem 4.3]}}\label{subsec:related_work}

The authors of \cite{BrooksJetchevWesolowski17} study a notion of isogeny graph that is closely related to the one we introduced before, as we now explain.
To facilitate a comparison, we slightly modify their notation.

Let $A_0$ be an ordinary and absolutely simple abelian variety defined over a finite field $\Fq$.
The geometric ($\overline{\Fq}$) and $\Fq$-endomorphism algebras of $A_0$ coincide, and they are both isomorphic to the number field $K=\QQ(\pi)$.
Let $K^+$ be the totally real subfield of $K$ and let $\frl^+$ be a maximal ideal of $\OO_{K^+}$ above a rational prime $\ell$ coprime with $q$. 
Combining \cite[Theorem 3.3]{Howe95} and \cite[proof of Proposition 6.1]{Mar22_extensions}, we see that
if $A$ and $B$ are two abelian varieties $\Fpbar$-isogenous to $A_0$ which can be defined over the extension $\FF_{q^n}$ then
\begin{equation}\label{eq:same_homs}
   \Hom_{\FF_{q^n}}(A,B) = \Hom_{\Fpbar}(A,B). 
\end{equation}

The main theorem \cite[Theorem 4.3]{BrooksJetchevWesolowski17} gives the structure of the so-called \emph{$\frl^+$-isogeny graphs}, originally defined in \cite[Section 4.1]{BrooksJetchevWesolowski17}:
\begin{itemize}
    \item The vertices are isomorphism classes abelian varieties $A$ which are $\Fq$-isogenous to $A_0$, satisfying $(\End(A)\cap K^+)_{\ell} = (\OO_{K^+})_{\ell}$. 
   
    \item The edges from the vertex with representative $A$ are $\frl^+$-isogenies with domain $A$, that is, isogenies defined over $\Fpbar$ whose kernel is a proper, $(\End(A) \cap K^+)$-stable subgroup of $A[\mathfrak{l^+}]$, where $A[\mathfrak{l^+}] = \{x\in A(\overline{\mathbb{F}}_q):\alpha x = 0\text{ for all }\alpha \in\mathfrak{l^+}\OO_K\cap \End(A)\}$.
\end{itemize}
Set
\begin{equation}\label{eq:Randfrl}
    R \coloneqq \ZZ + N_1\OO_{K^+} + \ell N_2\OO_K
    \qquad\text{and}\qquad \frl\coloneqq \frl^+\OO_K\cap R,
\end{equation}
where $[\OO_{K^+} : \ZZ[\pi + q/ \pi]]= N_1\cdot \ell^a$ for $\gcd(N_1, \ell) = 1$ and $[\OO_K : \ZZ[\pi, q/\pi]] = N_2$.
We stress that this choice of $R$ and $\frl$ is not the only one that allows us to compare with \cite[Theorem 4.3]{BrooksJetchevWesolowski17}. 
The ideal $\frl$ is a maximal ideal of $R$, with residue field isomorphic to $\OO_{K^+}/\frl^+$ and $R$ is Bass at $\frl$, see Lemma~\ref{lem:resfieldBJW}.
In order to compare \cite[Theorem 4.3]{BrooksJetchevWesolowski17} and Main Theorem~\ref{mainthm:B}, we truncate the $\frl^+$-isogeny graph, which is infinite, by taking the subgraph of abelian varieties defined over a fixed finite field $\Fq$. 
In view of Equation~\eqref{eq:same_homs}, the vertices of the subgraph coincide with the $\Fq$-isomorphism classes and the edges represent isogenies that are defined over $\Fq$.

Th following statement follows from Lemma~\ref{lem:compvert} and Lemma~\ref{lem:compedges}.
\begin{comp}\label{compBJW}
The truncated $\frl^+$-isogeny graph coincides with the $(R,\frl)$-isogeny graph.
\end{comp}
The two assumptions we use in our definition of $(R,\frl)$-isogeny graphs are implied, but not equivalent, to the assumptions in \cite{BrooksJetchevWesolowski17}.
    Firstly, the endomorphism algebra $K$ of an absolutely simple ordinary abelian variety is a number field, as can be seen by combining, for example, \cite[Thm.~3.3]{Howe95} and \cite[Thm.~2.(c)]{Tate66}.
    Secondly, orders in $K$ with maximal totally real multiplication locally at $\ell$ are Gorenstein at the maximal ideals above $\ell$ by \cite[Lemma 4.4]{BrooksJetchevWesolowski17}.

Our construction being more general leads to some interesting new `families' of volcanoes. 
The structure theorem \cite[Thm. 4.3]{BrooksJetchevWesolowski17} determines that the $\frl^+$-isogeny graph is a $N(\frl)$-volcano if and only if $\frl^+$ is principal in $\OO_0 \cap K^+$ and $\OO_0^\times\subseteq\OO_{K^+}^\times$, where $\OO_0$ is the endomorphism ring of any abelian variety on the surface. 
In particular, \cite[Rem. 4.14]{BrooksJetchevWesolowski17} shows that the latter condition implies $[\OO_i^\times:\OO_{i+1}^\times] = 1$ for all $i\geq 0$. 
This is not the case for the volcanoes we consider,
see Example~\ref{ex:different_units_ratios} where $[\OO_0^\times:\OO_1^\times] = 4$ and $[\OO_1^\times:\OO_2^\times] = 3$.

Moreover, Main Theorem~\ref{mainthm:C} recovers also the structure theorem of the isogeny graph in the case of supersingular elliptic curves over the prime field, see Example~\ref{ex:ssec}.

\subsection{Structure of the paper}

Our structural results for isogeny graphs of abelian varieties, Main Theorems~\ref{mainthm:B} and~\ref{mainthm:C}, rely on Main Theorem~\ref{mainthm:Bass_cond} and related results on Bass orders. The background in Section~\ref{sec:background} briefly recalls relevant results and notation for orders in \'etale algebras (Section~\ref{ssec:orders}) and abelian varieties over finite fields (Section~\ref{ssec:avs}). Section~\ref{sec:sing_primes_min_oo} continues with orders, describing the role of maximal ideals containing the conductor and the lattice of overorders. We introduce multiplicator ladders of orders in Section~\ref{sec:mult_ladd}, Definition~\ref{def:multiplicatorladder}, which are the backbone for the $(R,\frl)$-isogeny graph. 
From Theorem~\ref{thm:all_ladders} (Main Theorem~\ref{mainthm:Bass_cond}), we immediately recover Corollary~\ref{cor:number_of_ladders}, a result of independent interest giving explicit description of the overorders of a given order with given multiplicator ladder. We further exploit the understanding of the conductor to recover the sizes of the Picard groups of orders in a multiplicator ladder in Section~\ref{sec:classnum}. Theorem~\ref{thm:JoPo_mult_ladder} will be used to explicitly establish the regularity of the related isogeny graphs. 
Following Waterhouse~\cite{Waterhouse}, Section~\ref{sec:avs} describes isogenies of abelian varieties in terms of kernel ideals of their endomorphism rings. This perspective is crucial to our definition of $\frl$-isogeny, Definition~\ref{def:l_isogenies}.  
Section~\ref{sec:p-isog} contains Theorem~\ref{thm:pp-isog}, which unifies the theory of multiplicator ladders from Section~\ref{sec:mult_ladd} with the abelian varieties and isogenies from Section~\ref{sec:avs}. Section~\ref{sec:graph_structure} culminates in Theorem~\ref{thm:graph_structure} (Main Theorem~\ref{mainthm:B}), establishing the general structure of $(R,\frl)$-isogeny graphs whose endomorphism rings contain the order $R$ which is Bass at $\frl$. Section~\ref{sec:volcanoes} analyzes the conditions in which this structure results in a volcano, namely Theorems~\ref{thm:volcano_inertramified} and \ref{thm:volcano_split} (Main Theorem~\ref{mainthm:C}). 
We conclude with Section~\ref{sec:examples}, providing examples and comparison (Comparison~\ref{compBJW}) with \cite{BrooksJetchevWesolowski17} and classical results on elliptic curves.


\section*{Acknowledgments}
The first author was supported by a faculty fellowship with the Commonwealth Cybersecurity Initiative. A portion of this work was completed while the first author was a postdoctoral researcher at Universiteit Leiden supported by the Quantum Software Consortium. The second author was supported by NWO through grant VI.Veni.202.107, by Agence Nationale de la Recherche under the MELODIA project (grant number ANR-20-CE40-0013) and by Marie  Sk{\l}odowska-Curie Actions - Postdoctoral Fellowships 2023 (project 101149209 - AbVarFq).
The third author was partially supported by the Additional Funding Programme for Mathematical Sciences, delivered by EPSRC (EP/V521917/1) and the Heilbronn Institute for Mathematical Research.

\section{Background}\label{sec:background}

\subsection{Orders}\label{ssec:orders}
Let $K$ be an \emph{\'etale algebra} over $\QQ$, that is, a finite product of number fields.
By a \emph{lattice} in $K$, we mean a finitely generated free sub-$\ZZ$-module of the additive group of $K$ with rank equal to $\dim_\QQ(K)$.
Given two lattices $L,L'$ in $K$, the \emph{colon lattice} of $L$ and $L'$ is
\[(L:L') \coloneqq \{x\in K: xL'\subseteq L\}.\]
An \emph{order} $R$ in $K$ is a subring of $K$ which is also a lattice in $K$.
All orders are contained in a unique \emph{maximal order} $\OO_K$, which is the product of the ring of integers of the components of $K$.
An \emph{overorder} of $R$ is an order $T$ in $K$ such that $R\subseteq T$.
They can be computed using the algorithms described in \cite{HofmannSircana20}, or \cite{KirschmerKluners24} when $K$ is a number field.

A \emph{fractional $R$-ideal} is a sub-$R$-module of $K$ which is also a lattice in $K$.
Note that maximal ideals of $R$ and the overorders of $R$ are fractional $R$-ideals.
The colon lattice construction will be of particular importance in the case of ideals: if $I$ and $J$ are fractional $R$-ideals, then $(I:J)$ is also a fractional $R$-ideal. 
Given $I$ is a fractional $R$-ideal, we define its \emph{multiplicator ring} as the order $(I:I)$.
Informally, the order $(I:I)$ is the largest order for which $I$ is an ideal.
A fractional $R$-ideal $I$ is said to be \emph{principal} if there exists a regular element $x\in K$ such that $I=xR$.
If $\frl$ is a maximal ideal of $R$, we denote by $R_\frl$ the \emph{localization} of $R$ at $\frl$.
Similarly, if $I$ is a fractional $R$-ideal, we set $I_\frl \coloneqq I\otimes_R R_\frl$.

A fractional $R$-ideal $I$ is called \emph{invertible} if $I(I:R)=R$, or, equivalently, if it is \emph{locally principal}, that is, for every maximal ideal $\frl$ of $R$ there exists a regular element $x\in K$ such that $I_\frl = (xR)_\frl$. 
The \emph{conductor} $\mathfrak{f}_R=(R:\mathcal{O}_K)$ of an order $R$ in $K$ is the largest $\OO_K$-ideal contained in $R$.
A maximal ideal $\frl$ of an order $R$ is called \emph{singular} if it contains the conductor $\mathfrak{f}_R$ of $R$.
If $\frl$ is not singular we say that it is \emph{regular}.
The \emph{class group} $\Cl(R)$ of an order $R$ is the group (under ideal multiplication) of invertible ideals of $R$ modulo the subgroup of principal fractional $R$-ideals.
In other references, $\Cl(R)$ is called the Picard group of $R$.

\begin{remark}
    In Sections~\ref{sec:sing_primes_min_oo} and~\ref{sec:mult_ladd}, we study the overorders of a given order $R$ with the purpose of applying these results to endomorphism ring of abelian varieties over finite fields.
    For this reason, we develop our theory for lattices and orders over the integers $\ZZ$.
    We nevertheless stress that all our results are valid for lattices and orders over any residually finite Dedekind domain $Z$.
    In this more general context, $K$ would be an \'etale algebra over $Q$, the field of fractions of $Z$.
\end{remark}

\subsection{Abelian varieties over finite fields and isogeny volcanoes}\label{ssec:avs}
In what follows, all morphisms between abelian varieties over $\Fq$ are defined over $\Fq$.
For example, we will write $\End(A)$ for $\End_{\Fq}(A)$.

Throughout, let $A/\Fq$ be an abelian variety with commutative endomorphism ring $\End(A)$.
Let $h(x)$ be the Weil polynomial determining the isogeny class of $A$.
By \cite[Thm~2.(c)]{Tate66}, the polynomial $h(x)$ does not have repeated complex roots.
Hence, $K = \QQ[x]/(h(x))$ is an \'etale algebra over $\QQ$.
Denote by $\pi$ the class of $x$ in $K$.
For every abelian variety $B$ isogenous to $A$, sending the Frobenius endomorphism of $B$ to $\pi$ induces a canonical isomorphism $K \cong \End(B)\otimes_\ZZ\QQ$.
In particular, for every such $B$, we canonically identify $\End(B)$ with an order in $K$.
The \'etale algebra $K$ comes equipped with an involution $\pi\mapsto q/\pi$.
For every $B$ isogenous to $A$, the order $\End(B)$ contains the order $\ZZ[\pi,q/\pi]$ by \cite[Proposition 3.5]{Waterhouse}.

Recall that a $g$-dimensional abelian variety $A$ over a finite field is \emph{ordinary} if $A[p](\bar\FF_p)\cong (\ZZ/p\ZZ)^g$.
If $A[p](\bar\FF_p)\cong (\ZZ/p\ZZ)^{g-1}$ then we say that $A$ is \emph{almost ordinary}.
Unlike in this paper, some authors require an almost ordinary abelian variety to have dimension $>1$.

A convenient way to study how the abelian varieties in a given isogeny class are connected by isogenies is to use the so-called \emph{isogeny graphs}, which are directed graphs whose vertices represent the isomorphism classes of the abelian varieties and the edges represent a particular kind of isogeny.
The notion of an isogeny \emph{volcano} graph was first introduced by Fouquet and Mourain \cite{FouquetMorain02} to describe the structure found in the isogeny graphs of elliptic curves by Kohel \cite{Kohel94}. 
We include the definition as given in \cite{Sutherland_2013}.
The \emph{degree} of a vertex in an undirected graph is the number of edges incident to that vertex, with loops each adding one to this count. 
\begin{definition}[$r$-volcano, level, surface]\label{def:r-volcano}
    Let $r$ be a positive integer. 
    A connected undirected graph $G = (V,E)$ is an $r$-volcano if the set $V$ of vertices may be partitioned into levels $V_0, V_1,\ldots,V_d$ such that the following hold:
    \begin{enumerate}
        \item The subgraph of $G$ induced by the vertices in $V_0$ (the surface) is a regular graph of degree at most two;\label{it:r-volcano-1}
        \item For each $i>0$, each vertex in $V_i$ has precisely one edge connecting to level $V_{i-1}$, and these are the only edges not on the surface;\label{it:r-volcano-2}
        \item  Every vertex not in $V_d$ has degree $r+1$.\label{it:r-volcano-3}
    \end{enumerate} 
\end{definition}

\section{Singular ideals and minimal overorders}
\label{sec:sing_primes_min_oo}
In this section, we study the set of orders containing a given order $R$ in an \'etale algebra $K$. 
In the following lemma, we record a series of well-known facts about regular and singular maximal ideals of $R$.
\begin{lemma}\label{lem:reg_vs_sing}
    Let $R$ be an order and $\frl$ a maximal ideal of $R$.
    Then the following are equivalent:
    \begin{enumerate}
        \item $\frl$ is regular;
        \item $\frl$ is invertible in $R$, that is, $\frl(R:\frl)=R$;
        \item $R = (\frl:\frl)$;
        \item $R$ is $\frl$-maximal, that is, $R_\frl = \OO_{K,\frl}$;
        \item every overorder $T$ of $R$ is $\frl$-maximal.
    \end{enumerate}
    Moreover, if $\frl$ is singular then
    \begin{enumerate}[resume]
        \item \label{lem:reg_vs_sing:mult_rng} $R \subsetneq (\frl:\frl) = (R:\frl)$;
        \item $\frl(R:\frl)=\frl$;
        \item \label{lem:reg_vs_sing:cond} $(R:(\frl:\frl)) = \frl$.
    \end{enumerate}
\end{lemma}
\begin{proof}
    See for example~\cite[Sec.~2]{MarType}.
\end{proof}

\begin{lemma}\label{lem:cond_PP}
    Let $R$ be an order and let $\frl$ be a singular maximal ideal of $R$.
    Set $T = (\frl : \frl)$.
    Let $\ff_R=(R:\OO_K)$ and $\ff_T=(T:\OO_K)$ be the respective conductors. 
    Then $\ff_R = \frl\ff_T$.
\end{lemma}
\begin{proof}
    Note that $\ff_R$ (resp.~$\ff_T$) is maximal among the fractional $\OO_K$-ideals contained in $R$ (resp.~in $T$).
    Since $\frl\ff_T \subseteq \frl T = \frl \subset R $, we deduce that $\frl\ff_T \subseteq \ff_R$.

    We now show the reverse inclusion.
    Firstly, note that $\ff_T\frl = \ff_T \frl \OO_K$, since $\ff_T$ is an $\OO_K$-ideal.
    Secondly, observe that 
    \begin{equation}\label{eq:cond_PP:1}
        \ff_R (\frl\OO_K)^{-1} = (R : \OO_K)(\OO_K : \frl\OO_K) \subseteq (R : \frl\OO_K).
    \end{equation}
    Thirdly, we have
    \begin{equation}\label{eq:cond_PP:2}
    \ff_T = ((\frl : \frl) : \OO_K) = ((R : \frl) : \OO_K) = (R : \frl\OO_K).
    \end{equation}
    Combining Equations~\eqref{eq:cond_PP:1} and \eqref{eq:cond_PP:2}, and multiplying on both sides by $\frl\OO_K$ we obtain the desired inclusion $\ff_R \subseteq \frl\ff_T$.
\end{proof}

\begin{definition}[$\frl$-overorder]\label{def:frl_oo}
    Let $R$ be an order and $\frl$ a maximal ideal of $R$.
    An overorder $S$ of $R$ is an $\frl$-overorder of $R$ if $(R:S)$ is either equal to $R$ or an $\frl$-primary $R$-ideal, that is, $\frl$ is the unique maximal ideal of $R$ containing $(R:S)$.
\end{definition}

\begin{lemma}\label{lem:intermediate_pp_oo}
    Let $R$ be an order, $\frl$ a singular maximal ideal of $R$, and $S$ an $\frl$-overorder of $R$.
    If $S'$ is an order such that $R \subseteq S' \subseteq S$ then $S'$ is also $\frl$-overorder of $R$.
\end{lemma}
\begin{proof}
    Note that $(R:S') \subseteq R$ with equality if and only if $S'=R$.
    If $S'=R$ then $S'$ is $\frl$-primary and we are done.
    If $R\subsetneq S'$ then we have inclusions $(R:S)\subseteq (R:S') \subsetneq R$.
    Since $(R:S)$ is $\frl$-primary then so is $(R:S')$.
\end{proof}

\begin{lemma}
    Let $R$ be an order, $\frl$ a singular maximal ideal of $R$.
    Set $T=(\frl:\frl)$.
    Then $(R:T)=\frl$.
    In particular, $T$ is an $\frl$-overorder of $R$.
\end{lemma}
\begin{proof}
    The statement follows from Lemma~\ref{lem:reg_vs_sing}.\eqref{lem:reg_vs_sing:cond}.
\end{proof}

\begin{definition}[Minimal $\frl$-overorder]\label{def:min_oo}
    Let $R$ be an order and $\frl$ a singular maximal ideal of $R$. 
    Let $S_1$ be an order and $S_2$ be an $\frl$-overorder of $R$.
    We say that $S_2$ is a minimal $\frl$-overorder of $S_1$ if $S_1\subsetneq S_2$ and the inclusions is minimal, that is, if $S'$ an order such that $S_1\subseteq S' \subseteq S_2$ then either $S'=S_1$ or $S'=S_2$.
\end{definition}
Note that the order $S_1$ in Definition~\ref{def:min_oo} is automatically an $\frl$-overorder of $R$ by Lemma~\ref{lem:intermediate_pp_oo}.

\begin{remark}\label{rmk:comparison_defs}
    Let $R$ be an order in an \'etale algebra $K$ and let $\frl$ be a singular maximal ideal of $R$.
    By Definitions \ref{def:frl_oo} and \ref{def:min_oo}, we have that $R$ is an overorder of $R$ and an $\frl$-overorder of $R$, but not a minimal $\frl$-overorder of $R$.
    Hence, our definition of $\frl$-overorder differs from \cite[Def.~5.11]{HofmannSircana20} as it is written, but an author confirmed in private communication that our definition aligns with their intentions and usage in their paper.
\end{remark}

\begin{lemma}\label{lemma:mind_cond_PP}
    Let $R$ be an order and $\frl$ be a singular maximal ideal of $R$.
    Assume that $T=(\frl:\frl)$ is a minimal $\frl$-overorder of $R$.
    Let $\ff_T=(T:\OO_K)$ be the conductor of $T$.
	Then 
    \[ T = R + \ff_T.  \]
\end{lemma}
\begin{proof}
    We have inclusions of orders
    \[ R \subseteq R + \ff_T \subseteq T. \]
    By Lemma~\ref{lem:cond_PP}, the conductor $\ff_T$ of $T$ strictly contains the conductor $\ff_R$ of $R$.
    Since $\ff_R$ is maximal among the $\OO_K$-ideals inside $R$, it follows that $R \subsetneq R + \ff_T$.
    Since $R\subset T$ is minimal it follows that $T=R+\ff_T$.
\end{proof}

\begin{proposition}	\label{prop:min_ext}
    Let $R$ be an order and $\frl$ be a singular maximal ideal of $R$. 
    Assume that $T=(\frl:\frl)$ is a minimal $\frl$-overorder of $R$.
	Then 
    $T$ is the unique minimal $\frl$-overorder of $R$ and
    exactly one of the following statements holds:
	\begin{enumerate}[(i)]
		\item \label{prop:min_ext:inert}
            $\frl$ is a regular maximal ideal of $T$, and the finite field inclusion $R/\frl \hookrightarrow T/\frl$ has degree a prime number bounded by $\dim_\QQ K$;
		\item \label{prop:min_ext:split}
            There are exactly two maximal ideals $\frL_1$ and $\frL_2$ of $T$ containing $\frl$, both regular in $T$ and with residue fields $T/\frL_i$ isomorphic to $R/\frl$;
		\item \label{prop:min_ext:sing_ram}
            There exists $t\in T\setminus R$ such that $t^2\in\frl$ and $T = R[t] = R+tR$, and the unique maximal ideal $\frL$ of $T$ above $\frl$ is $\frL = \frl + tR$. Moreover, the following statements hold:
            \begin{enumerate}[label={(\alph*)}, ref={(\alph*)}]
                \item \label{prop:min_ext:sing_ram:PP} $\frL = \frl + tT$;
                \item \label{prop:min_ext:sing_ram:res_fld} $T/\frL \cong R/\frl$;
                \item \label{prop:min_ext:sing_ram:incl} $\frL^2\subseteq \frl$;
                \item \label{prop:min_ext:sing_ram:ext} $\frl(\frL:\frL) = \frL$;
                \item \label{prop:min_ext:sing_ram:reg} $\frL^2 =\frl$ if and only if $\frL$ is regular.
            \end{enumerate}
	\end{enumerate}
\end{proposition}
\begin{proof}
    By \cite[Prop.~5.2]{HofmannSircana20}, every minimal $\frl$-overorder of $R$ sits inside $T$.
    Hence $T$ is the unique minimal $\frl$-overorder of $R$. 
	It is clear that \ref{prop:min_ext:inert}, \ref{prop:min_ext:split} and \ref{prop:min_ext:sing_ram} are mutually exclusive.
    The inclusions of residue rings $R/\frl \hookrightarrow T/\frl$ is minimal in the sense of \cite[Def.~1.1]{FerrandOlivier70}.
	Hence, by \cite[Lemma 1.2]{FerrandOlivier70}, we have the following mutually exclusive possibilities:
	\begin{enumerate}[(1)]
		\item \label{lemma:min_ext:field} $T/\frl$ is a field, and $\dim_{(R/\frl)}(T/\frl)$ is a prime number;
		\item \label{lemma:min_ext:diag} the inclusion is the diagonal embedding into $T/\frl \cong R/\frl \times R/\frl$; or
		\item \label{lemma:min_ext:nilp} the inclusion is the canonical embedding into $T/\frl \cong (R/\frl)[X]/(X^2)$.
	\end{enumerate}

    Assume we are in case \ref{lemma:min_ext:field}.
    Then $\frl$ is a maximal ideal in its own multiplicator ring $T$.
    Hence, $\frl$ is a regular maximal ideal of $T$ by Lemma~\ref{lem:reg_vs_sing}.
    Moreover, since $R/\frl \hookrightarrow T/\frl$ is a minimal inclusion of finite fields, it has prime degree.
    Let $\ell$ be the rational prime contained in $\frl$.
    Note that $\dim_{(R/\frl)}(T/\frl) \leq \dim_{\FF_\ell}(T/\ell T) = \dim_\QQ K$.
    Therefore, \ref{prop:min_ext:inert} holds. 

    If \ref{lemma:min_ext:diag} holds then there are exactly two maximal ideals $\frL_1$ and $\frL_2$ of $T$ containing $\frl$ which corresponds to the kernels of the two natural surjections
    \[ T/\frl \cong R/\frl \times R/\frl \twoheadrightarrow R/\frl. \]  
    It follows that $T/\frL_1\cong T/\frL_2\cong R/\frl$.
    Hence, $T/\frl \cong T/\frL_1\times T/\frL_2$, which implies $\frl = \frL_1\frL_2$.
    Since the multiplicator ring of $\frL_1\frL_2$ contains both $(\frL_1:\frL_1)$ and $(\frL_2:\frL_2)$, and we assume $T = (\frl : \frl)$ is a minimal $\frl$-overorder, we get 
    \[ T = (\frL_1:\frL_1) = (\frL_2:\frL_2). \]
    By Lemma~\ref{lem:reg_vs_sing}, $\frL_1$ and $\frL_2$ are then regular maximal ideals of $T$.
    Therefore, \ref{prop:min_ext:split} holds.  

    Finally, assume that \ref{lemma:min_ext:nilp} holds.
    Let $t$ be the preimage of $X$ in $T$.
    Then $t^2\in \frl$ and $T=R[t]=R +tR$.
    Also, $T$ has a unique maximal ideal $\frL$ above $\frl$ which is generated as an $R$-module by $\frl$ and $t$. 
    Hence, we are in case \ref{prop:min_ext:sing_ram}.
    Note that $\frL = \frl + tR = \frl + tT$, proving \ref{prop:min_ext:sing_ram:PP}.
    Since $T/\frl$ is a $2$-dimensional vector space over $R/\frl$ and $\frL \neq \frl$, we get that $T/\frL \cong R/\frl$ as stated in \ref{prop:min_ext:sing_ram:res_fld}.
    From $t^2 \in \frl$, we get that $\frL^2 \subseteq \frl$, that is, \ref{prop:min_ext:sing_ram:incl} holds.
    Since 
    \[ 1\leq\dim_{R/\frl} \frL/\frl < \dim_{R/\frl} T/\frl = 2 \] 
    and 
    \[ \frl \subsetneq \frl(\frL:\frL) \subseteq \frL,\] 
    we then have $\frl(\frL:\frL) = \frL$, as in \ref{prop:min_ext:sing_ram:ext}.
    If $\frL^2=\frl$ then the multiplicator ring of $\frL$ is $T$, hence $\frL$ is regular by Lemma~\ref{lem:reg_vs_sing}.
	Conversely, if $\frL$ is regular then $\frL/\frL^2$ is a one dimensional $T/\frL$ vector space.
	Then either $\frL= \frl $ or $\frL^2 = \frl$.
	The first option cannot occur since $t\in T\setminus R$.
    Therefore, \ref{prop:min_ext:sing_ram:reg} holds.
    This completes the proof in case \ref{prop:min_ext:sing_ram}.
\end{proof}

\begin{definition}[Inert, split, ramified, singular maximal ideal]\label{def:min_ext_primes}
    Let $R$, $\frl$ and $T$ be as in Proposition~\ref{prop:min_ext}.
    If \ref{prop:min_ext:inert} holds then we say that $\frl$ is inert in $T$.
    If \ref{prop:min_ext:split} holds then we say that $\frl$ is split in $T$.
    If \ref{prop:min_ext:sing_ram} holds and $\frL^2=\frl$ then we say that $\frl$ is ramified in $T$.
    If \ref{prop:min_ext:sing_ram} holds and $\frL^2 \subsetneq \frl$ then we say that $\frl$ is singular in $T$.
\end{definition}

Example~\ref{ex:all_behaviours} exhibits an order with four singular maximal ideals covering the four possible cases of Definition~\ref{def:min_ext_primes}.

\begin{corollary}
    Let $R$ be an order and $\frl$ be a singular maximal ideal of $R$.
    Assume that $T=(\frl:\frl)$ is the unique minimal $\frl$-overorder.
    Then, $\frl$ is inert, split or ramified in $T$ if and only if $T$ is $\frl$-maximal, that is, $T_\frl = \OO_\frl$.
\end{corollary}
\begin{proof}
    By Proposition~\ref{prop:min_ext}, $\frl$ is inert, split or ramified in $T$ if and only if the maximal ideals of $T$ above $\frl$ are regular.
    The statement then follows by Lemma~\ref{lem:reg_vs_sing}.
\end{proof}

\begin{definition}[Cohen-Macaulay type, Gorenstein, Bass order]\label{def:Gor_Bass}
    Let $R$ be an order and $\frl$ a maximal ideal of $R$.
    The dimension of $(R:\frl)/R$ as a vector spaces over $R/\frl$ is the Cohen-Macaulay type of $R$ at $\frl$.
    We say that $R$ is Gorenstein at $\frl$ if the Cohen-Macaulay type of $R$ at $\frl$ is one.
    Let $S$ be an overorder of $R$.
    We say that $S$ is Gorenstein at $\frl$ if $S$ is Gorenstein at every maximal ideal $\frL$ of $S$ above $\frl$.
    We say that $R$ is Bass at $\frl$ if every overorder $S$ of $R$ is Gorenstein at $\frl$. 
    We say that $R$ is Gorenstein (resp.~Bass) if it is so at every maximal ideal $\frl$.
\end{definition}

Note that to check whether $R$ is Bass at $\frl$ it suffices to check if its $\frl$-overorders are Gorenstein at $\frl$, for the reason we now explain.
Let $\mm_1,\ldots,\mm_n$ be the singular maximal ideals of $R$ and $S$ be any overorder of $R$.
By \cite[Thm.~5.13]{HofmannSircana20}, we can write $S=S_1+\ldots+S_n$ where $S_i$ is a $\mm_i$-overorder of $R$.
Say that $\frl = \mm_1$.
Since $S_\frl=S_{1,\frl}$, we see that $S$ is Gorenstein at $\frl$ if and only if $S_1$ is Gorenstein at $\frl$.
Hence, if every $\frl$-overorder of $R$ are Gorenstein at $\frl$ then every overorder of $R$ is Gorenstein at $\frl$. The reverse implication is trivial.

For more details about the Cohen-Macaulay type and other equivalent definitions of Gorenstein and Bass, we refer the reader to \cite[Sec.~3 and 4]{MarType}.

\begin{lemma}\label{lem:Gor_T_minimal}
    Let $R$ be an order and $\frl$ be a singular maximal ideal of $R$.
    If $R$ is Gorenstein at $\frl$ then $T=(\frl:\frl)$ is the unique minimal $\frl$-overorder of $R$.
\end{lemma}
\begin{proof}
    See \cite[Prop.~5.21]{HofmannSircana20}.
\end{proof}
The next example shows that the reverse implication does not hold.
\begin{example}\label{ex:jumps}
    Consider the endomorphism rings of abelian varieties in the isogeny class \href{https://abvar.lmfdb.xyz/Variety/Abelian/Fq/3/5/c_ab_ae}{3.5.c\_ab\_ae}\footnote{Our LMFDB links direct to the `xyz' version of the LMFDB which currently (\today) has more information about the isogeny classes of abelian varieties than the standard LMFDB. This may change with time, at which point these links will be obsolete, but the labels will remain accurate.} \cite{lmfdb}, which is defined by the Weil polynomial $h(x) = x^6+2x^5-x^4-4x^3-5x^2+50x+125$. 
    Set $K\coloneqq \QQ[x]/(h(x)) = \QQ[\pi]$ and $R\coloneqq \ZZ[\pi,5/\pi]$. The lattice of overorders of $R$ is depicted in Figure~\ref{fig:non-Gorenstein}. 
    The order $\mathbf{T_1}$ of index $4$ in the maximal order $\mathcal{O}_K$ is not Gorenstein, but the maximal order $\OO_K$ is indeed the unique overorder of $\mathbf{T_1}$ and $\mathcal{O}_K$ is the multiplicator ring of the unique singular maximal ideal $\frl_{\mathbf{T_1}}$.

    This example illustrates other interesting phenomena, for example the order $\mathbf{S}$ is the unique minimal overorder of $\mathbf{T_3}$, but $\mathbf{S}$ is not the multiplicator ring of the unique singular maximal ideal $\frl_{\mathbf{T_3}}$ of $\mathbf{T_3}$. The orders $T_2, T_2', T_2''$ are all minimal overorders of $\mathbf{S}$, but only one of these orders is the multiplicator ring of a singular maximal ideal of an order below, namely the order we denote $T_2$ is $(\frl_{\mathbf{T_3}}:\frl_{\mathbf{T_3}})$.
    
    \begin{figure}
        \centering
	    \begin{tikzpicture}
            \node[] (R) at (-1,0) {$R$};
            \node[] (T3) at (-1,1) {$\mathbf{T_3}$};
            \node[] (S) at (-1,2) {$\mathbf{S}$};
            \node[] (T2) at (-1,3) {$T_2$};
            \node[] (T2') at (-2,3) {$T_2'$};
            \node[] (T2'') at (0,3) {$T_2''$};
            \node[] (T1) at (-1,4) {$\mathbf{T_1}$};
            \node[] (O) at (-1,5) {$\mathcal{O}_K$};
            \node[] (header) at (-1,6) {Order lattice};
            
            \node[] (pR) at (4,0) {$(R:\mathcal{O}_K)=\mathfrak{f}_R=\mathcal{L}^4$};
            \node[] (pT1) at (4,1) {$(\mathbf{T}_1:\mathcal{O}_K)=\mathfrak{f}_\mathbf{T_3}=\mathcal{L}^3$};
            \node[] (pS) at (4,2) {$(\mathbf{S}:\mathcal{O}_K)=\mathfrak{f}_\mathbf{S}=\mathcal{L}^2$};
            \node[] (pT2) at (4,3) {$(T_2:\mathcal{O}_K)=\mathfrak{f}_{T_2}=\mathfrak{f}_{T'_2}=\mathfrak{f}_{T_2''}=\mathcal{L}^2$};
            \node[] (pT3) at (4,4) {$(\mathbf{T_1}:\mathcal{O}_K)=\mathfrak{f}_\mathbf{T_1}=\mathcal{L} = \mathcal{O}_K\frl_{\mathbf{T_1}}$};
            \node[] (c) at (4,6) {Conductors};
            
            \node[] (pR) at (10,0) {$R$};
            \node[] (pT3) at (10,1) {$\mathbf{T_3}$};
            \node[] (pS) at (10,2) {$\mathbf{S}$};
            \node[] (pT2) at (10,3) {$T_2$};
            \node[] (pT2p) at (8,3) {$T'_2$};
            \node[] (pT2pp) at (9,3) {$T''_2$};
            \node[] (pT1) at (10,4) {$\mathbf{T_1}$};
            \node[] (OK) at (10,5) {$\mathcal{O}_K$};
            \node[] (s) at (10,6) {Mult. Ring Relations};
            
            \draw[-] (R) -- (T3) node[midway,left] {\small{$2$}};
            \draw[-] (T3) -- (S) node[midway,left] {\small{$2$}};
            \draw[-] (S) -- (T2) node[midway,left] {\small{$2$}};
            \draw[-] (S) -- (T2') node[midway,left] {\small{$2$}};
            \draw[-] (S) -- (T2'') node[midway,right] {\small{$2$}};
            \draw[-] (T2) -- (T1) node[midway,left] {\small{$2$}};
            \draw[-] (T2') -- (T1) node[midway,left] {\small{$2$}};
            \draw[-] (T2'') -- (T1) node[midway,right] {\small{$2$}};
            \draw[-] (T1) -- (O) node[midway,left] {\small{$4$}};
            
            \draw[->] (pT2) to (pT1);
            \draw[->] (pT2p) to (pT1);
            \draw[->] (pT2pp) to (pT1);
            \draw[->] (10.3,1) to [out = 0,in=0] (10.2,3);
            \draw[->] (10.3,2) to [out = 0,in=0] (10.2,4);
            \draw[->] (pT1) to (OK);
            \draw[->] (pR) to (pT3);
        \end{tikzpicture}
    \caption{Lattice of overorders of $R\coloneqq\ZZ[\pi,5/\pi]$ in the isogeny class \href{https://abvar.lmfdb.xyz/Variety/Abelian/Fq/3/5/c_ab_ae}{3.5.c\_ab\_ae} discussed in Example~\ref{ex:jumps}. Orders with Cohen-Macaulay type 2 are depicted in bold $\mathbf{T_1},\mathbf{S},\mathbf{T_3}$, and non-bold orders are Gorenstein. Arrows indicate multiplicator rings of singular maximal ideals, for example $\mathbf{T_1} = (\frl_{T_2}:\frl_{T_2}) = (\frl_{\mathbf{S}}: \frl_{\mathbf{S}})$.}
    \label{fig:non-Gorenstein}
    \end{figure}
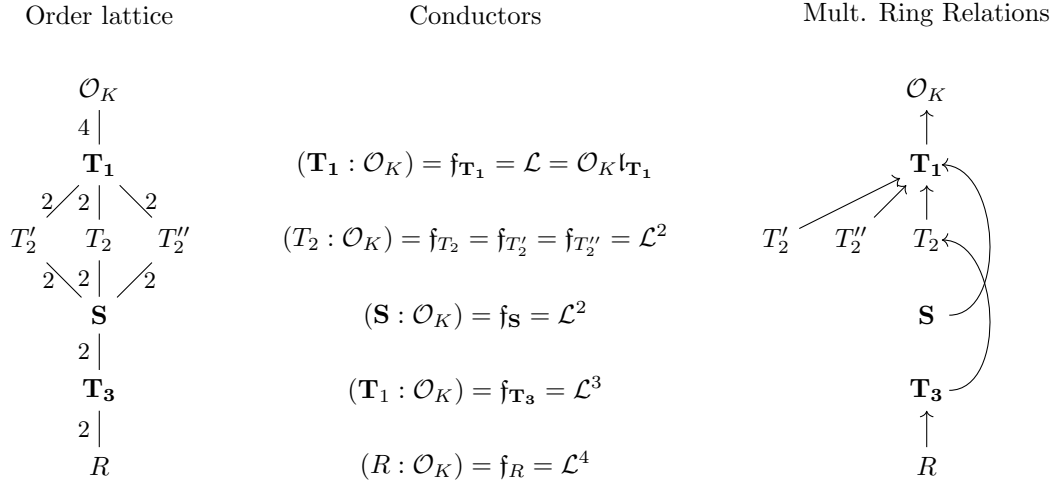
\end{example}

\begin{example}\label{ex:all_behaviours}
    Consider the endomorphism rings of abelian varieties in the isogeny class \href{https://abvar.lmfdb.xyz/Variety/Abelian/Fq/3/25/g_cg_ji}{3.25.g\_cg\_ji}, which is defined by the irreducible Weil polynomial
    $h(x)=x^6 + 6x^5 + 58x^4 + 242x^3 + 1450x^2 + 3750x + 15625$.
    Set $K\coloneqq\QQ[x]/(h(x))=\QQ[\pi]$.
    The order $R\coloneqq\ZZ[\pi,25/\pi]$ has four singular maximal ideals, at which $R$ is Gorenstein (in fact, Bass): 
    \begin{itemize}
        \item $\frl_2$ with $[R:\frl_2]=2$ which is singular in $(\frl_2:\frl_2)$;
        \item $\frl_4$ with $[R:\frl_4]=4$ which is split in $(\frl_4:\frl_4)$;
        \item $\frl_3$ with $[R:\frl_3]=3$ which is inert in $(\frl_3:\frl_3)$;
        \item $\frl_7$ with $[R:\frl_7]=7$ which is ramified in $(\frl_7:\frl_7)$.
    \end{itemize}
\end{example}

\section{Multiplicator ladders}\label{sec:mult_ladd}
Let $R$ be an order in an \'etale algebra $K$ and let $\frl$ be a maximal ideal of $R$.
In this section we study the set of $\frl$-overorders of $R$ when such orders form what we call the $\frl$-multiplicator ladder of $R$.

\begin{definition}[$\frl$-multiplicator ladder of $R$]\label{def:multiplicatorladder}
    We say that a chain of inclusion 
    \[ R=R_d \subsetneq R_{d-1} \subsetneq \ldots \subsetneq R_0 \]
    of overorders of $R$ is the $\frl$-multiplicator ladder of $R$ if the following conditions hold:
    \begin{itemize}
        \item $R=R_d, R_{d-1}, \ldots, R_0$ are all the $\frl$-overorders of $R$;
        \item for each $i=1,\ldots,d$, the order $R_i$ has a unique maximal ideal $\frL_i$ above $\frl$ and $\frL_i$ is singular;
        \item for each $i=1,\ldots,d$, the order $R_{i-1}$ equals $(\frL_i:\frL_i)$ and is the unique minimal $\frl$-overorder of $R_i$.
    \end{itemize}
\end{definition}

\begin{proposition}\label{prop:Bass_tot_ord}
    Let $R$ be an order and let $\frl$ be a maximal ideal of $R$.
    Consider the following statements.
    \begin{enumerate}[(i)]
        \item \label{prop:Bass_tot_ord:Bass} $R$ is Bass at $\frl$.
        \item \label{prop:Bass_tot_ord:tot_ord}
           $R$ has an $\frl$-multiplicator ladder.
        \item \label{prop:Bass_tot_ord:cond}
            The set of $\frl$-overorders of $R$ is totally ordered by inclusion: 
            \[ R=R_d \subsetneq R_{d-1} \subsetneq \ldots \subsetneq R_0. \]
            Moreover, for each $i=0,\ldots,d-1$, if $\ff_{R_i} = (R_i:\OO_K)$ is the conductor of $R_i$ then
            \[ R_i = R_{i+1} + \ff_{R_i} = R + \ff_{R_i}. \]
    \end{enumerate}
    We have the following implications
    \[ \ref{prop:Bass_tot_ord:Bass} \Longrightarrow \ref{prop:Bass_tot_ord:tot_ord} \Longrightarrow \ref{prop:Bass_tot_ord:cond}. \]
\end{proposition}
\begin{proof}
    Both implications are trivially true if $\frl$ is a regular ideal of $R$, so assume that $\frl$ is singular.
    The first implication is a recursive application of Lemma~\ref{lem:Gor_T_minimal} and Proposition~\ref{prop:min_ext}.
    We now prove the second implication.
    Let $R=R_d \subsetneq R_{d-1} \subsetneq \ldots \subsetneq R_0$ be the $\frl$-multiplicator ladder of $R$.
    Clearly, the set of $\frl$-overorders is totally ordered by inclusion.
    Fix an index $1\leq i \leq d$.
    By applying Lemma~\ref{lemma:mind_cond_PP} to the inclusion $R_{i+1}\subset R_i$, we get 
    \[ R_i = R_{i+1} + \ff_{R_i}. \]
    Inductively,
    \[ R_i = R + \ff_{R_d} + \ldots + \ff_{R_i}. \]
    To complete the proof, it is enough to observe that we have inclusions
    $\ff_{R_d} \subset  \ldots \subset \ff_{R_i}$.
\end{proof}

\begin{remark}\label{rmk:generalizeHofmannSircana}
    Proposition~\ref{prop:Bass_tot_ord} is a refinement of \cite[Cor.~5.25]{HofmannSircana20}.
    There the authors prove that, given an order $R$ which is Bass at a maximal ideal $\frl$, the chain of inclusions of $\frl$-overorders of $R$ splits into two totally ordered chains if there are two maximal ideals above $\frl$.
    Proposition~\ref{prop:Bass_tot_ord} tells us that if this is the case, then we are already at the top of the chain. 
    Hence the chain of $\frl$-overorders cannot not split.
\end{remark}

The following proposition shows that statements \ref{prop:Bass_tot_ord:cond} and \ref{prop:Bass_tot_ord:Bass} in Proposition~\ref{prop:Bass_tot_ord} are `almost' equivalent.
\begin{proposition}\label{prop:ladder_conseq}
    Let $\frl$ be a singular maximal ideal of $R$.
    Assume that $R$ has an $\frl$-multiplicator ladder, say $R=R_d \subsetneq R_{d-1} \subsetneq \ldots \subsetneq R_0$.
    Then: 
    \begin{enumerate}[(i)]
        \item \label{prop:ladder_conseq:primesattop} The $\frl$-maximal order $R_0$ has at most two maximal ideals above $\frl$, which are regular.
        \item \label{prop:ladder_conseq:Gor} For each $i=2,\ldots,d$, the order $R_i$ is Gorenstein at $\frl$.
        \item \label{prop:ladder_conseq:T1} The order $R_1$ is Gorenstein at $\frl$ if its unique maximal ideal $\frL_1$ above $\frl$ is ramified or split in $R_0$, while if $\frL_1$ is inert in $R_0$ then the Cohen-Macaulay type of $R_1$ at $\frL_1$ equals $\dim_{R_1/\frL_1}(R_0/\frL_1) - 1$.
    \end{enumerate}
\end{proposition}
\begin{proof}
    Part \ref{prop:ladder_conseq:primesattop}, follows from Proposition~\ref{prop:min_ext}.
    Fix an index $i>0$ and let $\frL_i$ be unique maximal ideal of $R_i$ above $\frl$.
    Then the Cohen-Macaulay type of $R_i$ at $\frL_i$ is the dimension of
    \[ \frac{(R_i:\frL_i)}{R_i} = \frac{R_{i-1}}{R_i}\cong \frac{R_{i-1}/\frL_i}{R_i/\frL_i} \]
    over $R_i/\frL_i$, where the first equality follows from Lemma~\ref{lem:reg_vs_sing}.\eqref{lem:reg_vs_sing:mult_rng}.
    If $i>1$ then by Proposition~\ref{prop:min_ext} we have that $R_{i-1}/\frL_i \cong (R_i/\frL_i)[X]/(X^2)$ which has dimension $2$ over $R_i/\frL_i$, which implies \ref{prop:ladder_conseq:Gor} by \cite[Main Theorem 1]{MarType}.
    Similarly, if $i=1$ and $\frL_1$ is ramified or split in $R_0$ then $R_0/\frL_1$ has dimension $2$ over $R_1/\frL_1$, showing that $R_1$ is Gorenstein at $\frL_1$.
    Finally, if $i=1$ and $\frL_1$ is inert in $R_0$ then $R_0/\frL_1$ is a field extension of $R_1/\frL_1$ of prime degree $r$.
    Then the Cohen-Macaulay type of $R_1$ at $\frL_1$ equals $r-1$.
\end{proof}

Examples~\ref{ex:tot_ord_not_imply_Bass} and~\ref{ex:cond_not_imply_tot_ord} show that the reverse implications in Proposition~\ref{prop:Bass_tot_ord} do not hold.

\begin{example}\label{ex:tot_ord_not_imply_Bass}
    Consider the isogeny class \href{https://abvar.lmfdb.xyz/Variety/Abelian/Fq/3/4/ab_d_ah}{3.4.ab\_d\_ah} which is defined by the irreducible Weil polynomial $h(x)=x^6 - x^5 + 3x^4 - 7x^3 + 12x^2 - 16x + 64$.
    Set $K\coloneqq\QQ[x]/(h(x)) = \QQ[\pi]$ and $R\coloneqq\ZZ[\pi,4/\pi]$.
    The order $R$ has two singular maximal ideals, both above $2$ and complex conjugate of each other,  that is, mapped to each other by the involution induced by $\pi\mapsto 4/\pi$.
    Let $\frl$ be one of them.
    The set of $\frl$-overorders of $R$ is totally ordered and consists of the orders: 
    \begin{equation}\label{eq:ex:tot_ord_not_imply_Bass:mult_ladder}
        R \subsetneq T \subsetneq \OO.
    \end{equation}
    The order $T$ (resp.~$\OO$) has a unique maximal ideal $\frl_T$ (resp.~$\frl_\OO$) above $\frl$.
    On the one hand, we have that 
    \[ T = (\frl:\frl) \quad \text{and} \quad \OO = (\frl_T : \frl_T), \]
    which means that Proposition~\ref{prop:Bass_tot_ord}.\ref{prop:Bass_tot_ord:tot_ord} holds for the order $R$ and the maximal ideal $\frl$, that is, \eqref{eq:ex:tot_ord_not_imply_Bass:mult_ladder} is the $\frl$-multiplicator ladder of $R$.
    On the other hand, the order $T$ is not Gorenstein at $\frl_T$ and hence the order $R$ is not Bass at $\frl$, that is, Proposition~\ref{prop:Bass_tot_ord}.\ref{prop:Bass_tot_ord:Bass} does not hold.
\end{example}

\begin{example}\label{ex:cond_not_imply_tot_ord}
    Consider the isogeny class \href{https://abvar.lmfdb.xyz/Variety/Abelian/Fq/2/3/a_ac}{2.3.a\_ac} which is defined by the irreducible Weil polynomial $h(x)=x^4 - 2x^2 + 9$.
    Set $K\coloneqq\QQ[x]/(h(x)) = \QQ[\pi]$ and $R\coloneqq\ZZ[\pi,\overline{\pi}]$.
    Let $\OO_K$ be the maximal order of $K$.
    The unique maximal ideal $\frl$  of $R$ above $2$ is the unique singular maximal ideal of $R$.
    The lattice of $\frl$-overorders of $R$ is totally ordered
    \begin{equation}\label{eq:ex:cond_not_imply_tot_ord:mult_ladder}
        R \subsetneq R_2 \subsetneq R_1 \subsetneq \OO_K
    \end{equation}
    On the one hand, we have that 
    \[ R_2 = R + \ff_{R_2} \quad \text{and} \quad R_1 = R + \ff_{R_1} = R_2 + \ff_{R_1}. \]
    This implies that Proposition~\ref{prop:Bass_tot_ord}.\ref{prop:Bass_tot_ord:cond} holds true for the order $R$.
    On the other hand, \eqref{eq:ex:cond_not_imply_tot_ord:mult_ladder} is not a multiplicator ladder, that is, Proposition~\ref{prop:Bass_tot_ord}.\ref{prop:Bass_tot_ord:tot_ord} doesn't hold for the order $R$ and maximal ideal $\frl$, since we have
    \[ R_2 = (\frl:\frl) \quad \text{and} \quad \OO_K = (\frl_2:\frl_2) = (\frl_1:\frl_1), \]
    where $\frl_2$ (resp.~$\frl_1$) is the unique singular maximal ideal of $R_2$ (resp.~$R_1$).
\end{example}

\begin{proposition}\label{prop:ladder}\
    Let $\frl$ be a singular maximal ideal of $R$. 
    Assume that $R$ has $\frl$-multiplicator ladder $R=R_d \subsetneq R_{d-1} \subsetneq \ldots \subsetneq R_0$.
    For $1\leq i\leq d$ let $\frl_i$ be the unique maximal ideal of $R_i$ above $\frl$.
    \begin{enumerate}[(i)]
        \item \label{prop:ladder:pp_ext} For $i=1,\ldots,d$, we have $\frl R_{i-1}  = \frl_i$.
        \item \label{prop:ladder:primes_ext} For $i=1,\ldots,d$, we have $\frl_i R_0 = \frl_1$.
        \item \label{prop:ladder:cond} For $i=0,\ldots,d$, the conductor $\ff_{R_i}=(R_i:\OO_K)$ of $R_i$ in $\OO_K$ satisfies $\ff_{R_i} = \frl_1^i \ff_{R_0} = \frl^i \ff_{R_0}$.
    \end{enumerate}
\end{proposition}
\begin{proof}
    For $i=1,\ldots,d-1$, we have $\frl_{i+1}(\frl_i:\frl_i) = \frl_{i+1}R_{i-1} = \frl_i$, by Proposition~\ref{prop:min_ext}.\ref{prop:min_ext:sing_ram:ext}.
    Since $R_{d-1} = (\frl:\frl)$ and hence $\frl = \frl R_{d-1}$, we inductively get $\frl_i = \frl R_{i-1}$ as claimed by~\ref{prop:ladder:pp_ext}.

    Applying \ref{prop:ladder:pp_ext}, for each $1\leq i \leq d$, we get 
    \[ \frl_1 = \frl R_{0} = \frl R_{i-1} R_0 = \frl_i R_0,\]
    as stated in \ref{prop:ladder:primes_ext}.

    We prove \ref{prop:ladder:cond} by induction.
    For $i=0$ the statement is clearly true.
    So we assume that it holds for $i-1$, that is,
    \[ \ff_{R_{i-1}} = (R_{i-1} : \OO_K ) = \frl_1^{i-1} \ff_{R_0} = \frl^{i-1} \ff_{R_0}. \]
    By Lemma~\ref{lem:cond_PP}, we have that $\ff_{R_i} = \frl_i\ff_{R_{i-1}}$.
    Since $\ff_{R_0} = R_0\ff_{R_0}$, we get
    \[ \ff_{R_i} = \frl_i \frl_1^{i-1} \ff_{R_0} = \frl_i R_0 \frl_1^{i-1} \ff_{R_0} = \frl_1^i\ff_{R_0} = \frl^i\ff_{R_0}, \]
    where the last two equalities are applications of~\ref{prop:ladder:primes_ext} and~\ref{prop:ladder:pp_ext}, respectively.
\end{proof}
\begin{corollary}\label{cor:order_at_level}
    Let $\frl$ be a singular maximal ideal of $R$. 
    Assume that $R$ has $\frl$-multiplicator ladder $R=R_d \subsetneq R_{d-1} \subsetneq \ldots \subsetneq R_0$.
    Then, for each $i=0,\ldots,d$, we have
    \[ R_i = R + \frl^i\ff_{R_0}.\]
\end{corollary}
\begin{proof}
    The statement follows from combining Proposition~\ref{prop:Bass_tot_ord}.\ref{prop:Bass_tot_ord:cond} and Proposition~\ref{prop:ladder}.\ref{prop:ladder:cond}.
\end{proof}
\begin{remark}\label{rmk:ladder_length}
    Let $\frl$ be a singular maximal ideal of $R$. 
    Assume that $R$ has $\frl$-multiplicator ladder $R=R_d \subsetneq R_{d-1} \subsetneq \ldots \subsetneq R_0$.
    The extension $\frl\OO_K$ of $\frl$ to the maximal order $\OO_K$ factors as $\pp^2$, $\mathfrak{r}\tilde{\mathfrak{r}}$, or $\qq$ depending on whether the unique maximal ideal of $R_1$ above $\frl$ is ramified, split or inert in $R_0$, see Definition~\ref{def:min_ext_primes}.
    In each of these three cases, we can relate the length of the ladder $d$ with the valuation of the conductor $\ff$ of $R$ in $\OO_K$ at the maximal ideals of $\OO_K$ above $\frl$: $d=\val_{\pp}(\ff)/2$ in the ramified case; $d=\val_{\mathfrak{r}}(\ff)=\val_{\tilde{\mathfrak{r}}}(\ff)$ in the split case; $d=\val_{\qq}(\ff)$ in the inert case.
\end{remark} 
The following result says that if $R$ has a multiplicator ladder at a singular maximal ideal $\frl$ then every overorder of $R$ belongs to a unique multiplicator ladder `above $\frl$'.
This is an expanded and more general version of Main Theorem~\ref{mainthm:Bass_cond} from Section~\ref{sec:intro}.

\begin{theorem}\label{thm:all_ladders}
    Let $\frl$ be a singular maximal ideal of $R$. 
    Assume that $R$ has an $\frl$-multiplicator ladder, denoted $R=R_d \subsetneq R_{d-1} \subsetneq \ldots \subsetneq R_0$.
    Let $T$ be an overorder of $R$.
    Then there exists a unique overorder $\OO$ of $R$ such that:
    \begin{enumerate}[(i)]
        \item \label{thm:all_ladders:base_loc} $\OO_\frl = R_\frl$;
        \item \label{thm:all_ladders:ladder} $\OO$ has an $\frL$-multiplicator ladder $\OO=\OO_d\subsetneq \OO_{d-1} \subsetneq \ldots \subsetneq \OO_0$ where $\frL=\frl\OO$ is the unique singular maximal ideal of $\OO$ above $\frl$;
        \item \label{thm:all_ladders:T} there exists a unique index $0\leq i\leq d$ such that $T=\OO_i$.
    \end{enumerate}
    Moreover,
    \begin{enumerate}[(a)]
        \item \label{thm:all_ladders:all_loc} for $i=0,\ldots,d$, we have $R_{i,\frl}=\OO_{i,\frl}$; 
        \item \label{thm:all_ladders:max_idl} for $i=1,\ldots,d$, the unique singular maximal ideal $\frL_i$ of $\OO_i$ above $\frL$ satisfies $\frL_i=\frl\OO_{i-1}$ and $\OO_i/\frL_i\cong R/\frl$;
        \item \label{thm:all_ladders:max_idl_OO1} $\frL_1$ is split (resp.~inert, resp.~ramified) in $\OO_0$ if and only if the unique maximal ideal $\frl_1$ of $R_1$ is split (resp.~inert, resp.~ramified) in $R_0$, and the corresponding residue fields are isomorphic;
        \item \label{thm:all_ladders:cond} for $i=0,\ldots,d$, the conductor $\ff_{\OO_i}=(\OO_i:\OO_K)$ satisfies $\ff_{\OO_i}=\frL^i\ff_{\OO_0} = \frl^i\ff_{\OO_0}$ and $\OO_i = \OO + \frl^i\ff_{\OO_0} = R +\frl^i\ff_{\OO_0}$.
    \end{enumerate}
\end{theorem}
\begin{proof}
    Let $\frl=\mm_1,\mm_2,\ldots,\mm_n$ be the singular maximal ideals of $R$.
    By \cite[Thm.~5.13]{HofmannSircana20}, for each $1\leq k \leq n$, there exists a unique $\mm_k$-overorder $S_k$ of $R$ such that
    $T = S_1 + \ldots + S_n$.
    Note that $S_{1,\frl}=T_\frl$, and that, for each $1 < j \leq n$, we have $S_{j,\frl} = R_\frl$.
    Also, there exist a unique index $0 \leq i \leq d$ such that $S_{1} = R_{i}$.
    For each $0\leq i\leq d$, define $\OO_i = R_i + S_2 + \ldots + S_n$.
    It follows that the order $\OO \coloneqq \OO_d$ satisfies the conditions \ref{thm:all_ladders:base_loc}, \ref{thm:all_ladders:ladder} and \ref{thm:all_ladders:T}, and that \ref{thm:all_ladders:all_loc} holds.
    Part \ref{thm:all_ladders:max_idl}, \ref{thm:all_ladders:max_idl_OO1} and \ref{thm:all_ladders:cond} follow from part~\ref{thm:all_ladders:all_loc}, Proposition \ref{prop:ladder} and Corollary~\ref{cor:order_at_level}.
\end{proof}

\begin{definition}[$\frl$-multiplicator ladder in the set of overorders]\label{def:frl_mult_ladd_in_oo}
    We call a chain 
    \[ \OO_d\subsetneq \OO_{d-1} \subsetneq \ldots \subsetneq \OO_0 \]
    as in Theorem~\ref{thm:all_ladders}.\ref{thm:all_ladders:ladder} an $\frl$-multiplicator ladder in the set of overorders of $R$.
    For each index $i$, we say that the order $\OO_i$ is at level $i$.
\end{definition}

Theorem~\ref{thm:all_ladders} tells us that every overorder of $R$ belongs to a unique $\frl$-multiplicator ladder.
We deduce that if two $\frl$-multiplicator ladders intersect, then they must coincide.
\begin{corollary}
    Let $\OO=\OO_d \subsetneq \ldots \subsetneq \OO_0$ and $\OO'=\OO'_d \subsetneq \ldots \subsetneq \OO'_0$ be $\frl$-multiplicator ladders in the set of overorders of $R$.
    If there are indices $0\leq i,j \leq d$ such that $\OO_i = \OO'_j$ then we have $\OO_k=\OO'_k$ for every index $0\leq k\leq d$.
\end{corollary}

The following corollary counts the number of distinct $\frl$-multiplicator ladders in the set of overorders of $R$.
\begin{corollary}\label{cor:parallel_ladders}
    Let $\mm_1,\ldots,\mm_n$ be the singular maximal ideals of $R$, and let $\frl\coloneqq\mm_1$.
    For each $i=1,\ldots,n$, let $\mathcal{S}_i$ be the set of $\mm_i$-overorders of $R$.  
    Assume that $R$ has an $\frl$-multiplicator ladder.
    The number of distinct $\frl$-multiplicator ladders in the set of overorders of $R$ equals $\prod_{i=2}^n \vert \mathcal{S}_i \vert$.
\end{corollary}
\begin{proof}
    By Theorem~\ref{thm:all_ladders}, the number of distinct $\frl$-multiplicator ladders in the set of overorders of $R$ equals the number $N$ of overorders $\OO$ of $R$ such that $\OO_\frl = R_\frl$.
    By \cite[Thm.~5.13]{HofmannSircana20}, each overorder $T$ of $R$ can be written as $T=T_1+\ldots+T_n$ for a unique tuple $(T_1,\ldots,T_n)$ in the cartesian product $\prod_{i=1}^n \mathcal{S}_i$.
    Hence, $T_\frl = R_\frl$ if and only if $T_1 = R$.
    It follows that $N = \prod_{i=2}^n \vert \mathcal{S}_i \vert$.
\end{proof}

\begin{lemma}\label{lem:not_inv_d_2}
    Let $\OO_d\subsetneq \ldots \subsetneq \OO_0 $ be an $\frl$-multiplicator ladder in the set of overorders of $R$ with $d\geq 1$.
    If $\frl\OO_1$ is not an invertible $\OO_1$-ideal then $d\in \{1,2\}$ and $R$ is not Bass at $\frl$.
\end{lemma}
\begin{proof}
    By Theorem~\ref{thm:all_ladders}.\ref{thm:all_ladders:ladder}, we know that $\OO_1$ is the multiplicator ring of $\frl\OO_1$.
    So the order $\OO_1$ is not Gorenstein at $\frl$, and $R$ is not Bass at $\frl$.
    
    We show that $d\leq 2$ by contraposition.
    If $d>2$ then for each $i=2,\ldots,d$ the order $\OO_i$ is Gorenstein at $\frl$ by Proposition~\ref{prop:ladder_conseq}.
    For $i=2,\ldots,d-1$, $\OO_i$-ideal $\frl\OO_i$ is invertible since it has multiplicator ring $\OO_i$.
    Hence, $\frl \OO_1 = (\frl \OO_2)\OO_1$ is an invertible $\OO_1$-ideal as well.
\end{proof}

We conclude this section with a result of independent interest which will not be used in the rest of the paper.
Let $R$ be an order and let $\mm_1,\ldots,\mm_n$ the be the singular maximal ideals of $R$.
Contrary to the rest of the paper where we work with one maximal ideal at the time, we assume for the rest of this section that $R$ has an $\mm_i$-multiplicator ladder $R=R_{\mm_i,d_i}\subsetneq \ldots\subsetneq R_{\mm_i,0}$ for each $0\leq i \leq n$.
The extension $\mm_i\OO_K$ of each singular maximal ideal $\mm_i$ to the maximal order $\OO_K$ factors as 
$\pp_i^2$, $\mathfrak{r}_i\tilde{\mathfrak{r}}_i$, or 
$\qq_i$ depending on  whether the unique maximal ideal of $R_{\mm_i,1}$ is ramified, split or inert in $R_{\mm_i,0}$, see Definition~\ref{def:min_ext_primes}.
After possibly relabelling, assume that $\mm_1,\ldots,\mm_{n_1}$ are ramified, $\mm_{n_1+1},\ldots,\mm_{n_2}$ are split and that $\mm_{n_2+1},\ldots,\mm_n$ are inert.
In the following corollary, we generalize \cite[Thm~6.13]{ChoHongLee24_arXiv}.
There the authors require $R$ to be Bass and an integral domain.
\begin{corollary}\label{cor:number_of_ladders}
    Assume that $R$ has an $\mm_i$-multiplicator ladder $R=R_{\mm_i,d_i}\subsetneq \ldots\subsetneq R_{\mm_i,0}$ for each $0\leq i \leq n$.
    Then each overorder $S$ of $R$ is of the form  
    \[ S = R + \left(\prod_{i=1}^n \mm_i^{k_i}\OO_K\right) =
    R + 
    \left(\prod_{i=1}^{n_1} \pp_i^{2k_i}\right)
    \left(\prod_{i=n_1+1}^{n_2} (\mathfrak{r}_i\tilde{\mathfrak{r}}_i)^{k_i}\right)
    \left(\prod_{i=n_2+1}^{n} \qq_i^{k_i}\right)
    , \]
    for unique integers $k_1,\ldots,k_n$ satisfying $0 \leq k_i \leq d_i$.
    Moreover, the conductor $(S:\OO_K)$ of $S$ is 
    \[ (S:\OO_K) = \left(\prod_{i=1}^n \mm_i^{k_i}\OO_K\right) = 
    \left(\prod_{i=1}^{n_1} \pp_i^{2k_i}\right)
    \left(\prod_{i=n_1+1}^{n_2} (\mathfrak{r}_i\tilde{\mathfrak{r}}_i)^{k_i}\right)
    \left(\prod_{i=n_2+1}^{n} \qq_i^{k_i}\right). \]
    In particular, the total number of overorders is $\prod_{i=1}^n (d_i+1)$ and, for each $j=1,\ldots,n$, the number of $\mm_j$-multiplicator ladders is $\prod_{i \neq j}(d_i+1)$.
\end{corollary}
\begin{proof}
    All statements follow from Theorem~\ref{thm:all_ladders}.\ref{thm:all_ladders:cond}.
\end{proof}

\begin{example}\label{ex:mult_ladders}
    We return to the setting of Example~\ref{ex:all_behaviours}: the isogeny class \href{https://abvar.lmfdb.xyz/Variety/Abelian/Fq/3/25/g_cg_ji}{\textnormal{3.25.g\_cg\_ji}}, defined by the irreducible Weil polynomial
    $h(x)=x^6 + 6x^5 + 58x^4 + 242x^3 + 1450x^2 + 3750x + 15625$.
    As before, let $K=\QQ[x]/(h(x))=\QQ(\pi)$ and $R=\ZZ[\pi,25/\pi]$. Recall that $R$ is Bass, has $4$ singular maximal ideals (namely $\mathfrak{l}_2$, $\mathfrak{l}_4$, $\mathfrak{l}_3$, and $\mathfrak{l}_7$), and $[\OO_K:R] = 2^4\cdot 3\cdot 7$. 
    In Figure~\ref{fig:MinimalOverOrders}, we depict the multiplicator ladders of $R$ at the four singular maximal ideals, as described by Corollary~\ref{cor:number_of_ladders}.
    In particular, we deduce that $R$ has $24$ overorders.
    \begin{figure}[h!]
        \begin{tikzpicture}
        \node[] (R) at (0,0) {$R$};
        \node[] (T7) at (-4.5,1) {$R + \frl_3\frl_4\frl_2^2\OO_K$};
        \node[] (T3) at (-1.5,1) {$R + \frl_7\frl_4\frl_2^2\OO_K$};
        \node[] (T4) at (1.5,1) {$R + \frl_7\frl_3\frl_2^2\OO_K$};
        \node[] (T2) at (4.5,1) {$R + \frl_7\frl_3\frl_4\frl_2\OO_K$};
        \node[] (T22) at (5,2) {$R + \frl_7\frl_3\frl_4\OO_K$};
        \draw[-] (R) -- (T7) node[midway,below] {\small{$\mathfrak{l}_7$}};
        \draw[-] (R) -- (T3) node[midway,right] {\small{$\mathfrak{l}_3$}};
        \draw[-] (R) -- (T4) node[midway,left] {\small{$\mathfrak{l}_4$}};
        \draw[-] (R) -- (T2) node[midway,below] {\small{$\mathfrak{l}_2$}};
        \draw[-] (T2) -- (T22) node[midway,below] {};
        \end{tikzpicture}
        \caption{Multiplicator ladders at the singular maximal ideals of $R$ as in Example~\ref{ex:mult_ladders}.}
        \label{fig:MinimalOverOrders}
    \end{figure}
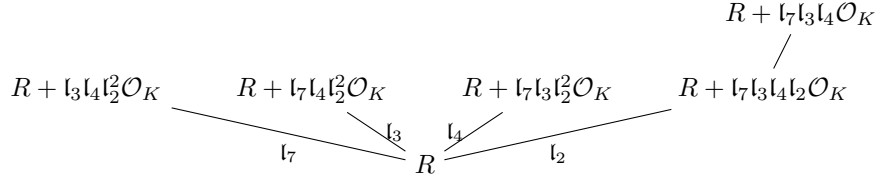
\end{example}

\section{Class numbers and levels}\label{sec:classnum}
As we work with orders in multiplicator ladders, we will be interested in how the class numbers change in the ladder. We recall the following generalization of classical class number formulas for orders in number fields; cf.~\cite[Chapter~1.12]{neukirch}. 
\begin{proposition}[Prop.~6.10, \cite{JP20}]
\label{prop:JoPo}
	Let $K$ be an \'etale algebra. If $R\subseteq \OO_K$ is an order in $K$, then
	$$
		\frac{\#\Cl(R)}{\#\Cl(\OO_K)}
			=\frac{[\OO_K : R]}{[\OO_K^\times : R^\times]} 
			\prod_{\pp} \frac{\prod_{\PP\mid \pp}(1 - \#(\OO_K/\PP)^{-1})}{1 - \#(R/\pp)^{-1}}
	$$
 \end{proposition}
 
 By exploiting Proposition~\ref{prop:min_ext}, we can give an analogue of Proposition~\ref{prop:JoPo} in the case of minimal $\frl$-overorders. We use the notions of inert, split, ramified and singular maximal ideals from Definition~\ref{def:min_ext_primes}.
 
\begin{proposition}
\label{prop:JoPo_min_overorder}
	Let $R$ be an order and $\frl$ a singular maximal ideal of $R$. Assume that $T = (\frl : \frl)$ is the minimal $\frl$-overorder of $R$. 
	\begin{enumerate}[(a)]
	\item If $\frl$ is inert in $T$ and $f$ is the prime extension degree of $T/\frl$ over $R/\frl$, then
	$$
		\frac{\#\Cl(R)}{\#\Cl(T)}
			=\frac{[T:R]}{[T^\times : R^\times]}\frac{(1- \#(R/\frl)^{-f})}{(1- \#(R/\frl)^{-1})}
	$$
	\item If $\frl$ is split in $T$, then
	$$
		\frac{\#\Cl(R)}{\#\Cl(T)}
			=\frac{[T:R]}{[T^\times : R^\times]}(1- \#(R/\frl)^{-1})
	$$
	\item If $\frl$ is ramified or singular in $T$, then
	$$
		\frac{\#\Cl(R)}{\#\Cl(T)}
			=\frac{ [T:R]}{[T^\times : R^\times]}
	$$
	\end{enumerate} 
\end{proposition}
\begin{proof}
    We apply Proposition~\ref{prop:JoPo} to both $T$ and $R$ and take the quotient. The only part which requires inspection is the term which is a product over maximal ideals $\pp$ of $T$ or $R$, respectively.
    Set
    \[ E(T,R)\coloneqq \frac{\prod_{\pp\subset R} \frac{\prod_{\PP\mid \pp}(1 - \#(\OO_K/\PP)^{-1})}{1 - \#(R/\pp)^{-1}}}{\prod_{\pp\subset T} \frac{\prod_{\PP\mid \pp}(1 - \#(\OO_K/\PP)^{-1})}{1 - \#(T/\pp)^{-1}}}. \]
    Since $T_\pp = R_\pp$ for all maximal ideals $\pp\neq\frl$ we have
    \[ E(T,R) = \frac{\prod_{\frL \mid \frl }(1 - \#(T/\frL)^{-1})}{1 - \#(R/\frl)^{-1}}. \]

    Recall the results of Proposition~\ref{prop:min_ext}, in the language of Definition~\ref{def:min_ext_primes}. If $\frl$ is inert in $T$, then $\frl$ is a maximal ideal of $T$ and $R/\frl \hookrightarrow T/\frl$ is a field extension of prime degree $f\leq \operatorname{rank}_\ZZ(R)$.
    If $\frl$ is split, there are two maximal ideals $\frL_1$ and $\frL_2$ of $T$ lying over $\frl$ and we have $T/\frL_1\cong T/\frL_2\cong R/\frl$.
    If $\frl$ is ramified or singular in $T$, then there is precisely one maximal ideal $\frL$ of $T$ lying over $\frl$ and we have $T/\frL \cong R/\frl$. 
    It follows that
    \[
        E(T,R) = 
        \begin{cases}
            (1- \#(R/\frl)^{-f})/ (1- \#(R/\frl)^{-1}) & \text{ if } \frl \text{ is inert};\\
            1-\#(R/\frl)^{-1} & \text{ if } \frl \text{ is split};\\
            1 & \text{ if } \frl \text{ is ramified or singular}.
        \end{cases}
    \]  
    The result now follows from Proposition~\ref{prop:JoPo}.
\end{proof}
 
In the next theorem, we compare the sizes of the class groups of two order at consecutive levels of the same $\frl$-multiplicator ladder in the set of overorders of a given order $R$.
\begin{theorem}\label{thm:JoPo_mult_ladder}
    Fix an order $R$ and a singular maximal ideal $\frl$ of $R$.

    Assume that $R$ has an $\frl$-multiplicator ladder.
    Let $\OO_d\subsetneq \dots \subseteq \OO_0$ be an $\frl$-multiplicator ladder in the set of overorders of $R$. 
    \begin{enumerate}[(a)]
        \item For $2 \leq i \leq d$, then we have
        $$
            \frac{\#\Cl(\OO_i)}{\#\Cl(\OO_{i-1})}
                =\frac{ \#(R/\frl)}{[\OO_{i-1}^\times : \OO_i^\times]}
        $$
        \item If $R$ is Bass at $\frl$, then
        $$
            \frac{\#\Cl(\OO_1)}{\#\Cl(\OO_{0})}
                =\frac{ \#(R/\frl)-\delta_\frl}{[\OO_{0}^\times : \OO_1^\times]}
        $$
        where $\delta_\frl =-1,0$, or $1$ if $\frl$ is inert, ramified, or split in $\OO_0$, respectively.
    \end{enumerate} 
\end{theorem}
 
\begin{proof} 
    First, recall that for $1 \leq i\leq d$, there is a unique (singular) maximal ideal $\frL_i$ of $\OO_i$ lying over $\frl$ and we have $\#(\OO_i/\frL_i) =\#(R/\frl)$, see Theorem~\ref{thm:all_ladders}.\ref{thm:all_ladders:max_idl}. Therefore, we deduce that
    $$
    [\OO_{i-1} : \OO_{i}] 
        = \#(\OO_{i-1}/\OO_{i})
        = \#((\frL_{i} : \frL_i)/\OO_{i})
        = \#(\OO_i/\frL_i) 
        =\#(R/\frl),
    $$
    where the penultimate equality follows from the fact that $\OO_i$ is Gorenstein at $\frL_i$.

    Assume that $\frL_1$ is inert in $\OO_0$.  By Proposition~\ref{prop:ladder_conseq}.\ref{prop:ladder_conseq:T1}, the Cohen-Macaulay type of $\OO_1$ at $\frL_1$ is equal to $f-1$, where $f = \dim_{\OO_1/\frL_1}(\OO_0/\frL_1)$ is the residue field extension degree. Because we assume that $R$ is Bass at $\frl$, it follows that $\OO_1$ has type 1 and $f = 2$. 
    Observe that 
    $$
    \frac{(1-\#(R/\frl)^{-2})}{(1-\#(R/\frl)^{-1})} = (1+\#(R/\frl)^{-1}).
    $$
    We now immediately obtain the desired result from Proposition~\ref{prop:JoPo_min_overorder}.
\end{proof}

\section{Abelian varieties over finite fields and kernel ideals}\label{sec:avs}
From now on, unless otherwise specified, when we write a morphism of abelian varieties which are defined over $\Fq$, we always mean an $\Fq$-morphism.
Consequently, we write $\End(A)$ for $\End_\Fq(A)$.
We begin by recalling the definition of kernel ideal introduced in \cite[Section 3.2]{Waterhouse}.
This notion is used to define well-behaved isogenies between abelian varieties.

Fix an isogeny class $\Ical$ of abelian varieties over $\Fq$ with commutative endomorphism algebra.
Identify the isomorphism algebras of each $A$ in $\Ical$ with $K=\QQ[x]/(h(x))=\QQ[\pi]$
where $h(x)$ is the Weil polynomial characterizing $\Ical$.
Then each $\End(A)$ is naturally identified with an order in $K$ by sending the Frobenius endomorphism of $A$ to $\pi$.
Fix an order $R$ in $K$ and denote by $\Ical_R$ the subset of $\Ical$ consisting of abelian varieties $A$ such that $R \subseteq \End(A)$.
For any fractional $R$-ideal $I\subseteq R$ and abelian variety $A \in \Ical_R$, define 
\[ A[I] \coloneqq \bigcap\ker\alpha, \]
where the (scheme-theoretic) intersection is taken over the elements $\alpha$ of $I\cdot \End(A)$, the extension of $I$ to $\End(A)$.

\begin{definition}[Kernel ideal]
    An ideal $I$ of $\End(A)$ is a kernel ideal for $A$ if and only if $I = \{\alpha\in\End(A):\alpha A[I] = 0\}$.
\end{definition}

The following lemma provides another way of characterizing kernel ideals.
\begin{lemma}\label{lemma:ker_idl_is_biggest}
    An ideal $I$ of $\End(A)$ is a kernel ideal for $A$ if and only if $I$ is the largest among all ideals $J$ of $\End(A)$ satisfying $A[J] = A[I]$.
\end{lemma}
\begin{proof}
    The largest ideal satisfying $A[J] = \cap_{\alpha\in I}\ker\alpha$ is $\{\alpha \in\End(A):\alpha A[I] = 0\}$, so $I$ is a kernel ideal if and only if $I=\{\alpha \in\End(A):\alpha A[I] = 0\}$.
\end{proof}

In general, kernel ideals are hard to classify.
Nevertheless, there are some classes of ideals that give rise to kernel ideals, as we now explain.
\begin{lemma}\label{lem:ker_idl_max}
    Every maximal ideal of $\End(A)$ is a kernel ideal for $A$.
\end{lemma}
\begin{proof}
    This follows from the characterization of kernel ideals in Lemma~\ref{lemma:ker_idl_is_biggest}.
\end{proof}

\begin{definition}[Divisorial ideal]
    An ideal $I$ of $R$ is divisorial if $(R:(R:I)) = I$.
\end{definition}

\begin{proposition}\label{prop:divisorial_kernel}
    Every divisorial ideal of $\End(A)$ is a kernel ideal for $A$.
    If $\End(A)$ is Gorenstein, then every fractional ideal $I$ with $(I:I)=\End(A)$ is a kernel ideal for $A$.
\end{proposition}
\begin{proof}
    The first statement is \cite[Remark 7(d)]{Kani11}.
    The second statement follows from the fact that every fractional ideal of a Gorenstein order is divisorial, see \cite[Section 2.6]{BuchmannLenstra94}.
\end{proof}

We now define and study a special kind of isogenies defined via ideal multiplication.
Under certain assumptions on the ideals, we are able to control the endomorphism ring of the target of the corresponding isogeny.
\begin{definition}[$I$-multiplication]\label{def:frl-multiplication}
    Let $A$ be an abelian variety in $\Ical_R$ with endomorphism ring $T \coloneqq \End(A)$.
    Let $I$ be a fractional $R$-ideal contained in $T$.
    The $I$-multiplication from $A$ is the projection $A\to A/A[IT]$.
\end{definition}

\begin{definition}[Ascending, descending, horizontal isogeny]
    Let $\varphi: A \to B$ be an isogeny of abelian varieties in $\Ical$.
    We say that $\varphi$ is
    \begin{itemize}
        \item ascending if $\End(A) \subsetneq \End(B)$;
        \item descending if $\End(B) \subsetneq \End(A)$;
        \item horizontal if $\End(A) = \End(B)$.
    \end{itemize}
\end{definition}

\begin{lemma}\label{lem:Imult}
    Let $A$ be an abelian variety in $\Ical_R$ with endomorphism ring $T\coloneqq \End(A)$.
    Let $I$ be a fractional $R$-ideal contained in $T$ and $\varphi_I : A\to A/A[IT]$ be the corresponding $I$-multiplication.
    Then:
    \begin{enumerate}[(1)]
        \item\label{lem:Imult:gen} if $IT$ is a kernel ideal for $A$ then $\varphi_I$ is horizontal when $(IT:IT)=T$ and ascending otherwise;
        \item\label{lem:Imult:inv} if $IT$ is invertible in $T$, that is, $IT(T:IT)=T$, then $\varphi_I$ is horizontal;
        \item\label{lem:Imult:sing_max} if $IT$ is a singular maximal ideal of $T$ then $\varphi_I$ is ascending and $\End(A/A[IT]) = (IT: IT)$.
    \end{enumerate}
\end{lemma}
\begin{proof}
    By  \cite[Prop.~3.9]{Waterhouse}, $(IT:IT)\subseteq\End(A/A[IT])$, with equality if $IT$ is a kernel ideal for $A$. This implies \ref{lem:Imult:gen}.
    If $IT$ is invertible in $T$ then it is divisorial and $(IT:IT)=T$.
    Then \ref{lem:Imult:inv} follows from Proposition~\ref{prop:divisorial_kernel} and \ref{lem:Imult:gen}.
    Finally, if $IT$ is a maximal ideal of $T$ then it is a kernel ideal by Lemma~\ref{lem:ker_idl_max}.
    If moreover $IT$ is singular then $T \subsetneq (IT:IT)$ by Lemma~\ref{lem:reg_vs_sing}.\eqref{lem:reg_vs_sing:mult_rng}.
    Then \ref{lem:Imult:sing_max} follows again from \cite[Prop.~3.9]{Waterhouse}.
\end{proof}

\begin{corollary}\label{cor:all_frl_oo_are_end}
    Let $\frl$ be a maximal ideal of $R$. 
    Assume that $R$ has an $\frl$-multiplicator ladder.
    If $A \in \Ical_R$ then for every $\frl$-overorder $T$ of $\End(A)$ there exists $B\in\Ical_R$ with $\End(B)=T$.
\end{corollary}
\begin{proof}
    Set $S\coloneqq\End(A)$.
    By Theorem~\ref{thm:all_ladders}, the order $S$ belongs to an $\frl$-multiplicator ladder in the set of overorders of $R$.
    If $S$ is at level $0$ of the $\frl$-multiplicator ladder then the statement is true.
    If $S$ is at level $>0$, then the unique maximal ideal $\frL$ of $S$ above $\frl$ is singular.
    Lemma~\ref{lem:Imult}.\ref{lem:Imult:sing_max} states that the endomorphism ring of $A/A[\frL]$ equals $(\frL:\frL)$, which is the unique minimal $\frl$-overorder of $S$.
    Iterating the procedure, we get the statement.
\end{proof}

The following well-known proposition connects ideals in the class group of the endomorphism ring of an abelian variety with isogenies to abelian varieties with isomorphic endomorphism rings.
\begin{proposition}\label{prop:Cl_acts}
    Let $A$ be an abelian variety in $\Ical$ and set $S\coloneqq\End(A)$.
    Then $\Cl(S)$ acts freely by $I$-multiplication on the set of isomorphism classes of abelian varieties $B$ in $\Ical$ with $S=\End(B)$.
\end{proposition}
\begin{proof}
    Let $I$ be an invertible $S$-ideal contained in $S$.
    Then $I$ is divisorial and hence a kernel ideal for $A$.
    Also, $A/A[I]$ has endomorphism ring $S$ by Lemma~\ref{lem:Imult}.\ref{lem:Imult:inv}.
    Moreover, if $J$ is an invertible $S$-ideal contained in $S$, then $A/A[I]\cong A/A[J]$ if and only if $I\cong J$ by \cite[Thm.~3.11]{Waterhouse}.
    This shows that $\Cl(S)$ acts on the set of isomorphism classes of abelian varieties with endomorphism ring $S$ and that this action is free.
\end{proof}

The following two propositions are similar in spirit, but logically independent.
They will be used in later sections to determine the size of the kernels of certain isogenies.
\begin{proposition}\label{prop:GorEnd}
    Let $\frl$ be a maximal ideal of $T \coloneqq \End(A)$ which is coprime to $q$.
    Assume that $T$ is Gorenstein at $\frl$.
    Then we have a $T$-linear isomorphism $A[\frl]\cong T/\frl$.
\end{proposition}
\begin{proof}
    Let $\ell$ be the rational prime contained in $\frl$.
    Because $\ell\nmid q$, multiplication by $\ell$ is separable and we identify $A[\ell]$ with the corresponding torsion subgroup in $A(\overline\FF_q)$.
    By \cite[Theorem~3.2]{MS25}, we have an isomorphism $A[\ell]\cong T/\ell T$ of $T$-modules which sends $A[\frl]$ to the biggest submodule of $T/\ell T$ annihilated by $\frl$, namely $(\ell T: \frl)/\ell T$.
    Hence
    \[ A[\frl] \cong \frac{(\ell T: \frl)}{\ell T} \cong \frac{(T: \frl)}{T} \cong \frac{T}{\frl},  \]
    where the last isomorphism is \cite[Thm.~6.3.(4)]{Bass63}.
\end{proof}

\begin{proposition}\label{prop:GorEnd_general}
    Let $I$ be a fractional ideal of $T \coloneqq \End(A)$ such that $I\subseteq T$ and which is coprime to $q$.
    Assume that $T$ is Gorenstein at all maximal ideals containing $I$ and that $(I:I)=T$.
    Then we have a $T$-linear isomorphism $A[I]\cong T/I$.
\end{proposition}
\begin{proof}
    Let $m$ be a rational integer inside $I$.
    Because $\gcd(m, q) = 1$, multiplication by $m$ is separable and we identify $A[m]$ with the corresponding torsion subgroup in $A(\overline\FF_q)$.
    By \cite[Theorem~3.2]{MS25}, we have an isomorphism $A[m]\cong T/m T$ of $T$-modules which sends $A[I]$ to the biggest submodule of $T/m T$ annihilated by $I$, which is $(m T: I)/m T$ by definition.
    Hence
    \[ A[I] \cong \frac{(m T: I)}{m T} \cong \frac{(T: I)}{T}. \]
    We will now show that $(T: I)/T \cong T/I$ concluding the proof.
    Consider the finite length $T$-modules $M \coloneqq (T: I)/T$ and $N \coloneqq  T/I$.
    Let $\mathcal{S}$ be the finite set of maximal ideals $\mm$ of $T$ such that $M_\mm \neq 0$ or $N_\mm \neq 0$.
    Observe that this is precisely the set of maximal ideals containing $I$.
    We have canonical $T$-linear isomorphisms $M\cong \oplus_{\mm\in\mathcal S}M_\mm$ and $N\cong \oplus_{\mm\in\mathcal S}N_\mm$.
    To conclude, it suffices to show that for each $\mm\in\mathcal S$ we have a $T$-linear isomorphism $M_\mm\cong N_\mm$.
    Fix $\mm\in \mathcal S$.
    By assumption $T$ is Gorenstein at $\mm$ and hence $I$ is locally principal at $\mm$, say $I_\mm = xT_\mm$, see \cite[Lemma~2.2 and Proposition~3.4]{MarType}.
    Therefore,
    \[ M_\mm = \frac{(T_\mm : I_\mm)}{T_\mm} = \frac{(T_\mm : xT_\mm)}{T_\mm} \cong \frac{(T_\mm : T_\mm)}{xT_\mm} = \frac{T_\mm}{I_\mm} = N_\mm, \]
    as required.
\end{proof}

\section{\texorpdfstring{$\frl$}{l}-isogenies and multiplicator ladders}\label{sec:p-isog}
Fix an isogeny class $\Ical$ of abelian varieties over $\Fq$ with commutative endomorphism algebra $K=\QQ[\pi]$.
Let $R$ be an order in $K$ and fix a maximal ideal $\frl$ of $R$.
In this section, we will assume that $R$ has an $\frl$-multiplicator ladder $R=R_d\subsetneq \ldots \subsetneq R_0$.
Note that $\frl$ is singular if and only if $d>0$.

\begin{proposition}\label{prop:level_AV}
    Let $A$ be an abelian variety in $\Ical_R$.
    Fix an integer $0\leq i \leq d$. 
    Then the following are equivalent:
    \begin{itemize}
        \item $\End(A)_\frl = R_{i,\frl}$;
        \item $\End(A) = \OO_i$ where $\OO_d\subsetneq \dots \subsetneq \OO_0$ is an $\frl$-multiplicator ladder in the set of overorders of $R$;
        \item $\left( \End(A) \right)_\frl = (R + \frl^i\ff_{R_0})_\frl$.
    \end{itemize} 
\end{proposition}
\begin{proof}
    It follows from Theorem~\ref{thm:all_ladders}, parts \ref{thm:all_ladders:all_loc}, \ref{thm:all_ladders:T}, \ref{thm:all_ladders:cond}, respectively. 
\end{proof}

\begin{definition}[Level]\label{def:level}
    Fix an integer $0\leq i \leq d$. 
    We say that an abelian variety $A$ in $\Ical_R$ is at the level $i$ (with respect to $R$ and $\frl$) if any of the equivalent statements in Proposition~\ref{prop:level_AV} hold.
\end{definition}

The next proposition describes how endomorphism rings distribute along a multiplicator ladder.
\begin{proposition}\label{prop:same_num_AV_along_ladder}
    If $\OO_d\subsetneq \dots \subsetneq \OO_0$ is any $\frl$-multiplicator ladder in the set of overorders of $R$, then exactly one of the following statements is true:
    \begin{enumerate}[(a)]
        \item \label{prop:same_num_AV_along_ladder:yesEnds} There exists an integer $\dmin$ satisfying $d\geq \dmin \geq 0$ such that the order $\OO_i$ arises as the endomorphism ring for an abelian variety in $\Ical_R$ if and only if $\dmin\geq i\geq 0$.
        If $\OO_{\dmin}$ is Bass at $\frl$, then there exists a positive integer $N$ so that, for every index $\dmin\geq i\geq 0$, the number $n(\OO_i)$ of isomorphism classes of abelian varieties $A \in \Ical_R$ with $\End(A)=\OO_i$   is
        \[
            n(\OO_i) = N\cdot \#\Cl(\OO_i).
        \]
        \item \label{prop:same_num_AV_along_ladder:noEnds} None of the orders $\OO_d,\dots, \OO_0$ arise as $\End(A)$ for any $A\in\Ical_R$.
    \end{enumerate}
    Moreover,   we have the implications
        $\ref{prop:same_num_AV_along_ladder:ord_alm_ord_Fp}\Longrightarrow
        \ref{prop:same_num_AV_along_ladder:Omin}\Longrightarrow
        \ref{prop:same_num_AV_along_ladder:laddersOmin}$  between the following statements.
    \begin{enumerate}[(i)]
        \item \label{prop:same_num_AV_along_ladder:ord_alm_ord_Fp} The isogeny class $\Ical$ is ordinary, almost ordinary, or over a prime field $\FF_p$.
        \item \label{prop:same_num_AV_along_ladder:Omin} There exists a unique overorder $\OO_{\min}$ of $R$ whose overorders are precisely the all endomorphism rings of the abelian varieties in $\Ical_R$.
        \item \label{prop:same_num_AV_along_ladder:laddersOmin} An $\frl$-multiplicator ladder satisfies \ref{prop:same_num_AV_along_ladder:yesEnds} if and only if it contains an overorder of the order $\OO_{\min}$ defined in \ref{prop:same_num_AV_along_ladder:Omin}, and the constant $\dmin$ is the same for all ladders satisfying \ref{prop:same_num_AV_along_ladder:yesEnds}.
    \end{enumerate}
\end{proposition}
\begin{proof}
    For the first part of the proposition, it suffices to show that if there exists an abelian variety in $A \in \Ical_R$ such that $\End(A)=\OO_i$ then \ref{prop:same_num_AV_along_ladder:yesEnds} holds.
    Assume there exists such an $A$ and let $\dmin$ be the largest integer in $d\geq d_{\min}\geq 0$ at which this occurs.
    Pick any abelian variety $A_{\dmin}$ such that $\End(A_{\dmin})=\OO_{\dmin}$.
    If $\dmin=0$ then we are done.
    If not, let $A_{\dmin-1}$ be the target of the ascending $\frl$-isogeny from $A_{\dmin}$.
    Then $A_{\dmin-1}$ is at level $\dmin-1$, that is, $\End(A_{\dmin-1})=\OO_{\dmin-1}$.
    Repeating this process shows that an order $\OO_i$ arises as the endomorphism ring for an abelian variety in $\Ical_R$ if and only if $d_{\min}\geq i\geq 0$, as claimed in \ref{prop:same_num_AV_along_ladder:yesEnds}.

    For every $0\leq i\leq d_{\min}$, let $n(\OO_i)$ be the number of abelian varieties in $\mathcal I_R$ with endomorphism ring $\OO_i$.
    Since $\Cl(\OO_i)$ acts freely on this set of abelian varieties by Proposition~\ref{prop:Cl_acts}, we can write $n(\OO_i) = N_i\cdot \#\Cl(\OO_i)$ for some integer $N_i$. 
    Our goal is to show that $N_i = N_0$ for all $d_{\min }\geq i\geq 0$.
    This is done by writing $N_i$ in terms of local information, using \cite[Prop.~8.1]{BergKarMar24_arXiv}, from which we borrow the following notation.
    If there is no maximal ideal of $R$ containing both $\pi$ and $q/\pi$ then set $d(\OO_i)=1$.
    If there is a maximal ideal $\pp$ of $R$ containing both $\pi$ and $q/\pi$, then let $d(\OO_i)$ be the number of isomorphism classes of the local-local parts of the Dieudonn\'e modules of the abelian varieties in $\Ical$ with endomorphism ring $\OO_i$.
    Note that if a maximal ideal $\pp$ as above exists, then it is unique by \cite[Prop.~4.14]{BergKarMar24_arXiv}.
    Let $w(\OO_i)$ be the product of number of isomorphism classes of the $\ell$-adic Tate modules of the abelian varieties in $\Ical$ with endomorphism ring $\OO_i$ for each $\ell\neq p$ multiplied with the number of isomorphism classes of the non-local-local parts of the Dieudonn\'e modules of the abelian varieties in $\Ical$ with endomorphism ring $\OO_i$.
    By \cite[Prop.~3.8, Cor.~4.18]{BergKarMar24_arXiv}, we can write
    \[ w(\OO_i) = \prod_{\mm \neq \pp } w_\mm(\OO_i), \]
    where, for each maximal ideal $\mm$ of $R$ distinct from $\pp$, we define $w_\mm(\OO_i)$ as the number of fractional $R$-ideals $I$ with $(I:I)=\OO_i$ modulo the relation $I\otimes_R R_\mathfrak{m} \cong_{R_\mathfrak{m}} J\otimes_R R_\mathfrak{m}$.
    Now, by \cite[Prop.~8.1]{BergKarMar24_arXiv}, we obtain
    \begin{equation}\label{eq:Ni}
        N_i =  d(\OO_i) \cdot w(\OO_i) = d(\OO_i) \cdot \prod_{\mm \neq \pp } w_\mm(\OO_i).        
    \end{equation}
    Observe that $d(\OO_i)$ (resp.~$w_\mm(\OO_i)$) depends only on the localization of $\OO_i$ at $\pp$ (resp.~$\mm$), see \cite[Lem.~3.5.(ii), Prop.~6.1, Lem.~6.4]{BergKarMar24_arXiv}.
    Since, $\OO_{\dmin}$ is Bass at $\frl$, it is well known that $w_\frl(\OO_i)=1$ for each $d_{\min} \geq i \geq 0$, see for example \cite[Prop.~2.7]{BuchmannLenstra94} for the case of integral domains or \cite[Prop.~3.4]{MarType} for the general case.
    For every maximal ideal $\mathfrak{m}\neq \frl$, we have that $\OO_i\otimes_R R_{\mm} = \OO_0\otimes_R R_{\mm}$ for each $d_{\min} \geq i \geq 0$ by the definition of an $\frl$-multiplicator ladder.
    Therefore, $w_\mm(\OO_i) = w_\mm(\OO_0)$ and $d(\OO_i) = d(\OO_0)$ for each $d_{\min} \geq i \geq 0$.
    Hence, $N_i=N_0$ thus completing the proof of the first part of the proposition.

    We now move to the second part of the proposition.
    The implication \ref{prop:same_num_AV_along_ladder:Omin}$\Rightarrow$\ref{prop:same_num_AV_along_ladder:laddersOmin} is a consequence of Theorem~\ref{thm:all_ladders}.
    Now assume that \ref{prop:same_num_AV_along_ladder:ord_alm_ord_Fp} holds.
    If $\Ical$ is ordinary or defined over $\FF_p$ the order $\ZZ[\pi,q/\pi]$ then the endomorphism rings of the abelian varieties in $\Ical$ are precisely the overorders of $\ZZ[\pi,q/\pi]$. 
    This is \cite[Thm.~7.4,Thm.~6.1]{Waterhouse} in the simple case, while it follows from \cite{Deligne,CentelegheStixI} in general, see also \cite[Thm.~4.3]{Mar_sqfree}.
    Hence, the order $R[\pi,q/\pi]$ satisfies the defining property of $\OO_{\min}$ stated in \ref{prop:same_num_AV_along_ladder:Omin}.
    If $\Ical$ is almost ordinary, by the equivalence of categories developed in \cite[Thm.~1.1]{OswalShankar} for the simple case and extended in \cite[Rmk.~2.11]{BergKarMar} to the commutative endomorphism algebra case, the unique overorder $S$ of $R[\pi,q/\pi]$ satisfying $\OO_{\pp} = \OO_{K,\pp}$, where $\pp$ is the unique maximal ideal of $R$ containing both $\pi$ and $q/\pi$, namely $\pp = (p,\pi,q/\pi)$, and $\OO_{\qq} = R_{\qq}$ for every other maximal ideal $\qq$, satisfies the defining property of $\OO_{\min}$ stated in \ref{prop:same_num_AV_along_ladder:Omin}.
    This completes the proof of the implication \ref{prop:same_num_AV_along_ladder:ord_alm_ord_Fp}$\Rightarrow$\ref{prop:same_num_AV_along_ladder:Omin}.
\end{proof}

\begin{corollary}
    \label{cor:O_min_Bass}
    Assume that the isogeny class $\Ical$ is ordinary or over the prime field.
    Define $\OO_{\min}$ as in Proposition~\ref{prop:same_num_AV_along_ladder}.\ref{prop:same_num_AV_along_ladder:Omin}.
    Then the following statements are equivalent:
    \begin{enumerate}[(i)]
        \item \label{cor:O_min_Bass:1orbit} For each overorder $\OO$ of $\OO_{\min}$, there is precisely one $\Cl(\OO)$-orbit of abelian varieties $A$ in $\Ical_R$ with $\End(A)=\OO$.
        \item \label{cor:O_min_Bass:Bass} $\OO_{\min}$ is Bass at all maximal ideals.
    \end{enumerate}
\end{corollary}
\begin{proof}
    Let $\OO$ be an overorder of $\OO_{\min}$.
    It follows from \cite[Theorem~5.2]{BergKarMar24_arXiv} that every fractional $\OO$-ideal $I$ with $(I:I)=\OO$ is locally principal.
    It is well-known that this is equivalent to $\OO$ being Gorenstein, see for example \cite[Sec.~3 and 4]{MarType}.
    It follows that $\OO_{\min}$ is Bass if and only if the set of isomorphism classes fractional $\OO_{\min}$-ideals equals
    \[ \bigsqcup_{\OO_{\min} \subseteq \OO} \Cl(\OO). \]
    The result now follows from \cite[Corollary~5.4]{BergKarMar24_arXiv}.
\end{proof}
The statement in Corollary \ref{cor:O_min_Bass} without the assumptions on the isogeny class $\Ical$ does not hold.
Indeed, assume that $\Ical$ is almost ordinary over $\Fq$.
Set $R=\ZZ[\pi,q/\pi]$ and let $\OO_{\min}$ be as in Proposition~\ref{prop:same_num_AV_along_ladder}.\ref{prop:same_num_AV_along_ladder:Omin}.
Assume furthermore that the completion $K_\mm$ of $K$ at the unique maximal ideal $\mm$ containing both $\pi$ and $q/\pi$ is an unramified extension of $\QQ_p$.
Then \cite[Theorem~1.1.(2)]{OswalShankar} implies that for every overorder $\OO$ of $\OO_{\min}$ there will be at least two $\Cl(\OO)$-orbits of abelian varieties $A$ in $\Ical_R$ with $\End(A)=\OO$, even if $\OO_{\min}$ is Bass.
\bigskip

We are now ready to define $\frl$-isogenies.
Ascending and horizontal $\frl$-isogenies will be defined as $\frl$-multiplications satisfying certain extra conditions.
By Lemma~\ref{lem:Imult}, such isogenies are necessarily non-descending.
In the interest of generalizing established results on isogeny graphs, we wish to have a notion of descending isogenies as well.
Conceptually, descending isogenies are defined to be the unique isogeny which completes a certain commutative diagram, which is in turn constructed via the aforementioned ascending $\frl$-multiplications and a similar ideal-theoretic construction which we call a virtual isogeny.
The word virtual indicates that these isogenies will not correspond to edges in the isogeny graphs we will consider in Section~\ref{sec:graph_structure}, see Definition~\ref{def:R_frl_graph}, but rather sit implicitly in the background as a part of the construction.

\begin{definition}[Virtual $\frl$-isogenies]\label{def:virt_isog}
    Let $A \in \Ical_R$ be an abelian variety at level $i$ for $0\leq i \leq d-1$.
    Let $\OO_d\subsetneq \dots \subsetneq \OO_0$ be the $\frl$-multiplicator ladder in the set of overorders of $R$ to which $\End(A) = \OO_i$ belongs.
    The virtual $\frl$-isogeny from $A$ is defined as the $\frl\OO_{i}$-multiplication
    \[ A \to A/A[\frl \OO_{i}]. \]
\end{definition}

\begin{lemma}\label{lem:virt_isog_is_hor}
    Let $A \in \Ical_R$ be an abelian variety at level $i$ for $0\leq i \leq d-1$.
    Assume that $R$ is Bass at $\frl$.
    Then the virtual $\frl$-isogeny $\psi$ from $A$ is horizontal.
\end{lemma}
\begin{proof}
    Let $\OO_d\subsetneq \dots \subsetneq \OO_0$ be the $\frl$-multiplicator ladder in the set of overorders of $R$ for which there is an index $i$ such that $\End(A) = \OO_i$.
    By Theorem~\ref{thm:all_ladders}.\ref{thm:all_ladders:max_idl}, we have $\frL_{i+1}\coloneqq \frl \OO_i$ is the unique singular maximal of $\OO_{i+1}$, which has then multiplicator ring $(\frL_{i+1}:\frL_{i+1}) = \OO_{i}$.
    Since $R$ is Bass at $\frl$, the order $\OO_i$ is Gorenstein at $\frl$.
    Hence, $\frL_{i+1}$ is an invertible $\OO_i$-ideal.
    Therefore, $\psi$ is horizontal by Lemma~\ref{lem:Imult}.\ref{lem:Imult:inv}.
\end{proof}

\begin{remark}\label{rmk:vit_is_hor}
    In Lemma~\ref{lem:virt_isog_is_hor}, we can conclude the virtual isogeny $\psi$ from $A$ at level $i$ for $0 \leq i \leq d-1$ is horizontal in a variety of situations even without assuming that $R$ is Bass at $\frl$.
    Let $\OO_d\subsetneq \dots \subsetneq \OO_0$ be the $\frl$-multiplicator ladder in the set of overorders of $R$ to which $\End(A) = \OO_i$ belongs.

    As we showed in the proof of Lemma~\ref{lem:virt_isog_is_hor}, the ideal $\frl \OO_i$ has multiplicator ring $\OO_{i}$.
    If $\frl \OO_i$ is invertible in $\OO_{i}$ then we can deduce that $\psi$ is horizontal by Lemma~\ref{lem:Imult}.\ref{lem:Imult:inv}.
    This is the case if $i=0$ or $i>1$ because $\OO_i$ is Gorenstein at $\frl$ by Proposition~\ref{prop:ladder_conseq}.
    For $i=1$, if $\OO_1$ is Gorenstein at $\frl$ or if $d\geq 3$ then by Lemma~\ref{lem:not_inv_d_2} again $\frl\OO_1$ is invertible and hence that $\psi$ is horizontal.
    If $\frl\OO_1$ is not invertible in $\OO_1$, like in Example~\ref{ex:frlOO1_not_inv} below, we cannot conclude whether $\psi$ is horizontal or ascending.
\end{remark}

\begin{example}\label{ex:frlOO1_not_inv}
    Consider the isogeny class of abelian threfolds over $\FF_{11}$ defined by the Weil polynomial 
    $ x^6 - 11x^5 + 64x^4 - 255x^3 + 704x^2 - 1331x + 1331 $.
    Consider the order $R\coloneqq \ZZ[\pi,11/\pi]$.
    It has a unique singular maximal ideal $\frl$, which is the unique maximal ideal containing $2$.
    The order $R$ has an $\frl$-multiplicator ladder:
    \[ R = \OO_2 \subsetneq \OO_1 \subsetneq \OO_0 = \OO_K. \]
    The maximal $R$-ideal $\frl=\frl \OO_1$ has multiplicator ring $\OO_1$, but it is not invertible in $R_1$.
\end{example}

\begin{definition}[$\frl$-isogeny]\label{def:l_isogenies}
    Let $A$, $B$ and $C$ be abelian varieties in $\Ical_R$, with $A$ at level $i$, with $0\leq i \leq d$.
    Let $\OO_d\subsetneq \dots \subsetneq \OO_0$ be the $\frl$-multiplicator ladder in the set of overorders of $R$ to which $\End(A) = \OO_i$ belongs.

    We say that an isogeny $\varphi: A \to B$ is an (ascending or horizontal) $\frl$-isogeny from $A$ if 
    $\ker\varphi \cong R/\frl$ as $R$-modules, and
    if $\varphi$ factors as an $\frL$-multiplication followed by an isomorphism, that is, 
    $\varphi: A \to A/A[\frL] \cong B$, where $\frL$ is a maximal ideal of $\OO_i=\End(A)$ above $\frl$.
    
    Assume in addition that $1\leq i \leq d-1$ (so that $\frL=\frl\OO_{i-1}$ by Theorem~\ref{thm:all_ladders}.\ref{thm:all_ladders:max_idl}) and let $\iota: B\to A/A[\frl \OO_{i}]$ be the unique isogeny that completes the diagram
    \[  
        \begin{tikzcd}
                                    & B \arrow[d, dotted, "\iota"] \arrow[dr,"\delta"]\\
           A \arrow[ur, "\varphi"] \arrow[r]   & A/A[\frl \OO_{i}]  \arrow[r,"\cong", no head]& C
        \end{tikzcd}
    \]
    where the horizontal arrow $A\to A/A[\frl \OO_{i}]$ is the virtual $\frl$-isogeny from $A$.
    We say that an isogeny $\delta:B \to C$ is a (descending) $\frl$-isogeny induced by $\varphi$ if 
    $\ker\delta \cong R/\frl$ as $R$-modules, and
    if $\delta$ factors as in the diagram. 
    For $A$ ranging in $\Ical_R$, the ascending and horizontal $\frl$-isogenies from $A$ and the descending $\frl$-isogenies they induce will be called $\frl$-isogenies.
\end{definition}
\begin{remark}\label{rmk:frlisog_levels}
    Let $\varphi:A \to B$ be as in Definition~\ref{def:l_isogenies}.
    Then $\varphi$ is horizontal if $i=0$ and ascending if $i>0$.
    In the latter case, $B$ is at level $i-1$.
    Say now that $1\leq i \leq d-1$ and consider the isogeny $\delta:B\to C$ from Definition~\ref{def:l_isogenies}.
    If $R$ is Bass at $\frl$, or more generally one of the situations described in Remark~\ref{rmk:vit_is_hor} holds, then $C$ is also at level $i$ and $\delta$ is descending.
    These two considerations are independent of the requirements on the kernel of the isogenies.
\end{remark}

\begin{theorem}
\label{thm:pp-isog}
    Assume that $R$ is Bass at a singular maximal ideal $\frl$ coprime to $q$.
    Let $A$ be an abelian variety in $\Ical_R$ at level $i$ and let $\OO_d\subsetneq\dots\subsetneq \OO_0$ be the  $\frl$-multiplicator ladder containing $\End(A)$, so that $\End(A)=\OO_i$.
    Denote by $\frL_1$ the unique maximal ideal of $\OO_1$ above $\frl$.
    The following statements hold.
    \begin{enumerate}[(i)]
        \item \label{it:prop:pp-isog:i_geq_1} If $i > 0$ then there exists a unique ascending $\frl$-isogeny $\varphi$ from $A$.
        \item \label{it:prop:pp-isog:hor} If $0\leq i \leq d-1$ then there exists a unique virtual horizontal isogeny from $A$, which has kernel $R$-linearly isomorphic to $R/\frl\times R/\frl$.
        \item \label{it:prop:pp-isog:desc} If $1\leq i \leq d-1$ then the isogeny $\varphi$ from \ref{it:prop:pp-isog:i_geq_1} induces a unique descending $\frl$-isogeny.
        \item \label{it:prop:pp-isog:i_eq_0} If $i = 0$ then there are no ascending $\frl$-isogenies from $A$ and exactly one of the following possibilities occurs:
        \begin{enumerate}
            \item If $\frL_1$ is inert in $\OO_0$ then there is no horizontal $\frl$-isogeny from $A$.
            \item If $\frL_1$ is split in $\OO_0$ then there are two horizontal $\frl$-isogenies from $A$, both of which are horizontal.
            \item If $\frL_1$ is ramified in $\OO_0$ then there exists a unique horizontal $\frl$-isogeny from $A$.
        \end{enumerate}
    \end{enumerate}
\end{theorem}
\begin{proof}
    For $i=1,\ldots,d$, the unique maximal ideal $\frL_{i}$ of $\OO_{i}$ above $\frl$ is $\frl \OO_{i-1}$ by Theorem~\ref{thm:all_ladders}.\ref{thm:all_ladders:max_idl}.
    Moreover, $\frL_{i}$ has multiplicator ring $\OO_{i-1}$ and residue field $\OO_{i}/\frL_{i} \cong R/\frl$.
    Also, $\OO_{i,\frl} = R_{i,\frl}$, which implies that $\End(A)=\OO_{i}$ is Gorenstein at $\frl$.
    Hence, we deduce \ref{it:prop:pp-isog:i_geq_1} from Proposition~\ref{prop:GorEnd}.
    Assume that $0\leq i \leq d-1$.
    The virtual isogeny from $A$ is horizontal by Lemma~\ref{lem:virt_isog_is_hor}.
    By Proposition~\ref{prop:GorEnd_general}, we have
    \[ A[\frl \OO_{i}] \cong \frac{\OO_i}{\frL_{i+1}}. \]
    As shown in the proof of Proposition~\ref{prop:min_ext}, we have that
    \[ \frac{\OO_i}{\frL_{i+1}} \cong \frac{(\OO_{i+1}/\frL_{i+1})[X]}{(X^2)}\cong \frac{\OO_{i+1}}{\frL_{i+1}} \times \frac{\OO_{i+1}}{\frL_{i+1}}. \]
    This shows that the virtual isogeny from $A$ has kernel $R$-linearly isomorphic to $R/\frl\times R/\frl$, thus completing the proof of \ref{it:prop:pp-isog:hor}.
    Part~\ref{it:prop:pp-isog:desc} follows, by comparing dimensions, from the fact that $\varphi$ has kernel $R$-linearly isomorphic to $R/\frl$ and the virtual horizontal isogeny has kernel $R$-linearly isomorphic to $R/\frl\times R/\frl$.
    Part~\ref{it:prop:pp-isog:i_eq_0} is an application of Theorem~\ref{thm:all_ladders}.\ref{thm:all_ladders:max_idl_OO1} and Proposition~\ref{prop:GorEnd}.
\end{proof}

Recall that $\Cl(\OO)$ acts freely on the set of isomorphism classes of abelian varieties in $\Ical$ with endomorphism ring $\OO$ by Proposition~\ref{prop:Cl_acts}.
Roughly speaking, the next lemma says that ascending $\frl$-isogenies don't mix the orbits.
\begin{lemma}\label{le:class_group_orbits_on_levels}
    Suppose that $R$ is Bass at $\frl$.
    Let $A \in \Ical_R$ be an abelian variety $A$ at level $i$ and let $\OO_d\subsetneq\dots\subsetneq \OO_0$ be the  $\frl$-multiplicator ladder containing $\End(A)$, so that $\End(A)=\OO_i$.
    If $i > 0$ and $B$ is the target of the unique ascending $\frl$-isogeny from $A$, then the following property holds:
    An abelian variety is in the $\Cl(\OO_{i-1})$-orbit of $B$ if and only if it is the target of an ascending $\frl$-isogeny from some abelian variety in the $\Cl(\OO_{i})$-orbit of $A$.
\end{lemma}
\begin{proof} 
    The multiplication of any two ideals corresponds to the composition of the isogenies they induce by \cite[Prop.~3.12]{Waterhouse}.
    An ascending $\frl$-isogeny $A\to B$ is by definition the $\frL$-multiplication from $A$, where $\frL$ is the unique ideal of $\End(A)$ above $\frl$.
    Since the endomorphism algebra is commutative, ideal multiplication is commutative.
    Hence the following diagram, where the vertical arrows denote ascending $\frl$-isogenies and the horizontal arrows denote class group actions (cf.~\ref{prop:Cl_acts}), is commutative.
    \[  
        \begin{tikzcd}
                                    B\cong A/A[\mathfrak L] \arrow[r, "I\OO_{i-1}"]  & B/B[I\OO_{i-1}]\cong A/A[I\mathfrak L] \\
        A \arrow[u] \arrow[r, "I"]   & A/A[I]  \arrow[u]
        \end{tikzcd}
    \]
    Then, the claim of the lemma boils down to show that every ideal class in $\Cl(\OO_{i-1})$ is represented $I\OO_{i-1}$ for some invertible fractional $\OO_i$-ideal $I$.
    This is precisely the well-known surjectivity of the extension map $\Cl(\OO_{i})\to \Cl(\OO_{i - 1})$.
\end{proof}

\section{\texorpdfstring{$(R,\frl)$}{(R,l)}-isogeny graphs}\label{sec:graph_structure}

\begin{figure}
\centering
    \begin{tikzpicture}
        \node[] () at (0,4) {level};
        \node[] () at (2,4) {$\frl$-isogenies};
        \node[] () at (5,4) {singular ideal};
        \node[] () at (8,4) {order};
        \node[] () at (11,4) {conductor};

        \node[] (i-1) at (0,3) {$i-1$};
        \node[] (i) at (0,2) {$i$};
        \node[] () at (0,1) {\vdots};
        \node[] (d) at (0,0) {$d$};

        \node[] (A) at (1.7,2) {$A$};
        \node[] (C) at (2.7,2) {$C$};
        \node[] (B) at (1.7,3) {$B$};
        \draw[->] (A) to node[left] {\tiny{$\frl\mathcal{O}_{i-1}$}} (B);
        \draw[->] (A) to node[above] {\tiny{$\frl\mathcal{O}_{i}$}} (C);

        \node[] (Li-1) at (5,3) {$\frL_{i-1} = \frl\mathcal{O}_{i-2}$};
        \node[] (Li) at (5,2) {$\frL_{i} = \frl\mathcal{O}_{i-1}$};
        \node[] () at (5,1) {\vdots};
        \node[] (Ld) at (5,0) {$\frL_{d} = \frl\mathcal{O}_{d-1}$};
        \node[] () at (5,-1) {\vdots};
        \node[] (l) at (5,-2) {$\frl$};

        \node[] (Oi-1) at (8,3) {$\mathcal{O}_{i-1} = [\frL_i:\frL_i]$};
        \node[] (Oi) at (8,2) {$\mathcal{O}_{i} = [\frL_{i+1}:\frL_{i+1}]$};
        \node[] () at (8,1) {\vdots};
        \node[] (Od) at (8,0) {$\mathcal{O}_d$};
        \node[] () at (8,-1) {\vdots};
        \node[] (R) at (8,-2) {$R$};

        \node[] (ci-1) at (11,3) {$\frl^{i-1}\mathfrak{f}_{\mathcal{O}_d}=\frL^{i-1}\mathfrak{f}_{\mathcal{O}_d}$};
        \node[] (ci) at (11,2) {$\frl^{i}\mathfrak{f}_{\mathcal{O}_d}=\frL^{i}\mathfrak{f}_{\mathcal{O}_d}$};
        \node[] () at (11,1) {\vdots};
        \node[] (cd) at (11,0) {$\mathfrak{f}_{\mathcal{O}_d}$};
    \end{tikzpicture}
\caption{Summary of ideals and orders, for $1<i<d$.}

\label{fig:index_checks}
\end{figure}

Here, we define the $(R,\frl)$-isogeny graph for an isogeny class of abelian varieties. 
In general, an $(R,\frl)$-isogeny graph is disconnected, and we mostly analyze the connected components separately. 
Throughout this section, let $\Ical$ be a squarefree isogeny class of $g$-dimensional abelian varieties over a finite field $\Fq$ with a commutative endomorphism algebra $K = \QQ[\pi]$.
Fix an order $R\subseteq K$ and a maximal ideal $\frl\subset R$ coprime to~$q$ for which $R$ has an $\frl$-multiplicator ladder $R=R_d\subsetneq \ldots \subsetneq R_0$.
As before, we write $\Ical_R$ for the subset of $\Ical$ consisting of abelian varieties in $\Ical$ whose endomorphism rings contain $R$.

\begin{definition}[$(R,\frl)$-isogeny graph]\label{def:R_frl_graph}
The $(R,\frl)$-isogeny graph is defined as follows.
\begin{itemize}
	\item The vertices are $\Fq$-isomorphism classes of abelian varieties in $\Ical_R$.
	\item For every $\frl$-isogeny from $A$ to $B$ for $A,B\in\Ical_R$ up to pre- and post-composition with automorphisms (see Appendix~\ref{app:auts}), we place a directed edge from the vertex of the isomorphism class of $A$ to the vertex of the isomorphism class of $B$.
\end{itemize}
\end{definition}

Recall that, given an abelian variety $B \in \Ical_R$, the virtual horizontal isogeny required for the definition of a descending isogeny to $B$ in Definition~\ref{def:l_isogenies} requires an overorder of $R$ which is strictly contained in $\End(B)$ locally at $\frl$. 
Therefore, if $\End(B)$ is locally equal to $R$ at $\frl$, then there are no descending edges with target $B$.
Sometimes it is possible to find these missing edges by replacing $R$ with a smaller order, see Lemma~\ref{lem:smallerR}, although it is possible that all candidate suborders will fail to have a multiplicator ladder, as needed.

\begin{lemma}\label{lem:smallerR}
    Consider an $(R,\frl)$-isogeny graph $\mathcal{G}$ with $d$ levels.
    Let $R'$ be a suborder of $R$ and define $\frl' = \frl\cap R'$.
    If $R'$ has an $\frl'$-multiplicator ladder, then the $(R,\frl)$-isogeny graph is a subgraph $\mathcal{G}'$ of the $(R',\frl')$-isogeny graph.
    Moreover, for $A,B \in \Ical_R$, if $A\to B$ is an edge $E$ in $\mathcal{G}'$ then $E$ is also an edge in $\mathcal{G}$, unless $E$ is corresponds to a descending $\frl'$-isogeny and $B$ is at level $d$ in $\mathcal{G}$.
 \end{lemma}

\begin{proof} 
    True by definition of the isogeny graph and the fact that $R$ is an $\frl'$-overorder of $R'$.
\end{proof}

\begin{remark}\label{rmk:Cl_graphs_aut}
    Consider an $(R,\frl)$-isogeny graph.
    For every $A \in \Ical_R$ representing a vertex, we have a natural surjective group homomorphism $\Cl(R)\to \Cl(\End(A))$.
    Hence, by Proposition~\ref{prop:Cl_acts}, the association $A\mapsto A/A[I]$ induces a level-preserving action of $\Cl(R)$ on the set of vertices which sends $\frl$-isogeny edges to $\frl$-isogeny edges.
    In other words, each ideal class in $\Cl(R)$ induces a graph automorphism.
\end{remark}

\begin{remark}\label{rmk:findRfromA}
    Fix an abelian variety $A$ in $\Ical$, set $\OO\coloneqq \End(A)$, and pick a maximal ideal $\frL$ of $\OO$.
    Observe that $A$ defines a vertex in some $(R,\frl)$-isogeny graph if and only if $R \subseteq \OO$.
    In this remark we discuss how to choose an order $R$ in $K$ together with a maximal ideal $\frl$ of $R$ such that the $(R,\frl)$-isogeny graph contains the vertex corresponding to $A$, $\frl$ is below $\frL$, $R$ is Bass at $\frl$, and we have descending $\frl$-isogenies with target $A$.
    
    If there exists a pair $(R,\frl)$ that satisfies all these requirements then the following three statements hold:
    every order $S$ such that $R\subseteq S \subseteq \OO$ is Bass at the maximal ideal $\frL \cap S$;
    $\OO$ is not at level $0$ of its $\frl$-multiplicator ladder, that is, $\frL$ is singular; we have a strict inclusion $R\subsetneq \OO$.
    It follows that, if there exists a pair $(R,\frl)$ that satisfies all the requirements then 
    there exists one such that the inclusion $R \subsetneq \OO$ is minimal and $(\frl:\frl)=\OO$.
    These two conditions imply that 
    \[ \frL^2 \subseteq \frl \subset R \subsetneq \OO, \]
    by Proposition~\ref{prop:min_ext}.\ref{prop:min_ext:sing_ram}.
    So, such a `minimal' pair $(R,\frl)$, if it exists, can be found by lifting all the maximal subrings of $\OO/\frL^2$ via the canonical projection $\OO\to \OO/\frL^2$ and check if, among these lifts, there is one, say $R$, which is Bass at the maximal ideal $\frl\coloneqq \frL \cap R$.
\end{remark}

\begin{example}\label{ex:non_Bass}
    Consider the (non-simple) isogeny class \href{https://abvar.lmfdb.xyz/Variety/Abelian/Fq/3/11/b_e_cv}{3.11.b\_e\_cv} of ordinary abelian $3$-folds over $\FF_{11}$.
    The order $R_2\coloneqq \ZZ[\pi,11/\pi]$ has two singular maximal ideals: one, $\frL$, above the rational prime $5$ and one, $\mm$, above $2$.
    Note that $R_2$ is Bass at $\frL$ and not Bass at $\mm$.
    The lattice of $\mm$-overorders consists of six orders one of which is not Gorenstein at~$\mm$.
    In fact, this order, which we denote by $T$ has Cohen-Macaulay type $2$ at the unique maximal ideal above~$\mm$, see Figure~\ref{fig:non_Bass}.

    \begin{figure}[h!]
        \centering
	    \begin{tikzpicture}
            \node[] (R) at (-1,0) {$R_2=\ZZ[\pi,11/\pi]$};
            \node[] (T) at (-1,1) {$T$};
            \node[] (T1) at (-1,2) {$T_2$};
            \node[] (T2) at (-2,2) {$T_1$};
            \node[] (T3) at (0,2) {$T_3$};
            \node[] (Omm) at (-1,3) {$T_0$};
            \node[] (header) at (-1,4) {$\mm$-overorders of $R$};   
           
            \draw[-] (R) -- (T) node[midway,left] {\small{$4$}};
            \draw[-] (T) -- (T1) node[midway,left] {\small{$4$}};
            \draw[-] (T) -- (T2) node[midway,left] {\small{$4$}};
            \draw[-] (T) -- (T3) node[midway,right] {\small{$4$}};
            \draw[-] (T1) -- (Omm) node[midway,left] {\small{$4$}};
            \draw[-] (T2) -- (Omm) node[midway,left] {\small{$4$}};
            \draw[-] (T3) -- (Omm) node[midway,right] {\small{$4$}};

            \node[] (R2) at (4,1) {$R_2=\ZZ[\pi,11/\pi]$};
            \node[] (R1) at (4,2) {$R_1$};
            \node[] (R0) at (4,3) {$R_0$};
            \node[] (header) at (4,4) {$\frL$-multiplicator ladder of $R_2$};
            
            \draw[-] (R2) -- (R1) node[midway,left] {\small{$5$}};
            \draw[-] (R1) -- (R0) node[midway,left] {\small{$5$}};
        \end{tikzpicture}
    \caption{The lattice of inclusions of all overorders of $R_2$ is the cartesian product of the lattice of $\mm$-overorders of $R_2$ and the $\frL$-multiplicator ladder of $R_2$.}
    \label{fig:non_Bass}
    \end{figure}
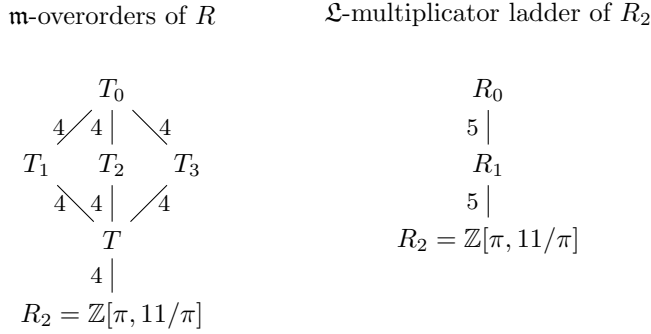

    By Corollary~\ref{cor:number_of_ladders}, there are six $\frL$-ladders in the set of overorders of $R_2$, one per each $\mm$-overorder.
    If we consider an order $S$ in the $\frL$-multiplicator ladder corresponding to $T$, then 
    the constant $N$ from Proposition~\ref{prop:same_num_AV_along_ladder} equals $2$, that is,
    for every order $S$ in such a ladder the number of isomorphism classes with endomorphism ring $S$ is $2\cdot \#Cl(S)$.
    The reason is that these orders $S$ are precisely the non-Gorenstein overorders of $R_2$, having Cohen-Macauly type $2$ and, using the same notation as in the proof of Proposition~\ref{prop:same_num_AV_along_ladder}, hence $w_\mm(S)=2$ by \cite[Thm.~6.2]{MarType}.
    For every other $\frL$-multiplicator ladder, we have $N=1$.

    In order to have descending isogenies to the bottom level, we want to extend the $\frL$-ladder containing $R_2$.
    We proceed following Remark~\ref{rmk:findRfromA}.
    One can verify that the only order $S$ satisfying $\frL^2 \subset S \subset R_2$ is $R_3\coloneqq\ZZ+\frL^2$.
    Moreover, $R_3$ is Bass at the maximal ideal $\frl$ below $\frL$.
    Hence, if we consider the $(R_3,\frl)$-isogeny graph we see that Proposition~\ref{prop:same_num_AV_along_ladder}.\ref{prop:same_num_AV_along_ladder:yesEnds} holds with $d_{\min}=2$ for every $\frl$-multiplicator ladder.
    In fact, since the isogeny class is ordinary, the endomorphism rings are precisely the overorders of $R_2$; cf.~Proposition~\ref{prop:same_num_AV_along_ladder}.
\end{example}

Before we inspect the connected components of $(R,\frl)$-isogeny graphs, we remark that there is subtlety in the definition of connected components for directed graphs. 
In this paper, we simply consider connected components of the underlying undirected graphs.
Nevertheless, in the literature about directed graphs one can find a stronger notion of connectedness: any two vertices $u$ and $v$ are in the same \emph{strongly connected component} if there are directed paths from both $u$ to $v$ and $v$ to $u$.
Although these two definitions of connectedness are inequivalent in general, Lemma~\ref{lem:stronglyconnected} shows that they coincide in our setting, unless we are missing descending edges at the bottom level.

\begin{lemma}\label{lem:stronglyconnected}
    Let $A$ and $B$ be vertices in the $(R,\frl)$-isogeny graph, say at levels $i_A$ and $i_B$ respectively.
    If $i_B<i_A$, assume further that $i_A<d$, where $d$ is the length of the $\frl$-multiplicator ladder in the set of overorders of $R$ containing $\End(A)$.
    If there is a path from $A$ to $B$ then there is also a path from $B$ to $A$.
    In particular, $A$ and $B$ belong to the same strongly connected component.
\end{lemma}
\begin{proof}
    It suffices to consider the case when the path from $A$ to $B$ consists of a single edge.
    Let $\varphi$ be the $\frl$-isogeny inducing the edge.
    We will distinguish three cases.
    Assume first that the edge $\varphi:A\to B$ is horizontal.
    By Remark~\ref{rmk:frlisog_levels}, $A$ and $B$ are at level $0$.
    Let $\OO_0=\End(A)=\End(B)$.
    Then $\varphi$ is induced by an $\frL$-multiplication where $\frL$ is an invertible maximal ideal of $\OO_0$ above $\frl$.
    Let $N$ be the order of the ideal class of $\frL$ in $\Cl(\OO_0)$.
    Then the $\frL^{N-1}$-multiplication from $B$ is a composition of $\frl$-isogeny with target $A$, thus giving a path from $B$ to $A$.

    Assume now that the edge $\varphi:A\to B$ is ascending. 
    Say that $A$ is at level $i$ and $B$ is at level $i-1$.
    Let $\delta: B \to C$ be the descending isogeny induced by $\varphi$, which exists by our assumption on the level $i_A$ of $A$, and let $\psi:A\to C$ be the corresponding virtual horizontal isogeny.
    Now, $\psi$ is an $\frl\OO_i$-multiplication and $\frl\OO_i$ is an invertible $\OO_i$-ideal.
    Let $N$ be the order of its class in $\Cl(\OO_i)$.
    Then the $(\frl\OO_i)^{N-1}$-multiplication from $C$ can be written as sequence of $N-1$ up-down paths in the graph and has target $A$.
    If we pre-compose it with $\delta$ we get a path $B\to A$, as required.

    Finally, assume $\varphi:A\to B$ is descending from $A$ at level $i-1$ to $B$ at level $i$.
    Let $\delta:D\to A$ be the ascending $\frl$-isogeny inducing $\varphi$ and let $\psi:D \to B$ the corresponding virtual horizontal isogeny.
    Then $\psi$ is induced by an $\frl\OO_i$-multiplication.
    Let $N$ be the order of $\frl\OO_i$ in $\Cl(\OO_i)$.
    Then the $(\frl\OO_i)^{N-1}$-multiplication from $B$ can be written as sequence of $N-1$ up-down paths from level $i$ to level $i-1$ in the graph and has target $D$.
    So, after composing with $\delta$, we obtain the required path from $B$ to $A$.  
\end{proof}

\begin{lemma}
\label{le:graph_ladder}
    Consider an $(R,\frl)$-isogeny graph. 
    Assume that $R$ is Bass at $\frl$.
    For every overorder $S$ of $R$, let $\mathcal C(S)$ be the set of all connected components of the $(R,\frl)$-isogeny graph where $S$ occurs as an endomorphism ring of some abelian variety.
    \begin{enumerate}[(a)]
        \item For every $\frl$-multiplicator ladder $\OO_d\subsetneq \dots \subsetneq \OO_0$ in the set of overorders of $R$ which 
        $$
            \mathcal C(\OO_{0}) =\mathcal C(\OO_{1}) =\dots = \mathcal C(\OO_{\dmin}) \neq \emptyset,
        $$
        where $\dmin$ is defined as in Proposition~\ref{prop:same_num_AV_along_ladder}.\ref{prop:same_num_AV_along_ladder:yesEnds}.
        \label{it:component_count} 
        \item If $G$ is any connected component of the $(R,\frl)$-isogeny graph then there is a unique $\frl$-multiplicator ladder $\OO_d\subsetneq \dots \subsetneq \OO_0$ in the set of overorders of $R$ so that $\OO_{d_{\min}},\dots, \OO_0$ are precisely the endomorphism rings of the abelian varieties in $G$, where $\dmin$ is defined as in Proposition~\ref{prop:same_num_AV_along_ladder}.\ref{prop:same_num_AV_along_ladder:yesEnds}. 
        \label{it:component_height}
    \end{enumerate} 
\end{lemma}
\begin{proof} 
    First consider part~\ref{it:component_count}.  
    By Proposition~\ref{prop:same_num_AV_along_ladder}.\ref{prop:same_num_AV_along_ladder:yesEnds}, we know that $\mathcal C(\OO_{\dmin}) \neq \emptyset$, 
    and, since $R$ is Bass at $\frl$, there exists and integer $N$ such that there are $N$ class group orbits of the abelian varieties in $\Ical_R$ with endomorphism ring $\OO_{i}$ for every $\dmin \geq i \geq 1$.
    Now, fix any $i$ satisfying $\dmin \geq i\geq 1$ and let $A_1,\dots, A_N$ be representatives of the $N$ class group orbits.
    For each $1\leq j\leq N$, let $B_j$ be the target of the ascending $\frl$-isogeny from $A_j$.
    By construction, the endomorphism ring of each $B_j$ is $ \OO_{i-1}$, proving $\mathcal C(\OO_{i - 1}) \supseteq \mathcal C(\OO_{i})$.
    Moreover, for every $j$, we see that an abelian variety is in the $\Cl(\OO_{i - 1})$-orbit of $B_j$ if and only if it is the target of an ascending $\frl$-isogeny from some abelian variety in the $\Cl(\OO_{i})$-orbit of $A_j$ by Lemma~\ref{le:class_group_orbits_on_levels}.
    Therefore, the union of orbits $\cup_{j = 1}^N \Cl(\OO_{i-1})\cdot B_j$ is a collection of $N\cdot \#\Cl(\OO_{i - 1})$ abelian varieties with endomorphism ring $\OO_{i-1}$ each one the target of an ascending $\frl$-isogeny from an abelian variety with endomorphism ring $\OO_i$.
    This number accounts for all abelian varieties in $\mathcal I_R$ with endomorphism ring $\OO_{i - 1}$ by Proposition~\ref{prop:same_num_AV_along_ladder}.\ref{prop:same_num_AV_along_ladder:yesEnds}, so we have $\mathcal C(\OO_{i}) = \mathcal C(\OO_{i-1})$ and we are done.
    The existence statement in part~\ref{it:component_height} follows from part~\ref{it:component_count} automatically, while the uniqueness is a consequence of Theorem~\ref{thm:all_ladders}.
\end{proof}

\begin{lemma}\label{lem:down_up}
    Assume that $R$ is Bass at $\frl$.
    Let $G$ be a connected component of the $(R,\frl)$-isogeny graph on $\Ical_R$,
    and let $\OO_d\subsetneq \dots \subsetneq \OO_0$ be the $\frl$-multiplicator ladder in the set of overorders of $R$ which contains precisely the endomorphism rings of the abelian varieties in $G$.
    Define $\dmin$ as in Proposition~\ref{prop:same_num_AV_along_ladder}.\ref{prop:same_num_AV_along_ladder:yesEnds}.
    \begin{enumerate}[(i)]
        \item \label{it:down_up} Assume that $1 \leq i < \dmin$.
        For each down-up path $A \overset{\delta}{\to} B \overset{\varepsilon}{\to} C$ in $G$, with $A,C$ at level $i-1$ and $B$ at level $i$, we have $C\cong A/A[\frl\OO_{i-1}]$.
        \item \label{it:down_up_seq}
        Let $A$ and $C$ be abelian varieties at level $j$ of $G$ with $0 \leq j < \dmin$.
        Assume that there exists a path $A \to B_1 \to \ldots \to B_n \to C$ in $G$
        with $B_1,\ldots,B_n$ all laying on levels $\geq j$.
        Then there exists an integer $m\leq 0$ such that $C\cong A/A[\frl^m\OO_{i-1}]$.
    \end{enumerate}
\end{lemma}
\begin{proof}
    We start by proving part~\ref{it:down_up}.
    Denote by $\delta: A\to B$ the downward arrow in the statements.
    The unique maximal ideal of $\OO_i$ above $\frl$ is $\frl\OO_{i-1}$ by Theorem~\ref{thm:all_ladders}.\ref{thm:all_ladders:ladder}.
    The diagram below depicts abelian varieties at level $i-1$ and $i$, identifying them with their isomorphism classes.
    In particular, we suppress the isomorphisms implicit in the definition of $\frl$-isogenies and virtual $\frl$-isogenies; cf.~Definitions~\ref{def:virt_isog} and \ref{def:l_isogenies}.
    Using this convention we have that $\varepsilon:B\to C$ is the $\frl\OO_{i-1}$-multiplication and $C = B/B[\frl\OO_{i-1}]$. 
    Moreover, the downward $\frl$-isogeny $\delta:A\to B$ comes from an upward edge $D\to A$ given by $\frl\OO_{i-1}$-multiplication, and $D/D[\frl\OO_i]=B$.
    Since we have the equality $(\frl\OO_{i-1})(\frl\OO_{i-1}) = (\frl\OO_{i-1})(\frl\OO_{i})$ and ideal multiplication corresponds to composition of isogenies, we can fill the parallelogram with an isomorphism $C \cong A/A[\frl\OO_{i-1}]$, as desired.
    \[  
        \begin{tikzcd}
            \text{level }i-1 &&      & A \arrow[d, dashed] \arrow[dr, "\delta"]  \arrow[r, "\frl\OO_{i-1}"] & A/A[\frl\OO_{i-1}]  \arrow[r,"\cong", no head] & C\arrow[r, equal]   & B/B[\frl\OO_{i-1}] \\
            \text{level }i  && D \arrow[ur, "\frl\OO_{i-1}"] \arrow[r,swap, "\frl\OO_{i}"]   & D/D[\frl \OO_{i}]  \arrow[r,equal] & B \arrow[urr,swap, "\frl\OO_{i-1}"] \arrow[ur,"\varepsilon"]
        \end{tikzcd}
    \]
    
    We now move to part~\ref{it:down_up_seq}.
    Split the given path connecting $A$ and $C$ at the points where it intersects level $j$, so that each part with the endpoints removed traverses only vertices in levels $>j$.
    It suffices to show the result for each part.
    So, we assume that the abelian varieties $B_1,\ldots,B_n$ occurring in the path from $A$ to $C$ are all at levels $>j$.
    Let $j_{\max}$ be the maximum of these levels.
    We proceed by induction on $m \coloneqq j_{\max} - j$.
    The case $m=1$ is treated by part~\ref{it:down_up}.
    So we assume that the statement holds for all sub-paths touching $<m$ levels.
    By recursively splitting the path where it touches the same level, and applying the inductive hypothesis, we are left to consider the situation where our path is of the form 
    \[ A \overset{\delta}{\longrightarrow} B_1 \overset{(\frl\OO_{j'})^s}{\longrightarrow} B_2 \overset{\varepsilon}{\longrightarrow} C\] 
    with $B_1,B_2$ at level $j'$, $s$ a non-negative integer, $\delta$ a composition of descending $\frl$-isogenies, and $\varepsilon$ a composition of ascending $\frl$-isogenies.
    Set $C' \coloneqq C/C[(\frl\OO_{j})^{-s}]$.
    We summarize our situation in the following diagram, where, as in the proof of part~\ref{it:down_up}, we suppress the isomorphisms implicit in the definition of $\frl$-isogenies and virtual $\frl$-isogenies:
    \[  
        \begin{tikzcd}
            \text{level }j  && A \arrow[dr, "\delta"] \arrow[rr, dotted] & & C' &  & C \arrow[ll, "(\frl\OO_{j})^{-s}"]\\
            \text{level }j' &&   & B_1 \arrow[rr, "(\frl\OO_{j'})^s"] \arrow[ur, dashed]&  & B_2 \arrow[ur, "\varepsilon"] & 
        \end{tikzcd}
    \]
    Composing the arrows $B_1\to B_2 \to C \to C'$ gives the dashed arrow $B_1\to C'$ in the diagram.
    By construction, this is an ascending $\frl$-isogeny.
    By inductive hypothesis, we see that the dotted arrow $A\to C'$ is an $(\frl\OO_j)^{s'}$-multiplication for some non-negative integer $s'$.
    Hence, we obtain that $C\cong A/A[(\frl\OO_j)^{s+s'}]$.
\end{proof}

We are now ready to state and prove the main structure theorem for $(R,\frl)$-isogeny graphs, which was anticipated in Section~\ref{sec:intro} as Main Theorem~\ref{mainthm:B}.
Assume that $R$ is Bass at $\frl$.
Let $G$ be a connected component of the $(R,\frl)$-isogeny graph on $\Ical_R$,
and let $\OO_d\subsetneq \dots \subsetneq \OO_0$ be the $\frl$-multiplicator ladder in the set of overorders of $R$ which contains the endomorphism rings of the abelian varieties in $G$.
Define $\dmin$ as in Proposition~\ref{prop:same_num_AV_along_ladder}.\ref{prop:same_num_AV_along_ladder:yesEnds}.
For each $0\leq i\leq d_{\min}$, let $G_i$ be the subgraph of $G$ whose vertices are isomorphism classes of abelian varieties at level $i$, meaning $\End(A) = \OO_i$.
We identify $A\in \Ical_R$ with the vertex it represents in the graph.
\begin{theorem} 
\label{thm:graph_structure}
    Assume that $R$ is Bass at $\frl$.
    \begin{enumerate}[(a)]
        \item \label{it:surface} 
        Let $\Cl_\frl(\OO_0)$ be  the subgroup of $\Cl(\OO_0)$ generated by the maximal $\OO_0$-ideals lying over $\frl$.
        The vertices of $G_0$ consists of a single orbit of $\Cl_\frl(\OO_0)$.
        \begin{itemize}
            \item If $\delta_\frl =-1$, that is, $\frl$ is inert in $\OO_0$ then $G_0$ is a totally disconnected graph, that is, the set of edges is empty.
            \item If $\delta_\frl =0$ or $1$, that is, $\frl$ is ramified, or split in $\OO_0$, respectively, then $G_0$ is isomorphic to the (directed) Cayley graph of $\Cl_\frl(\OO_0)$ with generators the maximal ideals of $\OO_0$ above $\frl$.
            In particular, $G_0$ is connected.
        \end{itemize}
        
        \item \label{it:edges_up_surface} If $i = 0$ and $A\in G_0$, then there are $\frac{ \#(R/\frl)-\delta_\frl}{[\OO_{0}^\times : \OO_1^\times]}$ vertices in $G_1$ whose (unique) ascending edge has target $A$.
        There are no edges from $G_j$ to $G_0$ for $j > 1$. 
        
        \item \label{it:edges_up_lower} If $1 \leq i < \dmin$ and $A\in G_i$, then there are $\frac{ \#(R/\frl)}{[\OO_{i}^\times : \OO_{i+1}^\times]}$ vertices in $G_{i+1}$ whose (unique) ascending edge has target $A$. 
        These are the only  ascending edges with targets in $G_i$.
                
        \item \label{it:one_desc_per_asc} The in-degree is equal to the out-degree for each vertex on $G_i$ for $1\leq i \leq \min\{\dmin, d-1\}$.  Moreover, if $\dmin = d$, then the vertices at level $d$ have out-degree $1$ and in-degree $0$.
    \end{enumerate}
\end{theorem}
\begin{proof}
    By Lemma~\ref{lem:down_up}, any two vertices of $G_0$ are connected by a composition of $\frl\OO_0$-multiplications, which in turn is a composition of $\frL$-multiplications where $\frL$ is a maximal ideal of $\OO_0$ above $\frl$.
    Since $\Cl(\OO_0)$ acts freely on the set of vertices of $G_0$ by Proposition~\ref{prop:Cl_acts}, we obtain that the set of vertices of $G_0$ consists of a single orbit of $\Cl_\frl(\OO_0)$.
    If $\delta_\frl=-1$ then the unique maximal ideal of $\OO_0$ above $\frl$ is $\frl\OO_0$ which does not induce an $\frl$-isogeny, see Theorem~\ref{thm:pp-isog},\ref{it:prop:pp-isog:i_eq_0}.
    So in this case $G_0$ is totally disconnected.
    If $\delta_\frl=0$ or $1$ and $\frL$ is a maximal ideal of $\OO_0$ above $\frl$ then each $\frL$-multiplication is an $\frl$-isogeny again by Theorem~\ref{thm:pp-isog},\ref{it:prop:pp-isog:i_eq_0}.
    Hence, $G_0$ is isomorphic to the Cayley graph of $\Cl_\frl(\OO_0)$ with generators the maximal ideals of $\OO_0$ above $\frl$.
    This concludes the proof of part~\ref{it:surface}.

    By Lemma~\ref{le:graph_ladder}, there is an $\frl$-multiplicator ladder $\OO_{d}\subsetneq \dots\subsetneq \OO_0$ in the set of overorders of $R$ so that the endomorphism rings of abelian varieties in $G$ are precisely $\{\OO_i \mid i\leq \dmin\}$.
    Let $\mathcal G = \cup_{I \in \Cl(\OO_{\dmin})} [I]\cdot G$ be the graph which is the union of all $\Cl(\OO_{\dmin})$-translates of the component $G$.
    Because every abelian variety in $G$ at level $0$ is in the same class group orbit by part~\ref{it:surface}, the same is true for all abelian varieties in $\mathcal G$ at every level $\dmin\geq i\geq 0$ by Lemma~\ref{le:class_group_orbits_on_levels}.
    Now, recall that by Remark~\ref{rmk:Cl_graphs_aut}, class groups actions induce graph automorphisms of the $(R,\frl)$-isogeny graph.
    In particular, any two vertices of $\mathcal G$ corresponding to abelian varieties $A$ and $A'$ at the same level have the same in-degree and out-degree.
    Hence, since the graph $\mathcal G$ has $\#\Cl(\OO_i)$ vertices at level $i$ for every $0 \leq i\leq \dmin$, the number of ascending edges between two consecutive levels is the ratio of the class numbers.
    This ratio is precisely the value computed in Proposition~\ref{thm:JoPo_mult_ladder}. 
    The nonexistence of other edges follows from Theorem~\ref{thm:pp-isog}.
    This completes the proofs of parts \ref{it:edges_up_surface} and \ref{it:edges_up_lower}.  

    Now, note that each descending edge is uniquely defined by an ascending edge and a unique virtual horizontal isogeny by  Theorem~\ref{thm:pp-isog}.
    Suppose two ascending isogenies $A_1 \overset{\varphi_1}{\to} B$, $A_2 \overset{\varphi_2}{\to} B$ from level $i$ to level $i-1$ do define the same descending isogeny $B\to C$. Then $C\cong A_1/A_1[\frl\mathcal{O}_i]\cong A_2/A_2[\frl\mathcal{O}_i]$ by definition. Since $\frl\mathcal{O}_i$ is invertible in $\mathcal{O}_i$, this gives $A_1\cong A_2$ and thus $\varphi_1 = \varphi_2$.
    This proves part~\ref{it:one_desc_per_asc}.
\end{proof}

Theorem~\ref{thm:graph_structure} shows that it is possible for a connected component of an $(R,\frl)$-isogeny graph to have a totally disconnected surface. We provide an explicit example below.

\begin{example}\label{ex:c_ex_BJW}
    Consider the absolutely simple ordinary isogeny class with label \href{https://www.lmfdb.org/Variety/Abelian/Fq/4/5/e_f_ax_adi}{4.5.e\_f\_ax\_adi} of abelian $4$-folds defined over $\FF_5$.
    Let $K\coloneqq \QQ[\pi]$ be the endomorphism algebra of the isogeny class.
    The order $R_2\coloneqq \ZZ[\pi,5/\pi]$ is Bass and has a unique singular maximal ideal $\frL$ above $2$.
    The procedure described in Remark~\ref{rmk:findRfromA} finds a single maximal underorder $R$ of $\ZZ[\pi,5/\pi]$ which is Bass at the unique maximal ideal $\frl$ below $\frL$.
    The set of overorders of $R$ fits into the $\frl$-multiplicator ladder of $R$, which is
    \[ R \subseteq R_2 \subseteq R_1 \subseteq \OO_K. \]
    So, $d=3$ and $\dmin=2$.
    For later use, we remark that $\OO_{K^+}$ is contained in $R$.
    So the same holds for every order in the ladder.
    We have that $\delta_\frl=-1$, that is, the ideal $\frl$ is inert in $\OO_K$.
    Also, $\frl\OO_K$ has order $2$ in $\Cl(\OO_K)$.
    The $(R,\frl)$-isogeny graph has two isomorphic connected components which are depicted in Figure~\ref{fig:c_ex_BJW}
    \begin{figure}[htbp]
    \centering
        \begin{tikzpicture}[scale=1, vertex/.style={circle, draw, minimum size=6pt, inner sep=2pt}]
            \node[] () at (-3,0) {level $0$};
            \node[] () at (-3,-1.5) {level $1$};
            \node[] () at (-3,-3) {level $2$};
            \node[vertex] (A1) at (-1,0) {};
            \node[vertex] (A2) at (1,0) {};
            \node[vertex] (B1) at (-1,-1.5) {};
            \node[vertex] (B2) at (1,-1.5) {};
            \node[vertex] (C1) at (-1,-3) {};
            \node[vertex] (C2) at (1,-3) {};
            \draw[->] (A1) -- (B2);
            \draw[->] (A2) -- (B1);
            \draw[->] (B2) -- (C1);
            \draw[->] (B1) -- (C2);
            \draw[->] (C1) -- (B1);
            \draw[->] (C2) -- (B2);
            \draw[->] (B1) -- (A1);
            \draw[->] (B2) -- (A2);
        \end{tikzpicture}
        \caption{A connected component of the $(R,\frl)$-isogeny graph described in Example~\ref{ex:c_ex_BJW} for the isogeny class \href{https://www.lmfdb.org/Variety/Abelian/Fq/4/5/e_f_ax_adi}{4.5.e\_f\_ax\_adi}.}
        \label{fig:c_ex_BJW}
    \end{figure}
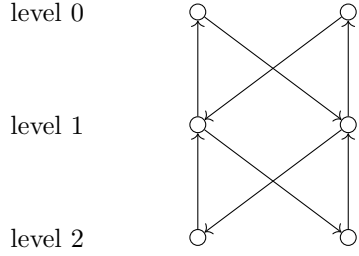
\end{example}

We conclude this section by deriving a formula for the number of connected componenets of an $(R,\frl)$-isogeny graph.
\begin{proposition}\label{prop:num_conn_comps}
    Consider an $(R,\frl)$-isogeny graph.
    Assume that $R$ is Bass at $\frl$.
	Write $\textsc{Max}(R,\frl)$ for the set of $\frl$-maximal overorders of $R$ belonging to an $\frl$-multiplicator ladder in the set of overorders of $R$.
    Given $\OO_0\in \textsc{Max}(R,\frl)$, let $\Cl_\frl(\OO_0)$ be the subgroup of $\Cl(\OO_0)$ generated by maximal ideals of $\OO_0$ over $\frl$.
    If $\OO_0$ is not an endomorphism ring of some abelian variety in $\Ical_R$ then set $N_{\OO_0}\coloneqq 0$; otherwise let $N_{\OO_0}$ be number of $\Cl(\OO_0)$-orbits of abelian varieties in $\Ical_R$ with endomorphism ring $\OO_0$.
	Then the number of connected components in the $(R,\frl)$-isogeny graph on $\mathcal{I}_R$ is 
	\[
		\sum_{\OO_0\in \textsc{Max}(R,\frl)} \left( N_{\OO_0}\cdot[\Cl(\OO_0) : \Cl_\frl(\OO_0)]\right).
	\]
\end{proposition}
\begin{proof}
  	An order $\OO_0$ occurs as the endomorphism ring of an abelian variety at the surface of a connected component if and only if $\OO_0\in \textsc{Max}(R,\frl)$ by Lemma~\ref{le:graph_ladder} and Theorem~\ref{thm:graph_structure}.
	Moreover, the same theorem shows that every abelian variety $A$ representing a vertex at the surface $G_0$ of a connected component $G$ has the same endomorphism ring $\End(A) = \OO_0$ for some $\OO_0\in \textsc{Max}(R,\frl)$, and we have $\#G_0 = \#\Cl_\frl(\OO_0)$.
    Also, by Proposition~\ref{prop:same_num_AV_along_ladder}, there are $N_{\OO_0}\cdot\#\Cl(\OO_0)$ abelian varieties in~$\mathcal{I}_R$ with endomorphism ring $\OO_0$.
    Therefore, for each fixed $\OO_0\in \textsc{Max}(R,\frl)$,  there are $N_{\OO_0}\cdot[\Cl(\OO_0) : \Cl_\frl(\OO_0)]$ connected components such that the endomorphism ring appearing at the surface is $\OO_0$. 
    The claim now follows from taking the sum over $ \textsc{Max}(R,\frl)$.
\end{proof}

\begin{remark}
    Maintain the notation and hypotheses of Proposition~\ref{prop:num_conn_comps}.
    Assume furthermore that $\Ical$ is ordinary, almost ordinary, or over the prime field $\FF_p$, and consider an $(R,\frl)$-isogeny graph.
    Let $\OO_{\min}$ be defined as in Proposition~\ref{prop:same_num_AV_along_ladder}.\ref{prop:same_num_AV_along_ladder:Omin}.
    Assume that $\OO_{\min}$ is Bass but not maximal.
    Then $\OO_{\min}$ has a multiplicator ladder at every singular prime by Proposition~\ref{prop:Bass_tot_ord}. 
    Let $\mm_1,\dots, \mm_n$ be the singular primes of $\OO_{\min}$ and write $d_i$ for the length of the $\mm_i$-multiplicator ladder for each $i$.
    Assume without loss of generality that $\mm_1$ is the unique maximal ideal or $\OO_{\min}$ above $\frl$.
    Then the set $\textsc{Max}(R,\frl)$ has $\prod_{i  \neq 1}(d_i - 1)$ elements described by Corollary~\ref{cor:number_of_ladders}.
    Moreover, if $\Ical$ is ordinary or defined over $\FF_p$ then for each $\OO_0$ in $\textsc{Max}(R,\frl)$ such that $N_{\OO_0}\neq 0$ we have $N_{\OO_0}=1$.
    If $\Ical$ is almost ordinary then for each $\OO_0$ in $\textsc{Max}(R,\frl)$ such that $N_{\OO_0}\neq 0$ we have $N_{\OO_0}=1$ if the unique maximal ideal $\pp$ of $\ZZ[\pi,q/\pi]$ containing both $\pi$ and $q/\pi$ is ramified in $O_K$ and $N_{\OO_0}=2$ otherwise, see~\cite[Theorem~1.1]{OswalShankar} and \cite[Rmk.~2.11]{BergKarMar}.
\end{remark}

\section{Volcanoes}\label{sec:volcanoes}
In this section we relate the $(R,\frl)$-isogeny graphs th structure of $r$-volcano defined in Definition~\ref{def:r-volcano}. 
Volcanoes are undirected graphs, whereas the $(R,\frl)$-isogeny graphs are directed. 
This issue has a natural solution for elliptic curves, as an $\ell$-isogeny $\varphi:E_1\to E_2$ has a unique dual $\ell$-isogeny $\widehat{\varphi}:E_2\to E_1$.
This is because of the fact that an elliptic curve is isomorphic to its dual, which is not true in general for higher dimensional abelian varieties.
By Definition~\ref{def:l_isogenies}, ascending $\frl$-isogenies uniquely determine descending $\frl$-isogenies, but not necessarily between the same pair of vertices, as we have seen in Example~\ref{ex:c_ex_BJW}.
To address this, we define separately undirected ascending and descending graphs in Definition~\ref{def:asc_desc_G}. 
We say that the connected component is a volcano when each of these graphs is a volcano and when they coincide, see Definition~\ref{def:G_volcano}.

Throughout this section, let $\Ical$ be a squarefree isogeny class of $g$-dimensional abelian varieties over a finite field $\Fq$ with a commutative endomorphism algebra $K = \QQ[\pi]$.
Fix an order $R\subseteq K$ and a maximal ideal $\frl\subset R$ coprime to~$q$.
Assume that $R$ is Bass at $\frl$.
Fix a connected component $G$ of the $(R,\frl)$-isogeny graph.
Denote by $\OO_d\subsetneq \dots \subsetneq \OO_0$ the $\frl$-multiplicator ladder in the set of overorders of $R$ containing all the endomorphism rings of the abelian varieties in $G$, see Lemma~\ref{le:graph_ladder}.\ref{it:component_height}.
Define $\dmin$ as in Proposition~\ref{prop:same_num_AV_along_ladder}.\ref{prop:same_num_AV_along_ladder:yesEnds}, so that each $\OO_i$ with $0\leq i\leq \dmin$ is an endomorphism ring.
For each maximal ideal $\frL$ of $\OO_0$ above $\frl$ we denote by $\ord(\frL)$ the order of the isomorphism class of $\frL$ in the class group $\Cl(\OO_0)$.

\begin{definition}[Ascending and descending graphs]\label{def:asc_desc_G}
    Define the \emph{undirected ascending graph} $\Gasc$ (resp.~\emph{undirected descending graph} $\Gdesc$) associated to $G$ as the graph with the same vertices as $G$ and whose edges are ascending (resp.~descending) edges of $G$ without their direction, and whose horizontal edges are determined as follows: draw $k$ undirected edges connecting a pair of vertices $(v_1,v_2)$ whenever the directed graph contains both $k$ horizontal edges of the form $(v_1,v_2)$ and $k$ horizontal edges of the form $(v_2,v_1)$ for $v_1\neq v_2$, and if a vertex has a loop then draw this loop as undirected (it may only be traversed in one direction). 
\end{definition}

\begin{definition}[$G$ is an $r$-volcano]\label{def:G_volcano}
    Let $r$ be a positive integer. 
    We say that $G$ is an $r$-volcano if $\Gasc=\Gdesc$ and $\Gasc$ is an $r$-volcano.
\end{definition}

We give a characterization of when $G$ is an $r$-volcano in Theorems~\ref{thm:volcano_inertramified} and \ref{thm:volcano_split} below, thus proving Main Theorem~\ref{mainthm:C} from Section~\ref{sec:intro}.

While $G$ is a connected graph, it is possible for $\Gasc$ and $\Gdesc$ to be disconnected, see the $(R,\frl)$-isogeny graph in Example~\ref{ex:c_ex_BJW}.

\begin{lemma}\label{lem:Gasc_conn_iff_Gasc0}
    The ascending graph $\Gasc$ is connected if and only if the surface $\Gasc_0$ is.
\end{lemma}
\begin{proof}
    Every vertex at level $i>0$ of $\Gasc$ is connected to a unique vertex at level $i-1$.
    It follows that the vertices at level $i$ are in the same connected components if and only if the vertices above them on level $i-1$ are in the same connected component.
    The statement follows readily.
\end{proof}

\begin{lemma}\label{lem:Gasc_eq_Gdesc}
    If $\dmin=0$ then the graphs $\Gasc$ and $\Gdesc$ are equal.
    If $\dmin \geq 1$ then $\Gasc$ and $\Gdesc$ are equal if and only if $\dmin < d$ and $\frl\OO_{\dmin-1}$ is principal as an $\OO_{\dmin-1}$-ideal.
\end{lemma}
\begin{proof}
    By definition, $\Gasc$ and $\Gdesc$ are automatically equal at the surface, so we only need to consider edges outside of the surface. 
    Thus, we assume $\dmin > 0$.   If $\dmin = d$, then there are no descending edges with targets at level $\dmin$, immediately proving $\Gasc$ and $\Gdesc$ are not equal. For the rest of the proof we assume that $0 < \dmin < d$.
    
    By Theorem~\ref{thm:graph_structure}, all edges outside the surface occur between consecutive levels.
    Moreover, in each of $\Gdesc$ and $\Gasc$, every abelian variety $B$ at level $i$ is connected to a unique vertex at level $i-1$ for every $1\leq i\leq \dmin$. 
    Therefore, $\Gdesc$ and $\Gasc$ are equal if and only if the unique vertex above $B$ in $\Gdesc$ is the same as the unique vertex above $B$ in $\Gasc$ for all abelian varieties $B$ below the surface.
    By Lemma~\ref{lem:down_up}, if $B$ is an abelian variety at level $i$ for $1\leq i\leq \dmin$ and $A$ is the unique vertex above $B$ in $\Gdesc$, then we identify the unique vertex above $B$ in $\Gasc$ to be $A/A[\frl\OO_{i -1}]$.
    Thus, $\Gdesc$ is equal to $\Gasc$ if and only if $A$ is isomorphic to $A/A[\frl\OO_{i -1}]$ for every abelian variety $A$ at level $i-1$ for every $1\leq i\leq \dmin$.
    Furthermore, the freeness of the class group action implies  $A\cong A/A[\frl\OO_{i-1}]$ if and only if  $\frl\OO_{i-1}$ is a principal $\OO_{i-1}$-ideal.
   Finally, $\frl\OO_{i-1}$ is a principal $\OO_{i-1}$-ideal for all $1\leq i\leq \dmin$ if and only if $\frl\OO_{\dmin-1}$ is a principal $\OO_{\dmin-1}$-ideal. 
   This concludes the proof.
\end{proof}

\begin{proposition}
    Let $r$ be a positive integer.
    If $\dmin=0$ or $d>\dmin$ then $\Gasc$ is an $r$-volcano if and only if $\Gdesc$ is an $r$-volcano.
\end{proposition}
\begin{proof}
    The case $\dmin=0$ follows directly from Lemma~\ref{lem:Gasc_eq_Gdesc}.
    Assume now that $d>\dmin>0$.
    By definition, there are three properties to check to establish whether or not $\Gasc$ and $\Gdesc$ are $r$-volcanoes, see Definition~\ref{def:r-volcano}.  
    We show that each of these three properties holds for $\Gasc$ if and only if it holds for $\Gdesc$.
    For \eqref{it:r-volcano-1}, we observe that the graphs $\Gasc$ and $\Gdesc$ have the same vertex set and the same edges at level 0 by definition.
    Property \eqref{it:r-volcano-2} holds for every $(R,\frl)$-isogeny graph automatically.
    For \eqref{it:r-volcano-3}, note that there is a one-to-one correspondence between the number of ascending edges from level $i$ to level $i-1$ and the number of descending edges from level $i-1$ to level $i$ for each $1\leq i\leq \dmin$ because $\dmin < d$.
    Because in-degree is equal to out-degree for every vertex at level $i$ and every $i > \dmin$ (Theorem~\ref{thm:graph_structure}.\ref{it:one_desc_per_asc}), every vertex of $\Gasc$ not at level $\dmin$ has degree $r +1$ if and only if the same is true for $\Gdesc$, completing the proof.
\end{proof}

\begin{lemma}\label{lem:r_volc_minus_r0}
    Assume $\dmin>0$.
    Then,
    $\Gasc$ is an $r$-volcano for a positive integer $r$ if and only if the following conditions hold:
    \begin{enumerate}[(i)]
        \item $\Gasc_0$ is a connected $r_0$-regular graph with $r_0\leq 2$;
    	\item\label{lem:gascvolcano_it1} $r +1 =\frac{\#(R/\frl)-\delta_\frl}{[\OO_0^\times:\OO_1^\times]}+r_0$;
	    \item\label{lem:gascvolcano_it2} $ r + 1 = \frac{\#(R/\frl)}{[\Ocal_i^\times : \Ocal_{i+1}^\times]} + 1$ for $1 \leq i < \dmin$.
    \end{enumerate}
\end{lemma}
\begin{proof}
    The result follows from combining Lemma~\ref{lem:Gasc_conn_iff_Gasc0} with Theorem~\ref{thm:graph_structure}.\ref{it:edges_up_surface} and \ref{it:edges_up_lower}.
\end{proof}

We study the surface $\Gasc_0$ of $\Gasc$.
Its structure depends on the splitting behaviour of $\frl\mathcal{O}_0$. 
First, we address the inert and ramified cases together, followed by the split case.
\begin{lemma}\label{lem:Gasc0conn_inert_ram}
    Assume that $\delta_\frl = -1$ or $0$.
    Then the following are equivalent:
    \begin{enumerate}[(i)]
        \item $\Gasc_0$ is connected.
        \item $\frl\OO_0$ is principal.
    \end{enumerate}
    If this is the case then $\Gasc$ is a $r_0$-regular graph with $r_0=\delta_\frl+1$.
\end{lemma}
\begin{proof}
    If $\delta_\frl=-1$ then $\Gasc_0$ is consists of $\ord(\frl\OO_0)$ vertices with no edges, by Theorem~\ref{thm:graph_structure}.\ref{it:surface}. 
    So it is connected if and only if $\frl\OO_0$ is principal.
    Assume now that $\delta_\frl=0$ and write $\frl\OO_0=\frL^2$ for the maximal ideal $\frL$ of $\OO_0$ above $\frl$.
    If $\frL$ is principal then $\Gasc_0$ is a single vertex with a loop.
    If $\ord(\frL)=2$ then $\Gasc_0$ consists of two vertices connected by one edge.
    If $\ord(\frL)=n>2$ then $G_0$ is a directed $n$-cycle, which implies that $\Gasc_0$ has no edges. 
    This concludes the proof.
\end{proof}

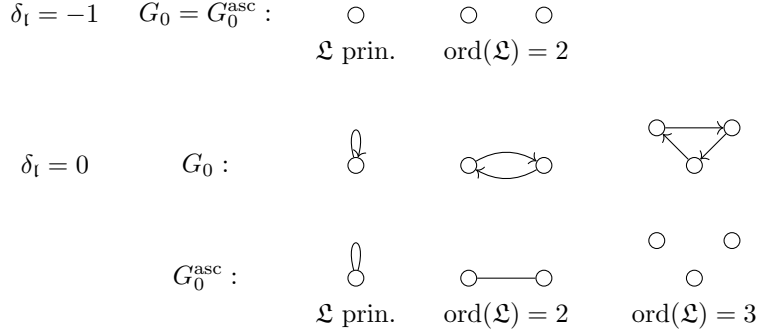
\begin{figure}[h!]
    \centering
    \begin{tikzpicture}[scale=1, vertex/.style={circle, draw, minimum size=6pt, inner sep=2pt}]
        \node[] (d-1) at (-2,8) {$\delta_\frl=-1$};
        \node[vertex] (d-10) at (2,8) {};
        \node[] () at (2,7.5) {$\frL$ prin.};
        \node[vertex] (d-11) at (3.5,8) {};
        \node[vertex] (d-12) at (4.5,8) {};
        \node[] () at (4,7.5) {$\ord(\frL)=2$};
        \node[] () at (0,8) {$G_0 = \Gasc_0:$};

        \node[] (d0) at (-2,6) {$\delta_\frl=0$}; 
        \node[vertex] (d00) at (2,6) {};
        \node[vertex] (d01) at (3.5,6) {};
        \node[vertex] (d02) at (4.5,6) {};
        \node[vertex] (d03) at (6,6.5) {};
        \node[vertex] (d04) at (7,6.5) {};
        \node[vertex] (d05) at (6.5,6) {};
        \draw[->,loop above] (d00) to (d00);
        \draw[->] (d01) to [out=30,in=150] (d02);
        \draw[->] (d02) to [out=210,in=330] (d01);
        \draw[->] (d03) to (d04);
        \draw[->] (d04) to (d05);
        \draw[->] (d05) to (d03);

        \node[vertex] (d00a) at (2,4.5) {};
        \node[vertex] (d01a) at (3.5,4.5) {};
        \node[vertex] (d02a) at (4.5,4.5) {};
        \node[vertex] (d03a) at (6,5) {};
        \node[vertex] (d04a) at (7,5) {};
        \node[vertex] (d05a) at (6.5,4.5) {};
        \draw[-,loop above] (d00a) to (d00a);
        \draw[-] (d01a) to (d02a);
        \node[] (G0) at (0,6) {$G_0:$};
        \node[] (Gasc) at (0,4.5) {$\Gasc_0:$};
        \node[] () at (2,4) {$\frL$ prin.};
        \node[] () at (4,4) {$\ord(\frL)=2$};
        \node[] () at (6.5,4) {$\ord(\frL)=3$};
    \end{tikzpicture}
    \caption{Example diagrams of the structures of $G_0$ and $\Gasc_0$ for $\delta_\frl=-1$ and $0$.}
    \label{fig:Level0s_inert_ram}
\end{figure}

\begin{theorem}\label{thm:volcano_inertramified}
    Assume that $\delta_\frl=-1$ or $0$. 
    Let $r$ be a positive integer.
    \begin{enumerate}[(i)]
        \item If $\dmin=0$ then $G$ is an $r$-volcano if and only if $\frl\OO_0$ is principal.
        \item \label{thm:volcano_inertramified:ii} If $\dmin>0$ then $G$ is an $r$-volcano if and only if the following three conditions hold:
        \begin{itemize}
            \item $\frl\OO_{\dmin-1}$ is principal;
            \item\label{lem:gascvolcano_it1} $r +1 =\frac{\#(R/\frl)-\delta_\frl}{[\OO_0^\times:\OO_1^\times]}+r_0$, where $r_0=\delta_\frl+1$;
            \item\label{lem:gascvolcano_it2} $ r + 1 = \frac{\#(R/\frl)}{[\Ocal_i^\times : \Ocal_{i+1}^\times]} + 1$ for $1 \leq i < \dmin$.
        \end{itemize}
    \end{enumerate}
\end{theorem}
\begin{proof}
	If $\dmin=0$, Lemma~\ref{lem:Gasc_eq_Gdesc} tells us that $\Gasc$ and $\Gdesc$ are equal. Further in this case, Lemma~\ref{lem:Gasc0conn_inert_ram} gives us a connected $r_0$-regular graph if and only if $\frl\OO_0$ is principal.
	If $\dmin>0$, Lemma~\ref{lem:Gasc_eq_Gdesc} tells us that $\Gasc$ and $\Gdesc$ are equal if and only if $\dmin<d$ and $\frl\OO_{\dmin-1}$ is principal. If $\frl\OO_{\dmin-1}$ is principal then $\frl\OO_0$ is principal as well. By Lemma~\ref{lem:r_volc_minus_r0}, $\Gasc$ is an $r$-volcano (for the specified $r$, in particular) if and only if $\Gasc_0$ is a connected $r_0$-regular graph, which is holds when $\frl\OO_0$ is principal by Lemma~\ref{lem:Gasc0conn_inert_ram}.
\end{proof}

To motivate the case where $\frl\OO_0$ is split (i.e., $\delta_\frl=1$.), we begin with an example to showcase how the surface structure can impact a potential volcano theorem.
A statement analogous to Lemma~\ref{lem:Gasc0conn_inert_ram} does not hold in this case, as Example~\ref{ex:Gasc0_split_frl_nonprinc} shows.

\begin{example}\label{ex:Gasc0_split_frl_nonprinc}
    Consider the isogeny class \href{https://www.lmfdb.org/Variety/Abelian/Fq/6/2/b_e_d_l_l_be}{6.2.b\_e\_d\_l\_l\_be} corresponding to the Weil polynomial
    $$
        h(x) = ( x^2 - 2 x + 2 )( x^{10} + 3 x^9 + 8 x^{8} + 13 x^{7} + 21 x^{6} + 27 x^{5} + 42 x^{4} + 52 x^{3} + 64 x^{2} + 48 x + 32 ).
    $$
    There are ten isomorphism classes in the isogeny class, which is neither ordinary nor supersingular and includes the product of a supersingular elliptic curve and an ordinary abelian $5$-fold over the finite field $\FF_2$.
    The order $\OO_1=\ZZ[\pi,2/\pi]$ has two singular maximal ideals, of the form $\frl_2$ and $\bar\frl_2$, both containing $5$.
    Since the maximal ideals are conjugate, we will discuss the situation only for one of them.
    The order $\OO_1$ is globally Bass and has a unique maximal suborder $R$ which is Bass at $\frl\coloneqq \frl_2$.
    One computes that the set of overorders of $\OO_1$ splits into to $\frl$-multiplicator ladders:
    \[ \OO_1 \subsetneq  \OO_0\quad\text{and}\quad
    \OO'_1 \subsetneq  \OO'_0. \]
    With this notation, $\OO'_0$ is the maximal order of $\QQ[\pi]$.
    We see that $\dmin=1$.
    The ideal $\frl$ splits into $\OO_0$ (and $\OO'_0$), that is, $\delta_\frl=1$.
    The $(R,\frl)$-isogeny graph consists of two connected components, one per ladder.
    We see that $\frl\OO_0 = \frL_1\frL_2$ with $\ord\frL_1 = 1$ and $\ord\frL_2 = 2$. 
    In particular, $\frl\OO_0$ is not principal, even if the surface $\Gasc_0$ of the ascending graph $\Gasc$ associated to each component $G$ is connected.
    See Figure~\ref{fig:6.2isograph}.
    
    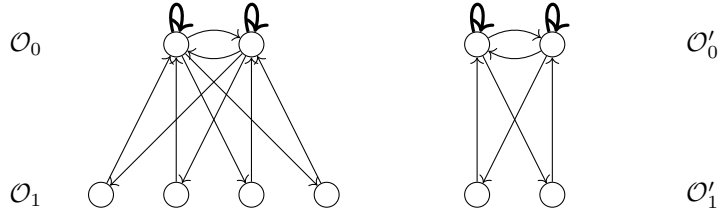
\begin{figure}[h!]
    \centering
    \begin{tikzpicture}
    	\node[circle,draw] (1) at (0,0) {};
    	\node[circle,draw] (2) at (1,0) {};
        \node[] () at (-2,0) {$\OO_0$}; 
        \draw[loop above,->,very thick] (1) to (1);
        \draw[loop above,->,very thick] (2) to (2);
        \draw[->] (1) to [out=30,in=150] (2);
        \draw[->] (2) to [out=210,in = -30] (1);
        \node[] () at (-2,-2) {$\OO_1$};
        \node[circle,draw] (3) at (-1,-2) {};
        \node[circle,draw] (4) at (0,-2) {};
        \node[circle,draw] (5) at (1,-2) {};
        \node[circle,draw] (6) at (2,-2) {};
        \draw[->] (3) to (1);
        \draw[->] (2) to (3);
        \draw[->] (4) to (1);
        \draw[->] (2) to (4);
        \draw[->] (5) to (2);
        \draw[->] (1) to (5);
        \draw[->] (6) to (2);
        \draw[->] (1) to (6);
        
        \node[] () at (7,0) {$\OO'_0$};
        \node[circle,draw] (1) at (4,0) {};
        \node[circle,draw] (2) at (5,0) {};
        \draw[loop above,->,very thick] (1) to (1);
        \draw[loop above,->,very thick] (2) to (2);
        \draw[->] (1) to [out=30,in=150] (2);
        \draw[->] (2) to [out=210,in = -30] (1);
        \node[] () at (7,-2) {$\OO'_1$};
        \node[circle, draw] (7) at (4,-2) {};
        \node[circle, draw] (8) at (5,-2) {};
        \draw[->] (7) to (1);
        \draw[->] (2) to (7);
        \draw[->] (8) to (2);
        \draw[->] (1) to (8);
    \end{tikzpicture}
    \caption{Connected components of the $(R,\frl)$-isogeny graph for the isogeny class \href{https://www.lmfdb.org/Variety/Abelian/Fq/6/2/b_e_d_l_l_be}{6.2.b\_e\_d\_l\_l\_be}. For each component, the surface vertices are depicted at the top, and the actions of $\frL_1$ and $\frL_2$ are distinguished by making $\frL_1$ bold and $\frL_2$ not bold.}
    \label{fig:6.2isograph}
\end{figure}
\end{example}

Producing a characterization of when $\Gasc_0$ is connected in the split case requires an analysis of the plethora of cases determined by how the subgroups in $\Cl_\frl(\OO_0)$ generated by the maximal $\frL_1$ and $\frL_2$ of $\OO_0$ above $\frl$ intersect.
In particular, it is possible for the orbit of $\frL_1$ to be contained in the orbit of $\frL_2$, or they could only partially overlap and in a variety of ways.
For this reason, we characterize when $G$ is a volcano for $\dmin = 0$ only under the additional assumption that $\frl\OO_0$ is a principal $\OO_0$-ideal, that is, when $\frL_1\cong\frL_2^{-1}$. In particular, for $\dmin = 0$ we do not obtain an `if and only if' characterization of volcano graphs.
If $\dmin >0$, we do not lose any of the generality of Theorem~\ref{thm:volcano_inertramified}, by Lemma~\ref{lem:Gasc_eq_Gdesc}.
\begin{lemma}\label{lem:Gasc0conn_split} 
    Assume that $\delta_\frl = 1$ and that $\frl\OO_0$ is principal.
    Then 
    $\Gasc_0$ is a connected regular graph with regularity $r_0=2$.
\end{lemma}
\begin{proof}
    Let $\frL_1$ and $\frL_2$ be the maximal ideal of $\OO_0$ above $\frl$.
    Since $\frl\OO_0$ is principal, we get that $\frL_1\cong\frL_2^{-1}$.
    If $\ord(\frL_1)=1$ then $\Gasc_0$ consists of a single vertex with $2$ loops.
    If $\ord(\frL_1)=n>1$ then $G_0$ consists of $n$ vertices giving rise to two directed $n$-cycles going with arrows going in opposite directions.
    Hence $\Gasc_0$ is an undirected $n$-cycle.
\end{proof}

\begin{figure}[h!]
    \centering
    \begin{tikzpicture}[scale=1, vertex/.style={circle, draw, minimum size=6pt, inner sep=2pt}]
        \node[] (d1) at (-2,2.5) {$\delta_\frl=1$};
        \node[vertex] (d10) at (2,2.5) {};
        \draw[->, loop above] (d10) to (d10); 
        \draw[->, loop left] (d10) to (d10); 

        \node[vertex] (d10a) at (2,0.5) {};
        \draw[-, loop above] (d10a) to (d10a); 
        \draw[-, loop left] (d10a) to (d10a); 

        \node[] () at (2,0) {$\frL_1,\frL_2$ prin.};

        \node[vertex] (d11) at (5,2.5) {};
        \node[vertex] (d12) at (6,2.5) {}; 
        \draw[->] [out=30,in=150] (d11) to (d12);
        \draw[->] [out=60,in=120] (d11) to (d12);
        \draw[->] [out=210,in=330] (d12) to (d11);
        \draw[->] [out=240,in=300] (d12) to (d11);

        \node[vertex] (d11a) at (5,0.5) {};
        \node[vertex] (d12a) at (6,0.5) {}; 
        \draw[-] [out=60,in=120] (d11a) to (d12a);
        \draw[-] [out=210,in=330] (d12a) to (d11a);

        \node[] () at (5.5,0) {$\ord\frL_i=2$};

        \node[vertex] (d31) at (9,2.) {};
        \node[vertex] (d32) at (10,2.) {};
        \node[vertex] (d33) at (9.5,2.7) {};
        \draw[->] (d31) to (d32);
        \draw[->] (d32) to (d33);
        \draw[->] (d33) to (d31);
        \draw[->] [out=330,in=210] (d31) to (d32);
        \draw[->] [out=90,in=330] (d32) to (d33);
        \draw[->] [out=210,in=90] (d33) to (d31);

        \node[vertex] (d31a) at (9,0.3) {};
        \node[vertex] (d32a) at (10,0.3) {};
        \node[vertex] (d33a) at (9.5,1) {};
        \draw[-] (d31a) to (d32a);
        \draw[-] (d32a) to (d33a);
        \draw[-] (d33a) to (d31a);
        \node[] () at (9.5,0) {$\ord\frL_i=3$};
    \end{tikzpicture}
    \caption{Example diagrams of the structures of $G_0$ and $\Gasc_0$ for $\delta_\frl=1$ when $\frl\OO_0$ is principal.}
    \label{fig:Level0s_split}
\end{figure}
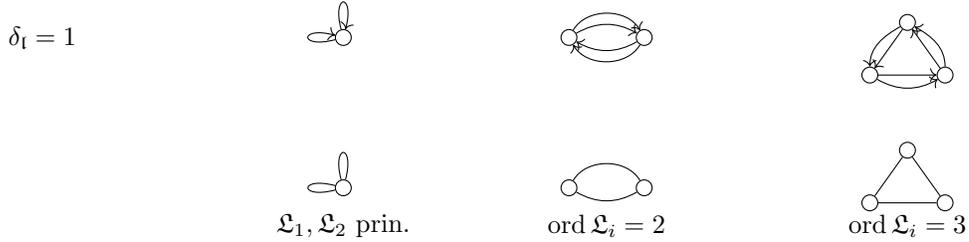

\begin{theorem}\label{thm:volcano_split}
    Assume that $\delta_\frl=1$. 
    Let $r$ be a positive integer.
    \begin{enumerate}[(i)]
        \item If $\dmin=0$ and $\frl\OO_0$ is principal then $G$ is an $r$-volcano.
        \item \label{thm:volcano_split:ii} If $\dmin>0$ then $G$ is an $r$-volcano if and only if the following three conditions hold:
        \begin{itemize}
            \item $\frl\OO_{\dmin-1}$ is principal;
            \item $r + 1 =\frac{\#(R/\frl)-1}{[\OO_0^\times:\OO_1^\times]} + 2$;
            \item $r + 1 = \frac{\#(R/\frl)}{[\Ocal_i^\times : \Ocal_{i+1}^\times]} + 1$ for $1 \leq i < \dmin$.
        \end{itemize}
    \end{enumerate}
\end{theorem}
\begin{proof}
	If $\dmin = 0$, then  $\Gasc$ and $\Gdesc$ are equal by Lemma~\ref{lem:Gasc_eq_Gdesc}. 
    The graph $\Gasc$ is a volcano if and only if it is connected, which is the case when $r_0 = 2$ by Lemma~\ref{lem:Gasc0conn_split}.
	
	If $\dmin >0$, then, again by Lemma~\ref{lem:Gasc_eq_Gdesc}, $\Gasc$ and $\Gdesc$ are equal if and only if $\frl\OO_{\dmin - 1}$ is principal, which implies that $\frl\OO_0$ is also principal. Further, by Lemma~\ref{lem:r_volc_minus_r0}, $\Gasc$ is an $r$-volcano for the specified $r$ if and only if $\Gasc_0$ is a connected $r_0$-regular graph, for some $r_0\leq 2$. By Lemma~\ref{lem:Gasc0conn_split}, since $\frl\OO_0$ is principal, $\Gasc_0$ is a connected 2-regular graph ($r_0=2$).
\end{proof}

Theorems~\ref{thm:volcano_inertramified} and \ref{thm:volcano_split} carefully record the impact of the unit groups $\OO_i^\times$ on the regularity of the graph. This impact was already evident in Theorem~\ref{thm:graph_structure}. The following example shows how changes in the unit group along the multiplicator ladder can balance out, resulting in an $r$-volcano structure.

\begin{example}\label{ex:different_units_ratios}
    Consider the isogeny class \href{https://www.lmfdb.org/Variety/Abelian/Fq/2/101/o\_dl}{2.101.o\_dl}, defined by the polynomial
    \[ h(x)=x^4 + 14 x^3 + 89 x^2 + 1414 x + 10201 \]
    consisting of absolutely simple ordinary abelian surfaces over $\FF_{101}$.
    The order $\OO_2=\ZZ[\pi,101/\pi]$ has a unique singular maximal ideal $\frl_2$, which is above the rational prime $3$.
    There exists a unique maximal suborder $R$ such that $\OO_2$ belongs to the $\frl$-multiplicator ladder of $R$, where $\frl = \frl_2 \cap R$:
    \[ R \subsetneq \OO_2 \subsetneq \OO_1 \subsetneq \OO_0. \]
    Note that $\dmin=2$
    and that $\OO_0$ is the maximal order of $K=\QQ[\pi]$.
    One computes that $\frl$ has residue field $\FF_9$, splits in $\OO_0$ giving $\delta_\frl=1$, and $\frl\OO_2$ is a principal $\OO_2$-ideal.
    Also, $[\OO_0^\times:\OO_1^\times] = 4$ and $[\OO_1^\times:\OO_2^\times] = 3$.
    The $54$ isomorphism classes of abelian varieties in the isogeny class distrubute into two isomorphic connected components $G$ of the $(R,\frl)$-isogeny graph.
    By Theorem~\ref{thm:volcano_split}, each $G$ is then a $3$-volcano.
    See Figure~\ref{fig:isog_graph_directed_non-const_units_rations}.

    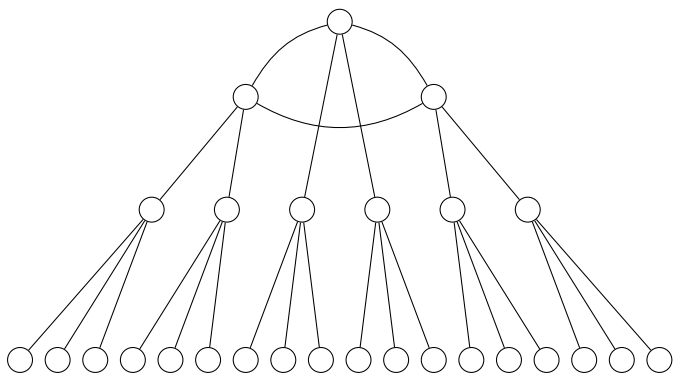
\begin{figure}[h!]
    \centering
    \begin{tikzpicture}
        \node[circle,draw] (01) at (-1.25,-.5) {};
    	\node[circle,draw] (02) at (0,0.5) {};
    	\node[circle,draw] (03) at (1.25,-.5) {};
        \node[circle,draw] (11) at (-2.5,-2) {};
    	\node[circle,draw] (12) at (-1.5,-2) {};
    	\node[circle,draw] (13) at (-0.5,-2) {};
        \node[circle,draw] (14) at (0.5,-2) {};
    	\node[circle,draw] (15) at (1.5,-2) {};
    	\node[circle,draw] (16) at (2.5,-2) {};
        \node[circle,draw] (21) at (-4.25,-4) {};
    	\node[circle,draw] (22) at (-3.75,-4) {};
    	\node[circle,draw] (23) at (-3.25,-4) {};
        \node[circle,draw] (24) at (-2.75,-4) {};
    	\node[circle,draw] (25) at (-2.25,-4) {};
    	\node[circle,draw] (26) at (-1.75,-4) {};
        \node[circle,draw] (27) at (-1.25,-4) {};
    	\node[circle,draw] (28) at (-0.75,-4) {};
    	\node[circle,draw] (29) at (-0.25,-4) {};
        \node[circle,draw] (210) at (0.25,-4) {};
    	\node[circle,draw] (211) at (0.75,-4) {};
    	\node[circle,draw] (212) at (1.25,-4) {};
        \node[circle,draw] (213) at (1.75,-4) {};
    	\node[circle,draw] (214) at (2.25,-4) {};
    	\node[circle,draw] (215) at (2.75,-4) {};
        \node[circle,draw] (216) at (3.25,-4) {};
    	\node[circle,draw] (217) at (3.75,-4) {};
        \node[circle,draw] (218) at (4.25,-4) {};
        \draw[-] (01) to [out=60, in=195](02);
        \draw[-] (02) to [out=345, in=120] (03);
        \draw[-] (03) to [out=210, in=330] (01);
        \draw[-] (11) to (01);
        \draw[-] (12) to (01);
        \draw[-] (13) to (02);
        \draw[-] (14) to (02);
        \draw[-] (15) to (03);
        \draw[-] (16) to (03);
        \draw[-] (21) to (11);
        \draw[-] (22) to (11);
        \draw[-] (23) to (11);
        \draw[-] (24) to (12);
        \draw[-] (25) to (12);
        \draw[-] (26) to (12);
        \draw[-] (27) to (13);
        \draw[-] (28) to (13);
        \draw[-] (29) to (13);
        \draw[-] (210) to (14);
        \draw[-] (211) to (14);
        \draw[-] (212) to (14);
        \draw[-] (213) to (15);
        \draw[-] (214) to (15);
        \draw[-] (215) to (15);
        \draw[-] (216) to (16);
        \draw[-] (217) to (16);
        \draw[-] (218) to (16);
    \end{tikzpicture}
    \caption{The ascending graph associated to a connected component of the $(R,\frl)$-isogeny graph in Example~\ref{ex:different_units_ratios}.}
    \label{fig:isog_graph_directed_non-const_units_rations}
    \end{figure}
\end{example}

\section{Comparisons with previous work and examples}\label{sec:examples}
We begin by proving that Theorem~\ref{thm:graph_structure} generalizes \cite[Theorem 4.3]{BrooksJetchevWesolowski17}, as claimed in Comparison~\ref{compBJW}. 
Recall the notation introduced for this comparison:
Let $A_0$ be an ordinary and absolutely simple abelian variety defined over a finite field $\Fq$.
The geometric and $\Fq$-endomorphism algebras of $A_0$ coincide, and they are both isomorphic to the number field $K=\QQ(\pi)$.
Let $K^+$ be the totally real subfield of $K$ and let $\frl^+$ be a maximal ideal of $\OO_{K^+}$ above a rational prime $\ell$ coprime with $q$. Set
\begin{equation}\label{eq:Randfrl}
    R \coloneqq \ZZ + N_1\OO_{K^+} + \ell N_2\OO_K
    \qquad\text{and}\qquad \frl\coloneqq \frl^+\OO_K\cap R,
\end{equation}
where $[\OO_{K^+} : \ZZ[\pi + q/\pi]]= N_1\cdot \ell^a$ for $\gcd(N_1, \ell) = 1$ and $[\OO_K : \ZZ[\pi, q/\pi]] = N_2$.

\begin{lemma}\label{lem:resfieldBJW}
    The ideal $\frl$ is a maximal ideal of $R$, with residue field isomorphic to $\OO_{K^+}/\frl^+$ and $R$ is Bass at $\frl$.
\end{lemma}
\begin{proof}
    By construction, the order $R$ contains $\OO_{K^+}$ locally at $\ell$. 
    Also, $\ell N_2 \OO_K \subseteq \ell\ZZ[\pi,q/\pi] \subseteq \frl^ +\OO_K$.
    Then $R/\frl$ has a natural structure of $\OO_{K^ +}/\frl^ +$-vector space and $(R+\frl^ +\OO_K)_\ell = (\OO_{K^ +}+\frl^ +\OO_K)_\ell$. In particular,
    \begin{equation}\label{eq:comp1}
        \frac{\OO_{K^+}}{\frl^+} \hookrightarrow\frac{R}{\frl} \cong \frac{R+\frl^+\OO_K}{\frl^+\OO_K} \cong \frac{\OO_{K^+}+\frl^+\OO_K}{\frl^+\OO_K}.
    \end{equation}
    
    Furthermore, 
    \begin{equation}\label{eq:comp2}
    \frac{\OO_{K^+}+\frl^+\OO_K}{\frl^+\OO_K} \hookrightarrow\frac{\OO_K}{\frl^+\OO_K}\cong \frac{\OO_{K^+}}{\frl^+}\oplus \frac{\OO_{K^+}}{\frl^+}.
    \end{equation}
    
    By \cite[Theorem 2.1]{BrooksJetchevWesolowski17}, the inclusion on the left of Equation~\eqref{eq:comp2} is not an isomorphism. It follows that the inclusion ${\OO_{K^+}}/{\frl^+} \hookrightarrow{R}/{\frl}$ in Equation~\eqref{eq:comp1} is an isomorphism.
    Hence, $\frl$ is a maximal ideal of $R$, with residue field isomorphic to $\OO_{K^+}/\frl^+$. 
    Moreover, the order $R$ is Bass at $\frl$ by \cite[Lemma 4.4]{BrooksJetchevWesolowski17}.
\end{proof}
    
The following two lemmas, Lemma~\ref{lem:compvert} and Lemma~\ref{lem:compedges}, give the statement in Comparison~\ref{compBJW}.
    
\begin{lemma}\label{lem:compvert}
    The set of vertices of the $(R,\frl)$-graph and of the $\frl^+$-graph truncated at $\Fq$ are the same and $R$ is strictly contained in the endomorphism ring of each $A\in\Ical_R$, that is, we have equalities:
    \[ \{ A\in \Ical : (\End(A))_\ell \supseteq (\OO_{K^+})_\ell \} = 
    \{ A\in \Ical : \End(A) \supsetneq R \}. \]
\end{lemma}
    \begin{proof}
    It suffices to show that the set $\Ical_R$ of overorders of $\ZZ[\pi,q/\pi]$ containing $R$ equals the set $\mathcal{Z}$ of overorders of $\ZZ[\pi,q/\pi]$ containing $\OO_{K^+}$ locally at $\ell$, and that $R$ is not an overorder of $\ZZ[\pi,q/\pi]$.
    If an overorder of $\ZZ[\pi,q/\pi]$ contains $R$, then it contains $\OO_{K^+}$ locally at $\ell$, so $\Ical_R\subseteq \mathcal{Z}$. 
    For the reverse inclusion, let $S\in\mathcal{Z}$.
    By localizing at $\ell$, we see that $S\supset N_1\OO_{K^+}$.
    Since $N_2\ell\OO_K \subseteq \ell\ZZ[\pi,q/\pi]$, we get
    \[ R \subseteq \ZZ + N_1\OO_{K^+} + \ell\ZZ[\pi,q/\pi] \subsetneq \ZZ[\pi,q/\pi] + N_1\OO_{K^+} \subseteq S.\]
\end{proof}
    
\begin{lemma}\label{lem:compedges}
    For every $A,B \in \Ical_R$, there is an equality of sets:
    \[\{ \frl^+ \text{-isogenies } A\to B \} = \{ \frl \text{-isogenies } A\to B \}.\]
\end{lemma}
\begin{proof}
    By definition the kernel of an $\frl$-isogeny is stable under the action of the order $R$ and it is killed by a power of $\ell$. 
    The order $R$ contains $\OO_{K^+}$ locally at $\ell$, so the kernel of an $\frl$-isogeny is also stable by the action of $\OO_{K^+,\ell}$. 
    Since $R/\frl$ is isomorphic to $\OO_{K^+}/\frl^+$ by Lemma~\ref{lem:resfieldBJW}, every $\frl$-isogeny is also an $\frl^+$-isogeny.
    Theorem~\ref{thm:graph_structure} and \cite[Prop.~4.10]{BrooksJetchevWesolowski17} show that
    we have the same number of $\frl^+$- and $\frl$-isogenies, after taking into account the contribution of automorphisms (cf.~Appendix~\ref{app:auts}).
\end{proof}

\begin{remark}\label{rmk:discrepancyBJW}
    There is a minor mistake in \cite[Theorem 4.3.(ii)]{BrooksJetchevWesolowski17}.
    Let $\OO_0$ be the endomorphism ring corresponding to the surface of a connected component of the $\frl^+$-isogeny graph.
    Assume that we are in the inert case, that is, $\frl^+\OO_0$ is a maximal ideal of $\OO_0$,  and that $\frl^+\OO_0$ is not a principal $\OO_0$-ideal.
    Then the surface consists of multiple vertices but no edges between them, as correctly predicted by \cite[Prop.~4.10.(i)]{BrooksJetchevWesolowski17}.
    Example~\ref{ex:c_ex_BJW} gives an isogeny class where this situation occurs.
    This explains the discrepancy between Theorem~\ref{thm:graph_structure}.\ref{it:surface} and \cite[Theorem 4.3.(ii)]{BrooksJetchevWesolowski17}.
\end{remark}

Because our work generalizes the results of \cite{BrooksJetchevWesolowski17}, it also generalizes the case of ordinary elliptic curves originally treated by \cite{Kohel94}.
Additionally, our work recovers the structure theorem for supersingular elliptic curves over a prime field $\FF_p$, as we explain in Example~\ref{ex:ssec}. 

\begin{example}[Supersingular elliptic curve isogeny volcanoes]\label{ex:ssec}
    We recover \cite[Theorem 2.7]{DelfsGalbraith} via Theorem~\ref{thm:graph_structure}.
    Let $E/\mathbb{F}_p$ be supersingular with endomorphism algebra $K = \mathbb{Q}(\sqrt{-p})$.
    The possible endomorphism rings of supersingular elliptic curves in the isogeny class of $E$ are $\mathbb{Z}[\frac{1 + \sqrt{-p}}{2}]$ and $\mathbb{Z}[\sqrt{-p}]$ if $p\equiv3\pmod{4}$, or only $\mathbb{Z}[\sqrt{-p}]$ if $p\not\equiv 3\pmod{4}$. There is a unique order of index $\ell$ inside any imaginary quadratic order.
    If $p\equiv 3\pmod{4}$ and $\ell = 2$, the connected components of 2-isogeny graphs have two levels. Otherwise, the $\ell$-isogeny graph connected components have one level. In the first case, the edges are degree-2 isogenies defined over $\FF_p$, determined by the ideals above $2\ZZ$. In the latter case, degree-$\ell$ isogenies are $\frl$-isogenies for primes $\frl$ above $\ell\ZZ$ when $\ell$ is split or ramified, and there are no edges otherwise.
    Theorem~\ref{thm:volcano_inertramified} and Theorem~\ref{thm:volcano_split} apply to recover the volcano structure of the graph, as orders in imaginary quadratic fields are all Bass and $\mathbb{Z}[\sqrt{-p}]$ is the minimal endomorphism ring. Furthermore, the unit groups of these orders are $\{\pm1\}$ by the congruence conditions on $p$.
    If $\ell = 2$ and $p\equiv 3\pmod{4}$, the $2$-isogeny graph connected components are volcanoes with two levels. Otherwise, the $\ell$-isogeny graph connected components are all one-level volcanoes. 
\end{example}

\appendix
\section{The Effect of Automorphisms}\label{app:auts}

Classical isogeny graphs take vertices to be isomorphisms classes of abelian varieties and directed edges to be isogenies, often of a given shape.
 In practice, an isogeny is defined as having a fixed abelian variety as the domain, not an isomorphism class. This leads to a choice for how to draw an edge in an isogeny graph:
\begin{enumerate}
	\item\label{aut_choice1} Fix an isomorphism class representative $A$ for a given vertex. Outgoing edges from that vertex are distinct isogenies with domain $A$.
	 When dealing with separable isogenies, the edges can be specified by the kernels of the different isogenies outgoing from $A$. In this case, we say that edges are up to \textbf{post-composition} with automorphisms.
	\item\label{aut_choice2} Fix an isomorphism class representative $A$ for a given vertex. Outgoing edges from that vertex are the distinct orbits of $\Aut(A)$ on the set of isogenies. In particular, two isogenies $\varphi,\psi$ from $A$ are in the same orbit if there exists $\eta\in\Aut(A)$ such that $\varphi = \psi\circ\eta$. In this case, we say that edges are up to \textbf{pre- and post-composition} with automorphisms. 
\end{enumerate} 

In this work, we make choice (\ref{aut_choice2}): we draw edges up to pre- and post-composition with automorphisms. Recall $\Aut(A) = \Aut_{\FF_q}(A)$.
\begin{remark}
When the automorphism group of $A$ is $\Aut(A)\cong \{[\pm 1]\}$, the choices (\ref{aut_choice1}) and (\ref{aut_choice2}) above are equivalent. This is due to the fact that kernels of isogenies are subgroups and as such are invariant under $[\pm 1]$. This is the case for any elliptic curve $E$ with $j(E)\not\in\{0,1728\}$. 
In particular, an elliptic curve $E/\mathbb{F}_q$ has more than two $\mathbb{F}_q$-automorphisms (that is, $|\Aut(E)| > 2$) if and only if either of the following holds:
\begin{itemize}
  \item $j(E) = 0$ and at least one of the following is true:
  \begin{itemize}
    \item $p = 2$ or $p = 3$;
    \item $p \equiv 1 \pmod{3}$;
    \item $q = p^k$ with $k > 1$;
  \end{itemize}
  \item $j(E) = 1728$ and at least one of the following is true:
  \begin{itemize}
    \item $p = 2$ or $p = 3$;
    \item $p \equiv 1 \pmod{4}$;
    \item $q = p^k$ with $k > 1$.
  \end{itemize}
\end{itemize}
\end{remark}

\begin{remark}
Let $A$ be a simple abelian variety over $\Fq$.
Then $\Aut(A)$ is finite if and only if $A$ has dimension one, that is, $A$ is an elliptic curve.
The product of two non-isogenous elliptic curves has $\Aut(E_1\times E_2) = \Aut(E_1)\times \Aut(E_2)$, which is finite. 
But if $E_1\sim E_2$, then $\End(E_1\times E_2)$ contains $\operatorname{Mat}_2(\ZZ)$, so the automorphism group $\Aut(E_1\times E_2)$ contains $\operatorname{GL}_2(\ZZ)$, and is therefore infinite. 
From these observations we deduce that the abelian varieties over $\Fq$ with finite automorphism groups are precisely the products of non-isogenous elliptic curves.
\end{remark}

\begin{example}\label{ex:auts}
Let $p = 19$, $E: y^2 = x^3 + 1/\mathbb{F}_p$. The (ordinary) elliptic curve $E$ has $j$-invariant $j(E) = 0$ and automorphism group of size $6$:
\[\Aut(E)= \{[\pm1],\pm\xi, \pm\xi^2: \xi(x,y)\coloneqq (7x,y), [-1](x,y)\coloneqq(x,-y)\}.\]
The full 2-torsion group of $E$ is defined over $\mathbb{F}_p$:
\[E[2] = \{O_E, (8,0),(12,0),(-1,0)\}.\]
Each nontrivial point of $E[2]$ defines a distinct $2$-isogeny, but $\Aut(E)$ permutes these points, so there is only one orbit in $E[2]$ under the action of $\Aut(E)$. This leads to two different 2-isogeny graphs, depending on the choice (\ref{aut_choice1}) or (\ref{aut_choice2}) above. See Figure~\ref{fig:main}.

\begin{figure}[h]
    \centering
    \begin{subfigure}{0.45\textwidth}
        \centering
        \begin{tikzpicture}
            \node[] (E) at (0,0) {$E:y^2 =x^3 + 1$};
            \node[] (E2) at (0,2) {$E_2:y^2 = x^3 +9x+3$};
            
            \draw[->] (E) to [out=60,in=300] (E2);
            \draw[->] (E) to [out = 30,in=330] (E2);
            \draw[->] (E) to [out=120,in=240](E2);
            \draw[->] (E2) to [out=210,in=150](E);
        \end{tikzpicture}
        \caption{Edges drawn up to post-composition with automorphisms (choice \ref{aut_choice1}).}
        \label{fig:sub1}
    \end{subfigure}
    \hfill
    \begin{subfigure}{0.45\textwidth}
        \centering
        \begin{tikzpicture}
        \node[] (E) at (0,0) {$E:y^2 = x^3 + 1$};
        \node[] (E2) at (0,2) {$E_2:y^2 = x^3 +9x+3$};
        
        \draw[->] (E) to [out=60,in=300] (E2);
        \draw[->] (E2) to [out=240,in=120] (E);
        \end{tikzpicture}
        \caption{Edges drawn up to pre- and post-composition with automorphisms (choice \ref{aut_choice2}).}
        \label{fig:sub2}
    \end{subfigure}
    \caption{The 2-isogeny graph component containing $E/\mathbb{F}_{19}$, as in Example~\ref{ex:auts}, drawn with the two choices of convention for drawing the edges.}
    \label{fig:main}
\end{figure}

\end{example}

We remark that both choices appear in the literature concerning isogeny graphs:
\begin{itemize}
	\item In \cite[p. 88]{Kohel94}, the author makes choice (\ref{aut_choice2}) and defines the edges of the elliptic curve $\ell$-isogeny graphs up to pre- and post-composition with automorphisms. 
	\item In \cite[p. 3]{BrooksJetchevWesolowski17}, the authors make choice (\ref{aut_choice1}), identifying isogenies by their kernels. They specify that there is an edge of multiplicity $m$ connecting abelian varieties $A$ and $B$ whenever there are $m$ distinct subgroups $\kappa$ of $A$ that are kernels of $\frl$-isogenies such that $A/\kappa\cong B$.
	\item In \cite[Remark 3.1.2]{Martindale18}, the author does not need to make a choice. Indeed the author is working with polarized objects, so only considering roots of unity, as opposed to the full (unpolarized) automorphsm group, and under the assumption that the only roots of unity in the ring of integers in the CM field are $\pm1$.
    Hence, the choices (\ref{aut_choice1}) and (\ref{aut_choice2}) are equivalent.
	\item In \cite[p. 391]{IonicaThome20}, the authors make choice (\ref{aut_choice2}), defining isogenies $\varphi,\psi$ to be equivalent if there exist automorphisms $\eta_1,\eta_2$ such that $\varphi = \eta_2\circ\psi\circ\eta_1$.
\end{itemize}

We make choice (\ref{aut_choice2}) to avoid dealing with inequivalent isogenies which might arise incident to vertices with infinite automorphism groups. The following example shows how even in the case of ordinary elliptic curves, extra automorphisms affect the regularity of the isogeny graph in such a way to prevent a volcano structure.

\begin{example}[Elliptic curves with extra automorphisms]\label{ex:ec_extra_aut}
Suppose $G$ is a connected component of an $\ell$-isogeny graph of ordinary elliptic curves or elliptic curves over $\FF_p$. If $G$ has more than one level, then $G$ is an $(\ell+1)$-volcano if and only if the automorphism groups of the elliptic curves are precisely $[\pm1]$.
To see this, suppose $A/\mathbb{F}_q$ is an elliptic curve over a field of characteristic $p$ with extra automorphisms. There are two possible endomorphism rings of elliptic curves with extra automorphisms, namely $\End(A)\cong \ZZ[\sqrt{-1}]$ or $\ZZ[\frac{1 + \sqrt{-3}}{2}]$ which necessarily are at level $0$ of any $\frl$-multiplicator ladder. In particular, $[\OO_0^\times:\OO_1^\times] = 2$ or $3$, and $[\OO_i^\times:\OO_{i+1}^\times] = 1$ for all $i\geq 1$. Let $R\subseteq\End(A)$ and consider the $(R,\frl)$-isogeny volcano component containing vertex $A$. The in-degree of the vertices at level $0$ are given by Theorem~\ref{thm:graph_structure} to be $\delta_\frl + 1 + \frac{\#(R/\frl) - \delta_\frl}{2}$. The in-degree of the vertices at levels $i\geq 1$ are likewise $(\#(R/\frl)+1)$. Since $\#(R/\frl)\geq 2$ and $\delta_\frl\in\{-1,0,1\}$, it is not possible for this to be a regular graph. 
\end{example}

\bibliographystyle{alpha}
\bibliography{biblio.bib}

\begin{thebibliography}{{LMF}25}

\bibitem[Bas63]{Bass63}
Hyman Bass.
\newblock On the ubiquity of {G}orenstein rings.
\newblock {\em Math. Z.}, 82:8--28, 1963.

\bibitem[BJW17]{BrooksJetchevWesolowski17}
Ernest~Hunter Brooks, Dimitar Jetchev, and Benjamin Wesolowski.
\newblock Isogeny graphs of ordinary abelian varieties.
\newblock {\em Res. Number Theory}, 3:Paper No. 28, 38, 2017.

\bibitem[BKM23]{BergKarMar}
Jonas Bergstr\"om, Valentijn Karemaker, and Stefano Marseglia.
\newblock Polarizations of abelian varieties over finite fields via canonical liftings.
\newblock {\em Int. Math. Res. Not. IMRN}, (4):3194--3248, 2023.

\bibitem[BKM24]{BergKarMar24_arXiv}
Jonas {Bergstr{\"o}m}, Valentijn {Karemaker}, and Stefano {Marseglia}.
\newblock {Abelian varieties over finite fields with commutative endomorphism algebra: theory and algorithms}.
\newblock {\em arXiv e-prints}, page arXiv:2409.08865, September 2024.

\bibitem[BL94]{BuchmannLenstra94}
J.~A. Buchmann and H.~W. Lenstra, Jr.
\newblock Approximating rings of integers in number fields.
\newblock {\em J. Th\'{e}or. Nombres Bordeaux}, 6(2):221--260, 1994.

\bibitem[CHL24]{ChoHongLee24_arXiv}
Sungmun {Cho}, Jungtaek {Hong}, and Yuchan {Lee}.
\newblock {Orbital integrals and Ideal class monoids for a Bass order}.
\newblock {\em arXiv e-prints}, page arXiv:2408.16199, August 2024.

\bibitem[CS15]{CentelegheStixI}
Tommaso~Giorgio Centeleghe and Jakob Stix.
\newblock Categories of abelian varieties over finite fields, {I}: {A}belian varieties over {$\Bbb{F}_p$}.
\newblock {\em Algebra Number Theory}, 9(1):225--265, 2015.

\bibitem[Del69]{Deligne}
Pierre Deligne.
\newblock Vari\'et\'es ab\'eliennes ordinaires sur un corps fini.
\newblock {\em Invent. Math.}, 8:238--243, 1969.

\bibitem[DG16]{DelfsGalbraith}
Christina Delfs and Steven~D. Galbraith.
\newblock Computing isogenies between supersingular elliptic curves over {$\Bbb{F}_p$}.
\newblock {\em Des. Codes Cryptogr.}, 78(2):425--440, 2016.

\bibitem[FM02]{FouquetMorain02}
Mireille Fouquet and Fran\c{c}ois Morain.
\newblock Isogeny volcanoes and the {SEA} algorithm.
\newblock In {\em Algorithmic number theory ({S}ydney, 2002)}, volume 2369 of {\em Lecture Notes in Comput. Sci.}, pages 276--291. Springer, Berlin, 2002.

\bibitem[FO70]{FerrandOlivier70}
D.~Ferrand and J.-P. Olivier.
\newblock Homomorphisms minimaux d'anneaux.
\newblock {\em J. Algebra}, 16:461--471, 1970.

\bibitem[How95]{Howe95}
Everett~W. Howe.
\newblock Principally polarized ordinary abelian varieties over finite fields.
\newblock {\em Trans. Amer. Math. Soc.}, 347(7):2361--2401, 1995.

\bibitem[HS20]{HofmannSircana20}
Tommy Hofmann and Carlo Sircana.
\newblock On the computation of overorders.
\newblock {\em Int. J. Number Theory}, 16(4):857--879, 2020.

\bibitem[IT20]{IonicaThome20}
Sorina Ionica and Emmanuel Thom\'{e}.
\newblock Isogeny graphs with maximal real multiplication.
\newblock {\em J. Number Theory}, 207:385--422, 2020.

\bibitem[JP20]{JP20}
Bruce~W. Jordan and Bjorn Poonen.
\newblock The analytic class number formula for 1-dimensional affine schemes.
\newblock {\em Bull. Lond. Math. Soc.}, 52(5):793--806, 2020.

\bibitem[Kan11]{Kani11}
Ernst Kani.
\newblock Products of {CM} elliptic curves.
\newblock {\em Collect. Math.}, 62(3):297--339, 2011.

\bibitem[KK24]{KirschmerKluners24}
Markus {Kirschmer} and J{\"u}rgen {Kl{\"u}ners}.
\newblock {Enumerating orders in number fields}.
\newblock {\em arXiv e-prints}, page arXiv:2411.08568, November 2024.

\bibitem[Koh96]{Kohel94}
David Kohel.
\newblock {\em Endomorphism rings of elliptic curves over finite fields}.
\newblock PhD thesis, University of California, Berkely, 1996.

\bibitem[{LMF}25]{lmfdb}
The {LMFDB Collaboration}.
\newblock The {L}-functions and modular forms database.
\newblock \url{https://www.lmfdb.org}, 2025.
\newblock [Online; accessed 11 July 2025].

\bibitem[Mar18]{Martindale18}
Chloe Martindale.
\newblock {\em Isogeny Graphs, Modular Polynomials, and Applications}.
\newblock Doctoral dissertation, Universiteit Leiden and Universit\'e de Bordeaux, 2018.

\bibitem[Mar21]{Mar_sqfree}
Stefano Marseglia.
\newblock Computing square-free polarized abelian varieties over finite fields.
\newblock {\em Math. Comp.}, 90(328):953--971, 2021.

\bibitem[Mar22]{Mar22_extensions}
Stefano Marseglia.
\newblock Computing base extensions of ordinary abelian varieties over finite fields.
\newblock {\em Int. J. Number Theory}, 18(9):1957--1974, 2022.

\bibitem[Mar24]{MarType}
Stefano Marseglia.
\newblock Cohen-{M}acaulay type of orders, generators and ideal classes.
\newblock {\em J. Algebra}, 658:247--276, 2024.

\bibitem[MS25]{MS25}
Stefano Marseglia and Caleb Springer.
\newblock Abelian varieties over finite fields and their groups of rational points.
\newblock {\em Algebra Number Theory}, 19(3):521--550, 2025.

\bibitem[Neu99]{neukirch}
J\"{u}rgen Neukirch.
\newblock {\em Algebraic number theory}, volume 322 of {\em Grundlehren der mathematischen Wissenschaften [Fundamental Principles of Mathematical Sciences]}.
\newblock Springer-Verlag, Berlin, 1999.
\newblock Translated from the 1992 German original and with a note by Norbert Schappacher, With a foreword by G. Harder.

\bibitem[OS20]{OswalShankar}
Abhishek Oswal and Ananth~N. Shankar.
\newblock Almost ordinary abelian varieties over finite fields.
\newblock {\em J. Lond. Math. Soc. (2)}, 101(3):923--937, 2020.

\bibitem[Piz80]{Pizer}
Arnold Pizer.
\newblock An algorithm for computing modular forms on {$\Gamma \sb{0}(N)$}.
\newblock {\em J. Algebra}, 64(2):340--390, 1980.

\bibitem[Sut13]{Sutherland_2013}
Andrew Sutherland.
\newblock Isogeny volcanoes.
\newblock {\em The Open Book Series}, 1(1):507–530, November 2013.

\bibitem[Tat66]{Tate66}
John Tate.
\newblock Endomorphisms of abelian varieties over finite fields.
\newblock {\em Invent. Math.}, 2:134--144, 1966.

\bibitem[Wat69]{Waterhouse}
William~C. Waterhouse.
\newblock Abelian varieties over finite fields.
\newblock {\em Ann. Sci. \'{E}cole Norm. Sup. (4)}, 2:521--560, 1969.

\end{thebibliography}

\end{document}